\documentstyle[12pt]{article}
\def\less{{\Bbb l}} 
\def\conj{{\Bbb j}} 
\def\more{{\cdot\!\! > }}
\def\Bbb{\bf}
\def\br{{\Bbb R}}
\def\bcp{{\Bbb  CP}}

\def\vol{~d\mu}

\def\eel#1{\label{#1}\end{equation}}

\newtheorem{main}{Theorem}

\def\<{\langle}

\def\bea{ \begin{eqnarray*} }
\def\eea{ \end{eqnarray*} }
\def\be{ \begin{equation} }
\def\ee{ \end{equation} }
\def\BC{{\Bbb C}}

\def\CO{{\cal O}}

\def\ztil{{\tilde Z}}

\def\d{{D_0}}
\def\i{{\cal I}}
\def\dbar{\bar{D}_0}
\def\z{{Z_0}}
\def\ddbar{{D_0\bar{D}_0}}

\def\llbar{{\ell\bar{\ell}}}
\def\ext{{\rm Ext}}
\def\Hom{{\rm Hom}}
\def\lext{{\cal E} xt}
\def\taurr{\tau_{Z_0}^0|_{D_0\bar{D}_0}}
\def\taur{\tau_{Z_0}^0|_{D_0\bar{D}_0}}

\newtheorem{thm}{Theorem}

\newtheorem{lem}{Lemma}
\newtheorem{propn}{Proposition}
\newtheorem{cor}{Corollary}

\newcounter{exam}[section]
\renewcommand{\theexam}{\thesection.\arabic{exam}}

\newenvironment{remark}{\medskip
\noindent {\bf Remark.}}{\hfill $\Box$}

\newenvironment{proof}{\medskip
\noindent {\bf Proof.}}{\hfill \rule{.5em}{1em}\mbox{}\bigskip}

\catcode`@=11 \catcode`!=11

\expandafter\ifx\csname fiverm\endcsname\relax
  \let\fiverm\fivrm
\fi

\let\!latexendpicture=\endpicture
\let\!latexframe=\frame
\let\!latexlinethickness=\linethickness
\let\!latexmultiput=\multiput
\let\!latexput=\put

\def\@picture(#1,#2)(#3,#4){%
  \@picht #2\unitlength
  \setbox\@picbox\hbox to #1\unitlength\bgroup
  \let\endpicture=\!latexendpicture
  \let\frame=\!latexframe
  \let\linethickness=\!latexlinethickness
  \let\multiput=\!latexmultiput
  \let\put=\!latexput
  \hskip -#3\unitlength \lower #4\unitlength \hbox\bgroup}

\catcode`@=12 \catcode`!=12

\catcode`!=11 

\def\PiC{P\kern-.12em\lower.5ex\hbox{I}\kern-.075emC}
\def\PiCTeX{\PiC\kern-.11em\TeX}

\def\!ifnextchar#1#2#3{%
  \let\!testchar=#1%
  \def\!first{#2}%
  \def\!second{#3}%
  \futurelet\!nextchar\!testnext}
\def\!testnext{%
  \ifx \!nextchar \!spacetoken
    \let\!next=\!skipspacetestagain
  \else
    \ifx \!nextchar \!testchar
      \let\!next=\!first
    \else
      \let\!next=\!second
    \fi
  \fi
  \!next}
\def\\{\!skipspacetestagain}
  \expandafter\def\\ {\futurelet\!nextchar\!testnext}
\def\\{\let\!spacetoken= } \\  

\def\!tfor#1:=#2\do#3{%
  \edef\!fortemp{#2}%
  \ifx\!fortemp\!empty
    \else
    \!tforloop#2\!nil\!nil\!!#1{#3}%
  \fi}
\def\!tforloop#1#2\!!#3#4{%
  \def#3{#1}%
  \ifx #3\!nnil
    \let\!nextwhile=\!fornoop
  \else
    #4\relax
    \let\!nextwhile=\!tforloop
  \fi
  \!nextwhile#2\!!#3{#4}}

\def\!etfor#1:=#2\do#3{%
  \def\!!tfor{\!tfor#1:=}%
  \edef\!!!tfor{#2}%
  \expandafter\!!tfor\!!!tfor\do{#3}}

\def\!cfor#1:=#2\do#3{%
  \edef\!fortemp{#2}%
  \ifx\!fortemp\!empty
  \else
    \!cforloop#2,\!nil,\!nil\!!#1{#3}%
  \fi}
\def\!cforloop#1,#2\!!#3#4{%
  \def#3{#1}%
  \ifx #3\!nnil
    \let\!nextwhile=\!fornoop
  \else
    #4\relax
    \let\!nextwhile=\!cforloop
  \fi
  \!nextwhile#2\!!#3{#4}}

\def\!ecfor#1:=#2\do#3{%
  \def\!!cfor{\!cfor#1:=}%
  \edef\!!!cfor{#2}%
  \expandafter\!!cfor\!!!cfor\do{#3}}

\def\!empty{}
\def\!nnil{\!nil}
\def\!fornoop#1\!!#2#3{}

\def\!ifempty#1#2#3{%
  \edef\!emptyarg{#1}%
  \ifx\!emptyarg\!empty
    #2%
  \else
    #3%
  \fi}

\def\!getnext#1\from#2{%
  \expandafter\!gnext#2\!#1#2}%
\def\!gnext\\#1#2\!#3#4{%
  \def#3{#1}%
  \def#4{#2\\{#1}}%
  \ignorespaces}

\def\!getnextvalueof#1\from#2{%
  \expandafter\!gnextv#2\!#1#2}%
\def\!gnextv\\#1#2\!#3#4{%
  #3=#1%
  \def#4{#2\\{#1}}%
  \ignorespaces}

\def\!copylist#1\to#2{%
  \expandafter\!!copylist#1\!#2}
\def\!!copylist#1\!#2{%
  \def#2{#1}\ignorespaces}

\def\!wlet#1=#2{%
  \let#1=#2
  \wlog{\string#1=\string#2}}

\def\!listaddon#1#2{%
  \expandafter\!!listaddon#2\!{#1}#2}
\def\!!listaddon#1\!#2#3{%
  \def#3{#1\\#2}}

\def\!rightappend#1\withCS#2\to#3{\expandafter\!!rightappend#3\!#2{#1}#3}
\def\!!rightappend#1\!#2#3#4{\def#4{#1#2{#3}}}

\def\!leftappend#1\withCS#2\to#3{\expandafter\!!leftappend#3\!#2{#1}#3}
\def\!!leftappend#1\!#2#3#4{\def#4{#2{#3}#1}}

\def\!lop#1\to#2{\expandafter\!!lop#1\!#1#2}
\def\!!lop\\#1#2\!#3#4{\def#4{#1}\def#3{#2}}

\def\!loop#1\repeat{\def\!body{#1}\!iterate}
\def\!iterate{\!body\let\!next=\!iterate\else\let\!next=\relax\fi\!next}

\def\!!loop#1\repeat{\def\!!body{#1}\!!iterate}
\def\!!iterate{\!!body\let\!!next=\!!iterate\else\let\!!next=\relax\fi\!!next}

\def\!removept#1#2{\edef#2{\expandafter\!!removePT\the#1}}
{\catcode`p=12 \catcode`t=12 \gdef\!!removePT#1pt{#1}}

\def\placevalueinpts of <#1> in #2 {%
  \!removept{#1}{#2}}

\def\!mlap#1{\hbox to 0pt{\hss#1\hss}}
\def\!vmlap#1{\vbox to 0pt{\vss#1\vss}}

\def\!not#1{%
  #1\relax
    \!switchfalse
  \else
    \!switchtrue
  \fi
  \if!switch
  \ignorespaces}

\let\!!!wlog=\wlog              
\def\wlog#1{}

\newdimen\headingtoplotskip     
\newdimen\linethickness         
\newdimen\longticklength        
\newdimen\plotsymbolspacing     
\newdimen\shortticklength       
\newdimen\stackleading          
\newdimen\tickstovaluesleading  
\newdimen\totalarclength        
\newdimen\valuestolabelleading  

\newbox\!boxA                   
\newbox\!boxB                   
\newbox\!picbox                 
\newbox\!plotsymbol             
\newbox\!putobject              
\newbox\!shadesymbol            

\newcount\!countA               
\newcount\!countB               
\newcount\!countC               
\newcount\!countD               
\newcount\!countE               
\newcount\!countF               
\newcount\!countG               
\newcount\!fiftypt              
\newcount\!intervalno           
\newcount\!npoints              
\newcount\!nsegments            
\newcount\!ntemp                
\newcount\!parity               
\newcount\!scalefactor          
\newcount\!tfs                  
\newcount\!tickcase             

\newdimen\!Xleft                
\newdimen\!Xright               
\newdimen\!Xsave                
\newdimen\!Ybot                 
\newdimen\!Ysave                
\newdimen\!Ytop                 
\newdimen\!angle                
\newdimen\!arclength            
\newdimen\!areabloc             
\newdimen\!arealloc             
\newdimen\!arearloc             
\newdimen\!areatloc             
\newdimen\!bshrinkage           
\newdimen\!checkbot             
\newdimen\!checkleft            
\newdimen\!checkright           
\newdimen\!checktop             
\newdimen\!dimenA               
\newdimen\!dimenB               
\newdimen\!dimenC               
\newdimen\!dimenD               
\newdimen\!dimenE               
\newdimen\!dimenF               
\newdimen\!dimenG               
\newdimen\!dimenH               
\newdimen\!dimenI               
\newdimen\!distacross           
\newdimen\!downlength           
\newdimen\!dp                   
\newdimen\!dshade               
\newdimen\!dxpos                
\newdimen\!dxprime              
\newdimen\!dypos                
\newdimen\!dyprime              
\newdimen\!ht                   
\newdimen\!leaderlength         
\newdimen\!lshrinkage           
\newdimen\!midarclength         
\newdimen\!offset               
\newdimen\!plotheadingoffset    
\newdimen\!plotsymbolxshift     
\newdimen\!plotsymbolyshift     
\newdimen\!plotxorigin          
\newdimen\!plotyorigin          
\newdimen\!rootten              
\newdimen\!rshrinkage           
\newdimen\!shadesymbolxshift    
\newdimen\!shadesymbolyshift    
\newdimen\!tenAa                
\newdimen\!tenAc                
\newdimen\!tenAe                
\newdimen\!tshrinkage           
\newdimen\!uplength             
\newdimen\!wd                   
\newdimen\!wmax                 
\newdimen\!wmin                 
\newdimen\!xB                   
\newdimen\!xC                   
\newdimen\!xE                   
\newdimen\!xM                   
\newdimen\!xS                   
\newdimen\!xaxislength          
\newdimen\!xdiff                
\newdimen\!xleft                
\newdimen\!xloc                 
\newdimen\!xorigin              
\newdimen\!xpivot               
\newdimen\!xpos                 
\newdimen\!xprime               
\newdimen\!xright               
\newdimen\!xshade               
\newdimen\!xshift               
\newdimen\!xtemp                
\newdimen\!xunit                
\newdimen\!xxE                  
\newdimen\!xxM                  
\newdimen\!xxS                  
\newdimen\!xxloc                
\newdimen\!yB                   
\newdimen\!yC                   
\newdimen\!yE                   
\newdimen\!yM                   
\newdimen\!yS                   
\newdimen\!yaxislength          
\newdimen\!ybot                 
\newdimen\!ydiff                
\newdimen\!yloc                 
\newdimen\!yorigin              
\newdimen\!ypivot               
\newdimen\!ypos                 
\newdimen\!yprime               
\newdimen\!yshade               
\newdimen\!yshift               
\newdimen\!ytemp                
\newdimen\!ytop                 
\newdimen\!yunit                
\newdimen\!yyE                  
\newdimen\!yyM                  
\newdimen\!yyS                  
\newdimen\!yyloc                
\newdimen\!zpt                  

\newif\if!axisvisible           
\newif\if!gridlinestoo          
\newif\if!keepPO                
\newif\if!placeaxislabel        
\newif\if!switch                
\newif\if!xswitch               

\newtoks\!axisLaBeL             
\newtoks\!keywordtoks           

\newwrite\!replotfile           

\newhelp\!keywordhelp{The keyword mentioned in the error message in unknown.
Replace NEW KEYWORD in the indicated response by the keyword that
should have been specified.}    

\!wlet\!!origin=\!xM                   
\!wlet\!!unit=\!uplength               
\!wlet\!Lresiduallength=\!dimenG       
\!wlet\!Rresiduallength=\!dimenF       
\!wlet\!axisLength=\!distacross        
\!wlet\!axisend=\!ydiff                
\!wlet\!axisstart=\!xdiff              
\!wlet\!axisxlevel=\!arclength         
\!wlet\!axisylevel=\!downlength        
\!wlet\!beta=\!dimenE                  
\!wlet\!gamma=\!dimenF                 
\!wlet\!shadexorigin=\!plotxorigin     
\!wlet\!shadeyorigin=\!plotyorigin     
\!wlet\!ticklength=\!xS                
\!wlet\!ticklocation=\!xE              
\!wlet\!ticklocationincr=\!yE          
\!wlet\!tickwidth=\!yS                 
\!wlet\!totalleaderlength=\!dimenE     
\!wlet\!xone=\!xprime                  
\!wlet\!xtwo=\!dxprime                 
\!wlet\!ySsave=\!yM                    
\!wlet\!ybB=\!yB                       
\!wlet\!ybC=\!yC                       
\!wlet\!ybE=\!yE                       
\!wlet\!ybM=\!yM                       
\!wlet\!ybS=\!yS                       
\!wlet\!ybpos=\!yyloc                  
\!wlet\!yone=\!yprime                  
\!wlet\!ytB=\!xB                       
\!wlet\!ytC=\!xC                       
\!wlet\!ytE=\!downlength               
\!wlet\!ytM=\!arclength                
\!wlet\!ytS=\!distacross               
\!wlet\!ytpos=\!xxloc                  
\!wlet\!ytwo=\!dyprime                 

\!zpt=0pt                              
\!xunit=1pt
\!yunit=1pt
\!arearloc=\!xunit
\!areatloc=\!yunit
\!dshade=5pt
\!leaderlength=24in
\!tfs=256                              
\!wmax=5.3pt                           
\!wmin=2.7pt                           
\!xaxislength=\!xunit
\!xpivot=\!zpt
\!yaxislength=\!yunit
\!ypivot=\!zpt
  \!dimenA=50pt \!fiftypt=\!dimenA     

\!rootten=3.162278pt                   
\!tenAa=8.690286pt                     
\!tenAc=2.773839pt                     
\!tenAe=2.543275pt                     

\def\!cosrotationangle{1}      
\def\!sinrotationangle{0}      
\def\!xpivotcoord{0}           
\def\!xref{0}                  
\def\!xshadesave{0}            
\def\!ypivotcoord{0}           
\def\!yref{0}                  
\def\!yshadesave{0}            
\def\!zero{0}                  

\let\wlog=\!!!wlog
\def\normalgraphs{%
  \longticklength=.4\baselineskip
  \shortticklength=.25\baselineskip
  \tickstovaluesleading=.25\baselineskip
  \valuestolabelleading=.8\baselineskip
  \linethickness=.4pt
  \stackleading=.17\baselineskip
  \headingtoplotskip=1.5\baselineskip
  \visibleaxes
  \ticksout
  \nogridlines
  \unloggedticks}
\def\setplotarea x from #1 to #2, y from #3 to #4 {%
  \!arealloc=\!M{#1}\!xunit \advance \!arealloc -\!xorigin
  \!areabloc=\!M{#3}\!yunit \advance \!areabloc -\!yorigin
  \!arearloc=\!M{#2}\!xunit \advance \!arearloc -\!xorigin
  \!areatloc=\!M{#4}\!yunit \advance \!areatloc -\!yorigin
  \!initinboundscheck
  \!xaxislength=\!arearloc  \advance\!xaxislength -\!arealloc
  \!yaxislength=\!areatloc  \advance\!yaxislength -\!areabloc
  \!plotheadingoffset=\!zpt
  \!dimenput {{\setbox0=\hbox{}\wd0=\!xaxislength\ht0=\!yaxislength\box0}}
     [bl] (\!arealloc,\!areabloc)}
\def\visibleaxes{%
  \def\!axisvisibility{\!axisvisibletrue}}

\def\!fixkeyword#1{%
  \errhelp=\!keywordhelp
  \errmessage{Unrecognized keyword `#1': \the\!keywordtoks{NEW KEYWORD}'}}

\!keywordtoks={enter `i\fixkeyword}

\def\fixkeyword#1{%
  \!nextkeyword#1 }

\def\axis {%
  \def\!nextkeyword##1 {%
    \expandafter\ifx\csname !axis##1\endcsname \relax
      \def\!next{\!fixkeyword{##1}}%
    \else
      \def\!next{\csname !axis##1\endcsname}%
    \fi
    \!next}%
  \!offset=\!zpt
  \!axisvisibility
  \!placeaxislabelfalse
  \!nextkeyword}

\def\!axisbottom{%
  \!axisylevel=\!areabloc
  \def\!tickxsign{0}%
  \def\!tickysign{-}%
  \def\!axissetup{\!axisxsetup}%
  \def\!axislabeltbrl{t}%
  \!nextkeyword}

\def\!axistop{%
  \!axisylevel=\!areatloc
  \def\!tickxsign{0}%
  \def\!tickysign{+}%
  \def\!axissetup{\!axisxsetup}%
  \def\!axislabeltbrl{b}%
  \!nextkeyword}

\def\!axisleft{%
  \!axisxlevel=\!arealloc
  \def\!tickxsign{-}%
  \def\!tickysign{0}%
  \def\!axissetup{\!axisysetup}%
  \def\!axislabeltbrl{r}%
  \!nextkeyword}

\def\!axisright{%
  \!axisxlevel=\!arearloc
  \def\!tickxsign{+}%
  \def\!tickysign{0}%
  \def\!axissetup{\!axisysetup}%
  \def\!axislabeltbrl{l}%
  \!nextkeyword}

\def\!axisshiftedto#1=#2 {%
  \if 0\!tickxsign
    \!axisylevel=\!M{#2}\!yunit
    \advance\!axisylevel -\!yorigin
  \else
    \!axisxlevel=\!M{#2}\!xunit
    \advance\!axisxlevel -\!xorigin
  \fi
  \!nextkeyword}

\def\!axisvisible{%
  \!axisvisibletrue
  \!nextkeyword}

\def\!axisinvisible{%
  \!axisvisiblefalse
  \!nextkeyword}

\def\!axislabel#1 {%
  \!axisLaBeL={#1}%
  \!placeaxislabeltrue
  \!nextkeyword}

\expandafter\def\csname !axis/\endcsname{%
  \!axissetup 
  \if!placeaxislabel
    \!placeaxislabel
  \fi
  \if +\!tickysign 
    \!dimenA=\!axisylevel
    \advance\!dimenA \!offset 
    \advance\!dimenA -\!areatloc 
    \ifdim \!dimenA>\!plotheadingoffset
      \!plotheadingoffset=\!dimenA 
    \fi
  \fi}

\def\grid #1 #2 {%
  \!countA=#1\advance\!countA 1
  \axis bottom invisible ticks length <\!zpt> andacross quantity {\!countA} /
  \!countA=#2\advance\!countA 1
  \axis left   invisible ticks length <\!zpt> andacross quantity {\!countA} / }

\def\plotheading#1 {%
  \advance\!plotheadingoffset \headingtoplotskip
  \!dimenput {#1} [B] <.5\!xaxislength,\!plotheadingoffset>
    (\!arealloc,\!areatloc)}

\def\!axisxsetup{%
  \!axisxlevel=\!arealloc
  \!axisstart=\!arealloc
  \!axisend=\!arearloc
  \!axisLength=\!xaxislength
  \!!origin=\!xorigin
  \!!unit=\!xunit
  \!xswitchtrue
  \if!axisvisible
    \!makeaxis
  \fi}

\def\!axisysetup{%
  \!axisylevel=\!areabloc
  \!axisstart=\!areabloc
  \!axisend=\!areatloc
  \!axisLength=\!yaxislength
  \!!origin=\!yorigin
  \!!unit=\!yunit
  \!xswitchfalse
  \if!axisvisible
    \!makeaxis
  \fi}

\def\!makeaxis{%
  \setbox\!boxA=\hbox{
    \beginpicture
      \!setdimenmode
      \setcoordinatesystem point at {\!zpt} {\!zpt}
      \putrule from {\!zpt} {\!zpt} to
        {\!tickysign\!tickysign\!axisLength}
        {\!tickxsign\!tickxsign\!axisLength}
    \endpicturesave <\!Xsave,\!Ysave>}%
    \wd\!boxA=\!zpt
    \!placetick\!axisstart}

\def\!placeaxislabel{%
  \advance\!offset \valuestolabelleading
  \if!xswitch
    \!dimenput {\the\!axisLaBeL} [\!axislabeltbrl]
      <.5\!axisLength,\!tickysign\!offset> (\!axisxlevel,\!axisylevel)
    \advance\!offset \!dp  
    \advance\!offset \!ht  
  \else
    \!dimenput {\the\!axisLaBeL} [\!axislabeltbrl]
      <\!tickxsign\!offset,.5\!axisLength> (\!axisxlevel,\!axisylevel)
  \fi
  \!axisLaBeL={}}

\def\arrow <#1> [#2,#3]{%
  \!ifnextchar<{\!arrow{#1}{#2}{#3}}{\!arrow{#1}{#2}{#3}<\!zpt,\!zpt> }}

\def\!arrow#1#2#3<#4,#5> from #6 #7 to #8 #9 {%
  \!xloc=\!M{#8}\!xunit
  \!yloc=\!M{#9}\!yunit
  \!dxpos=\!xloc  \!dimenA=\!M{#6}\!xunit  \advance \!dxpos -\!dimenA
  \!dypos=\!yloc  \!dimenA=\!M{#7}\!yunit  \advance \!dypos -\!dimenA
  \let\!MAH=\!M
  \!setdimenmode
  \!xshift=#4\relax  \!yshift=#5\relax
  \!reverserotateonly\!xshift\!yshift
  \advance\!xshift\!xloc  \advance\!yshift\!yloc
  \!xS=-\!dxpos  \advance\!xS\!xshift
  \!yS=-\!dypos  \advance\!yS\!yshift
  \!start (\!xS,\!yS)
  \!ljoin (\!xshift,\!yshift)
  \!Pythag\!dxpos\!dypos\!arclength
  \!divide\!dxpos\!arclength\!dxpos
  \!dxpos=32\!dxpos  \!removept\!dxpos\!!cos
  \!divide\!dypos\!arclength\!dypos
  \!dypos=32\!dypos  \!removept\!dypos\!!sin
  \!halfhead{#1}{#2}{#3}
  \!halfhead{#1}{-#2}{-#3}
  \let\!M=\!MAH
  \ignorespaces}
  \def\!halfhead#1#2#3{%
    \!dimenC=-#1%
    \divide \!dimenC 2 
    \!dimenD=#2\!dimenC
    \!rotate(\!dimenC,\!dimenD)by(\!!cos,\!!sin)to(\!xM,\!yM)
    \!dimenC=-#1
    \!dimenD=#3\!dimenC
    \!dimenD=.5\!dimenD
    \!rotate(\!dimenC,\!dimenD)by(\!!cos,\!!sin)to(\!xE,\!yE)
    \!start (\!xshift,\!yshift)
    \advance\!xM\!xshift  \advance\!yM\!yshift
    \advance\!xE\!xshift  \advance\!yE\!yshift
    \!qjoin (\!xM,\!yM) (\!xE,\!yE)
    \ignorespaces}

\def\betweenarrows #1#2 from #3 #4 to #5 #6 {%
  \!xloc=\!M{#3}\!xunit  \!xxloc=\!M{#5}\!xunit%
  \!yloc=\!M{#4}\!yunit  \!yyloc=\!M{#6}\!yunit%
  \!dxpos=\!xxloc  \advance\!dxpos by -\!xloc
  \!dypos=\!yyloc  \advance\!dypos by -\!yloc
  \advance\!xloc .5\!dxpos
  \advance\!yloc .5\!dypos
  \let\!MBA=\!M
  \!setdimenmode
  \ifdim\!dypos=\!zpt
    \ifdim\!dxpos<\!zpt \!dxpos=-\!dxpos \fi
    \put {\!lrarrows{\!dxpos}{#1}}#2{} at {\!xloc} {\!yloc}
  \else
    \ifdim\!dxpos=\!zpt
      \ifdim\!dypos<\!zpt \!dypos=-\!zpt \fi
      \put {\!udarrows{\!dypos}{#1}}#2{} at {\!xloc} {\!yloc}
    \fi
  \fi
  \let\!M=\!MBA
  \ignorespaces}

\def\!lrarrows#1#2{
  {\setbox\!boxA=\hbox{$\mkern-2mu\mathord-\mkern-2mu$}%
   \setbox\!boxB=\hbox{$\leftarrow$}\!dimenE=\ht\!boxB
   \setbox\!boxB=\hbox{}\ht\!boxB=2\!dimenE
   \hbox to #1{$\mathord\leftarrow\mkern-6mu
     \cleaders\copy\!boxA\hfil
     \mkern-6mu\mathord-$%
     \kern.4em $\vcenter{\box\!boxB}$$\vcenter{\hbox{#2}}$\kern.4em
     $\mathord-\mkern-6mu
     \cleaders\copy\!boxA\hfil
     \mkern-6mu\mathord\rightarrow$}}}

\def\!udarrows#1#2{
  {\setbox\!boxB=\hbox{#2}%
   \setbox\!boxA=\hbox to \wd\!boxB{\hss$\vert$\hss}%
   \!dimenE=\ht\!boxA \advance\!dimenE \dp\!boxA \divide\!dimenE 2
   \vbox to #1{\offinterlineskip
      \vskip .05556\!dimenE
      \hbox to \wd\!boxB{\hss$\mkern.4mu\uparrow$\hss}\vskip-\!dimenE
      \cleaders\copy\!boxA\vfil
      \vskip-\!dimenE\copy\!boxA
      \vskip\!dimenE\copy\!boxB\vskip.4em
      \copy\!boxA\vskip-\!dimenE
      \cleaders\copy\!boxA\vfil
      \vskip-\!dimenE \hbox to \wd\!boxB{\hss$\mkern.4mu\downarrow$\hss}
      \vskip .05556\!dimenE}}}

\def\putbar#1breadth <#2> from #3 #4 to #5 #6 {%
  \!xloc=\!M{#3}\!xunit  \!xxloc=\!M{#5}\!xunit%
  \!yloc=\!M{#4}\!yunit  \!yyloc=\!M{#6}\!yunit%
  \!dypos=\!yyloc  \advance\!dypos by -\!yloc
  \!dimenI=#2
  \ifdim \!dimenI=\!zpt 
    \putrule#1from {#3} {#4} to {#5} {#6} 
  \else 
    \let\!MBar=\!M
    \!setdimenmode 
    \divide\!dimenI 2
    \ifdim \!dypos=\!zpt
      \advance \!yloc -\!dimenI 
      \advance \!yyloc \!dimenI
    \else
      \advance \!xloc -\!dimenI 
      \advance \!xxloc \!dimenI
    \fi
    \putrectangle#1corners at {\!xloc} {\!yloc} and {\!xxloc} {\!yyloc}
    \let\!M=\!MBar 
  \fi
  \ignorespaces}

\def\setbars#1breadth <#2> baseline at #3 = #4 {%
  \edef\!barshift{#1}%
  \edef\!barbreadth{#2}%
  \edef\!barorientation{#3}%
  \edef\!barbaseline{#4}%
  \def\!bardobaselabel{\!bardoendlabel}%
  \def\!bardoendlabel{\!barfinish}%
  \let\!drawcurve=\!barcurve
  \!setbars}
\def\!setbars{%
  \futurelet\!nextchar\!!setbars}
\def\!!setbars{%
  \if b\!nextchar
    \def\!!!setbars{\!setbarsbget}%
  \else
    \if e\!nextchar
      \def\!!!setbars{\!setbarseget}%
    \else
      \def\!!!setbars{\relax}%
    \fi
  \fi
  \!!!setbars}
\def\!setbarsbget baselabels (#1) {%
  \def\!barbaselabelorientation{#1}%
  \def\!bardobaselabel{\!!bardobaselabel}%
  \!setbars}
\def\!setbarseget endlabels (#1) {%
  \edef\!barendlabelorientation{#1}%
  \def\!bardoendlabel{\!!bardoendlabel}%
  \!setbars}

\def\!barcurve #1 #2 {%
  \if y\!barorientation
    \def\!basexarg{#1}%
    \def\!baseyarg{\!barbaseline}%
  \else
    \def\!basexarg{\!barbaseline}%
    \def\!baseyarg{#2}%
  \fi
  \expandafter\putbar\!barshift breadth <\!barbreadth> from {\!basexarg}
    {\!baseyarg} to {#1} {#2}
  \def\!endxarg{#1}%
  \def\!endyarg{#2}%
  \!bardobaselabel}

\def\!!bardobaselabel "#1" {%
  \put {#1}\!barbaselabelorientation{} at {\!basexarg} {\!baseyarg}
  \!bardoendlabel}

\def\!!bardoendlabel "#1" {%
  \put {#1}\!barendlabelorientation{} at {\!endxarg} {\!endyarg}
  \!barfinish}

\def\!barfinish{%
  \!ifnextchar/{\!finish}{\!barcurve}}

\def\putrectangle{%
  \!ifnextchar<{\!putrectangle}{\!putrectangle<\!zpt,\!zpt> }}
\def\!putrectangle<#1,#2> corners at #3 #4 and #5 #6 {%
  \!xone=\!M{#3}\!xunit  \!xtwo=\!M{#5}\!xunit%
  \!yone=\!M{#4}\!yunit  \!ytwo=\!M{#6}\!yunit%
  \ifdim \!xtwo<\!xone
    \!dimenI=\!xone  \!xone=\!xtwo  \!xtwo=\!dimenI
  \fi
  \ifdim \!ytwo<\!yone
    \!dimenI=\!yone  \!yone=\!ytwo  \!ytwo=\!dimenI
  \fi
  \!dimenI=#1\relax  \advance\!xone\!dimenI  \advance\!xtwo\!dimenI
  \!dimenI=#2\relax  \advance\!yone\!dimenI  \advance\!ytwo\!dimenI
  \let\!MRect=\!M
  \!setdimenmode
  \!shaderectangle
  \!dimenI=.5\linethickness
  \advance \!xone  -\!dimenI
  \advance \!xtwo   \!dimenI
  \putrule from {\!xone} {\!yone} to {\!xtwo} {\!yone}
  \putrule from {\!xone} {\!ytwo} to {\!xtwo} {\!ytwo}
  \advance \!xone   \!dimenI
  \advance \!xtwo  -\!dimenI%
  \advance \!yone  -\!dimenI
  \advance \!ytwo   \!dimenI
  \putrule from {\!xone} {\!yone} to {\!xone} {\!ytwo}
  \putrule from {\!xtwo} {\!yone} to {\!xtwo} {\!ytwo}
  \let\!M=\!MRect
  \ignorespaces}

\def\shaderectanglesoff{%
  \def\!shaderectangle{}%
  \ignorespaces}

\shaderectanglesoff

\def\!!shaderectangle{%
  \!dimenA=\!xtwo  \advance \!dimenA -\!xone
  \!dimenB=\!ytwo  \advance \!dimenB -\!yone
  \ifdim \!dimenA<\!dimenB
    \!startvshade (\!xone,\!yone,\!ytwo)
    \!lshade      (\!xtwo,\!yone,\!ytwo)
  \else
    \!starthshade (\!yone,\!xone,\!xtwo)
    \!lshade      (\!ytwo,\!xone,\!xtwo)
  \fi
  \ignorespaces}

\def\frame{%
  \!ifnextchar<{\!frame}{\!frame<\!zpt> }}
\long\def\!frame<#1> #2{%
  \beginpicture
    \setcoordinatesystem units <1pt,1pt> point at 0 0
    \put {#2} [Bl] at 0 0
    \!dimenA=#1\relax
    \!dimenB=\!wd \advance \!dimenB \!dimenA
    \!dimenC=\!ht \advance \!dimenC \!dimenA
    \!dimenD=\!dp \advance \!dimenD \!dimenA
    \let\!MFr=\!M
    \!setdimenmode
    \putrectangle corners at {-\!dimenA} {-\!dimenD} and {\!dimenB} {\!dimenC}
    \!setcoordmode
    \let\!M=\!MFr
  \endpicture
  \ignorespaces}

\def\rectangle <#1> <#2> {%
  \setbox0=\hbox{}\wd0=#1\ht0=#2\frame {\box0}}

\def\plot{%
  \!ifnextchar"{\!plotfromfile}{\!drawcurve}}
\def\!plotfromfile"#1"{%
  \expandafter\!drawcurve \input #1 /}

\def\setquadratic{%
  \let\!drawcurve=\!qcurve
  \let\!!Shade=\!!qShade
  \let\!!!Shade=\!!!qShade}

\def\setlinear{%
  \let\!drawcurve=\!lcurve
  \let\!!Shade=\!!lShade
  \let\!!!Shade=\!!!lShade}

\def\sethistograms{%
  \let\!drawcurve=\!hcurve}

\def\!qcurve #1 #2 {%
  \!start (#1,#2)
  \!Qjoin}
\def\!Qjoin#1 #2 #3 #4 {%
  \!qjoin (#1,#2) (#3,#4)             
  \!ifnextchar/{\!finish}{\!Qjoin}}

\def\!lcurve #1 #2 {%
  \!start (#1,#2)
  \!Ljoin}
\def\!Ljoin#1 #2 {%
  \!ljoin (#1,#2)                    
  \!ifnextchar/{\!finish}{\!Ljoin}}

\def\!finish/{\ignorespaces}

\def\!hcurve #1 #2 {%
  \edef\!hxS{#1}%
  \edef\!hyS{#2}%
  \!hjoin}
\def\!hjoin#1 #2 {%
  \putrectangle corners at {\!hxS} {\!hyS} and {#1} {#2}
  \edef\!hxS{#1}%
  \!ifnextchar/{\!finish}{\!hjoin}}

\def\vshade #1 #2 #3 {%
  \!startvshade (#1,#2,#3)
  \!Shadewhat}

\def\hshade #1 #2 #3 {%
  \!starthshade (#1,#2,#3)
  \!Shadewhat}

\def\!Shadewhat{%
  \futurelet\!nextchar\!Shade}
\def\!Shade{%
  \if <\!nextchar
    \def\!nextShade{\!!Shade}%
  \else
    \if /\!nextchar
      \def\!nextShade{\!finish}%
    \else
      \def\!nextShade{\!!!Shade}%
    \fi
  \fi
  \!nextShade}
\def\!!lShade<#1> #2 #3 #4 {%
  \!lshade <#1> (#2,#3,#4)                 
  \!Shadewhat}
\def\!!!lShade#1 #2 #3 {%
  \!lshade (#1,#2,#3)
  \!Shadewhat}
\def\!!qShade<#1> #2 #3 #4 #5 #6 #7 {%
  \!qshade <#1> (#2,#3,#4) (#5,#6,#7)      
  \!Shadewhat}
\def\!!!qShade#1 #2 #3 #4 #5 #6 {%
  \!qshade (#1,#2,#3) (#4,#5,#6)
  \!Shadewhat}

\setlinear

\def\setdashpattern <#1>{%
  \def\!Flist{}\def\!Blist{}\def\!UDlist{}%
  \!countA=0
  \!ecfor\!item:=#1\do{%
    \!dimenA=\!item\relax
    \expandafter\!rightappend\the\!dimenA\withCS{\\}\to\!UDlist%
    \advance\!countA  1
    \ifodd\!countA
      \expandafter\!rightappend\the\!dimenA\withCS{\!Rule}\to\!Flist%
      \expandafter\!leftappend\the\!dimenA\withCS{\!Rule}\to\!Blist%
    \else
      \expandafter\!rightappend\the\!dimenA\withCS{\!Skip}\to\!Flist%
      \expandafter\!leftappend\the\!dimenA\withCS{\!Skip}\to\!Blist%
    \fi}%
  \!leaderlength=\!zpt
  \def\!Rule##1{\advance\!leaderlength  ##1}%
  \def\!Skip##1{\advance\!leaderlength  ##1}%
  \!Flist%
  \ifdim\!leaderlength>\!zpt
  \else
    \def\!Flist{\!Skip{24in}}\def\!Blist{\!Skip{24in}}\ignorespaces
    \def\!UDlist{\\{\!zpt}\\{24in}}\ignorespaces
    \!leaderlength=24in
  \fi
  \!dashingon}

\def\!dashingon{%
  \def\!advancedashing{\!!advancedashing}%
  \def\!drawlinearsegment{\!lineardashed}%
  \def\!puthline{\!putdashedhline}%
  \def\!putvline{\!putdashedvline}%
  \ignorespaces}%
\def\!dashingoff{%
  \def\!advancedashing{\relax}%
  \def\!drawlinearsegment{\!linearsolid}%
  \def\!puthline{\!putsolidhline}%
  \def\!putvline{\!putsolidvline}%
  \ignorespaces}

\def\setdots{%
  \!ifnextchar<{\!setdots}{\!setdots<5pt>}}
\def\!setdots<#1>{%
  \!dimenB=#1\advance\!dimenB -\plotsymbolspacing
  \ifdim\!dimenB<\!zpt
    \!dimenB=\!zpt
  \fi
\setdashpattern <\plotsymbolspacing,\!dimenB>}

\def\setdotsnear <#1> for <#2>{%
  \!dimenB=#2\relax  \advance\!dimenB -.05pt
  \!dimenC=#1\relax  \!countA=\!dimenC
  \!dimenD=\!dimenB  \advance\!dimenD .5\!dimenC  \!countB=\!dimenD
  \divide \!countB  \!countA
  \ifnum 1>\!countB
    \!countB=1
  \fi
  \divide\!dimenB  \!countB
  \setdots <\!dimenB>}

\def\setdashes{%
  \!ifnextchar<{\!setdashes}{\!setdashes<5pt>}}
\def\!setdashes<#1>{\setdashpattern <#1,#1>}

\def\setdashesnear <#1> for <#2>{%
  \!dimenB=#2\relax
  \!dimenC=#1\relax  \!countA=\!dimenC
  \!dimenD=\!dimenB  \advance\!dimenD .5\!dimenC  \!countB=\!dimenD
  \divide \!countB  \!countA
  \ifodd \!countB
  \else
    \advance \!countB  1
  \fi
  \divide\!dimenB  \!countB
  \setdashes <\!dimenB>}

\def\setsolid{%
  \def\!Flist{\!Rule{24in}}\def\!Blist{\!Rule{24in}}%
  \def\!UDlist{\\{24in}\\{\!zpt}}%
  \!dashingoff}
\setsolid

\def\!divide#1#2#3{%
  \!dimenB=#1
  \!dimenC=#2
  \!dimenD=\!dimenB
  \divide \!dimenD \!dimenC
  \!dimenA=\!dimenD
  \multiply\!dimenD \!dimenC
  \advance\!dimenB -\!dimenD
  \!dimenD=\!dimenC
    \ifdim\!dimenD<\!zpt \!dimenD=-\!dimenD
  \fi
  \ifdim\!dimenD<64pt
    \!divstep[\!tfs]\!divstep[\!tfs]%
  \else
    \!!divide
  \fi
  #3=\!dimenA\ignorespaces}

\def\!!divide{%
  \ifdim\!dimenD<256pt
    \!divstep[64]\!divstep[32]\!divstep[32]%
  \else
    \!divstep[8]\!divstep[8]\!divstep[8]\!divstep[8]\!divstep[8]%
    \!dimenA=2\!dimenA
  \fi}

\def\!divstep[#1]{
  \!dimenB=#1\!dimenB
  \!dimenD=\!dimenB
    \divide \!dimenD by \!dimenC
  \!dimenA=#1\!dimenA
    \advance\!dimenA by \!dimenD%
  \multiply\!dimenD by \!dimenC
    \advance\!dimenB by -\!dimenD}

\def\Divide <#1> by <#2> forming <#3> {%
  \!divide{#1}{#2}{#3}}

\def\circulararc{%
  \ellipticalarc axes ratio 1:1 }

\def\ellipticalarc axes ratio #1:#2 #3 degrees from #4 #5 center at #6 #7 {%
  \!angle=#3pt\relax
  \ifdim\!angle>\!zpt
    \def\!sign{}
  \else
    \def\!sign{-}\!angle=-\!angle
  \fi
  \!xxloc=\!M{#6}\!xunit
  \!yyloc=\!M{#7}\!yunit
  \!xxS=\!M{#4}\!xunit
  \!yyS=\!M{#5}\!yunit
  \advance\!xxS -\!xxloc
  \advance\!yyS -\!yyloc
  \!divide\!xxS{#1pt}\!xxS 
  \!divide\!yyS{#2pt}\!yyS 
  \let\!MC=\!M
  \!setdimenmode
  \!xS=#1\!xxS  \advance\!xS\!xxloc
  \!yS=#2\!yyS  \advance\!yS\!yyloc
  \!start (\!xS,\!yS)%
  \!loop\ifdim\!angle>14.9999pt
    \!rotate(\!xxS,\!yyS)by(\!cos,\!sign\!sin)to(\!xxM,\!yyM)
    \!rotate(\!xxM,\!yyM)by(\!cos,\!sign\!sin)to(\!xxE,\!yyE)
    \!xM=#1\!xxM  \advance\!xM\!xxloc  \!yM=#2\!yyM  \advance\!yM\!yyloc
    \!xE=#1\!xxE  \advance\!xE\!xxloc  \!yE=#2\!yyE  \advance\!yE\!yyloc
    \!qjoin (\!xM,\!yM) (\!xE,\!yE)
    \!xxS=\!xxE  \!yyS=\!yyE
    \advance \!angle -15pt
  \repeat
  \ifdim\!angle>\!zpt
    \!angle=100.53096\!angle
    \divide \!angle 360 
    \!sinandcos\!angle\!!sin\!!cos
    \!rotate(\!xxS,\!yyS)by(\!!cos,\!sign\!!sin)to(\!xxM,\!yyM)
    \!rotate(\!xxM,\!yyM)by(\!!cos,\!sign\!!sin)to(\!xxE,\!yyE)
    \!xM=#1\!xxM  \advance\!xM\!xxloc  \!yM=#2\!yyM  \advance\!yM\!yyloc
    \!xE=#1\!xxE  \advance\!xE\!xxloc  \!yE=#2\!yyE  \advance\!yE\!yyloc
    \!qjoin (\!xM,\!yM) (\!xE,\!yE)
  \fi
  \let\!M=\!MC
  \ignorespaces}

\def\!rotate(#1,#2)by(#3,#4)to(#5,#6){%
  \!dimenA=#3#1\advance \!dimenA -#4#2
  \!dimenB=#3#2\advance \!dimenB  #4#1
  \divide \!dimenA 32  \divide \!dimenB 32
  #5=\!dimenA  #6=\!dimenB
  \ignorespaces}
\def\!sin{4.17684}
\def\!cos{31.72624}

\def\!sinandcos#1#2#3{%
 \!dimenD=#1
 \!dimenA=\!dimenD
 \!dimenB=32pt
 \!removept\!dimenD\!value
 \!dimenC=\!dimenD
 \!dimenC=\!value\!dimenC \divide\!dimenC by 64 
 \advance\!dimenB by -\!dimenC
 \!dimenC=\!value\!dimenC \divide\!dimenC by 96 
 \advance\!dimenA by -\!dimenC
 \!dimenC=\!value\!dimenC \divide\!dimenC by 128 
 \advance\!dimenB by \!dimenC%
 \!removept\!dimenA#2
 \!removept\!dimenB#3
 \ignorespaces}

\def\putrule#1from #2 #3 to #4 #5 {%
  \!xloc=\!M{#2}\!xunit  \!xxloc=\!M{#4}\!xunit%
  \!yloc=\!M{#3}\!yunit  \!yyloc=\!M{#5}\!yunit%
  \!dxpos=\!xxloc  \advance\!dxpos by -\!xloc
  \!dypos=\!yyloc  \advance\!dypos by -\!yloc
  \ifdim\!dypos=\!zpt
    \def\!!Line{\!puthline{#1}}\ignorespaces
  \else
    \ifdim\!dxpos=\!zpt
      \def\!!Line{\!putvline{#1}}\ignorespaces
    \else
       \def\!!Line{}
    \fi
  \fi
  \let\!ML=\!M
  \!setdimenmode
  \!!Line%
  \let\!M=\!ML
  \ignorespaces}

\def\!putsolidhline#1{%
  \ifdim\!dxpos>\!zpt
    \put{\!hline\!dxpos}#1[l] at {\!xloc} {\!yloc}
  \else
    \put{\!hline{-\!dxpos}}#1[l] at {\!xxloc} {\!yyloc}
  \fi
  \ignorespaces}

\def\!putsolidvline#1{%
  \ifdim\!dypos>\!zpt
    \put{\!vline\!dypos}#1[b] at {\!xloc} {\!yloc}
  \else
    \put{\!vline{-\!dypos}}#1[b] at {\!xxloc} {\!yyloc}
  \fi
  \ignorespaces}

\def\!hline#1{\hbox to #1{\leaders \hrule height\linethickness\hfill}}
\def\!vline#1{\vbox to #1{\leaders \vrule width\linethickness\vfill}}

\def\!putdashedhline#1{%
  \ifdim\!dxpos>\!zpt
    \!DLsetup\!Flist\!dxpos
    \put{\hbox to \!totalleaderlength{\!hleaders}\!hpartialpattern\!Rtrunc}
      #1[l] at {\!xloc} {\!yloc}
  \else
    \!DLsetup\!Blist{-\!dxpos}
    \put{\!hpartialpattern\!Ltrunc\hbox to \!totalleaderlength{\!hleaders}}
      #1[r] at {\!xloc} {\!yloc}
  \fi
  \ignorespaces}

\def\!putdashedvline#1{%
  \!dypos=-\!dypos
  \ifdim\!dypos>\!zpt
    \!DLsetup\!Flist\!dypos
    \put{\vbox{\vbox to \!totalleaderlength{\!vleaders}
      \!vpartialpattern\!Rtrunc}}#1[t] at {\!xloc} {\!yloc}
  \else
    \!DLsetup\!Blist{-\!dypos}
    \put{\vbox{\!vpartialpattern\!Ltrunc
      \vbox to \!totalleaderlength{\!vleaders}}}#1[b] at {\!xloc} {\!yloc}
  \fi
  \ignorespaces}

\def\!DLsetup#1#2{
  \let\!RSlist=#1
  \!countB=#2
  \!countA=\!leaderlength
  \divide\!countB by \!countA
  \!totalleaderlength=\!countB\!leaderlength
  \!Rresiduallength=#2%
  \advance \!Rresiduallength by -\!totalleaderlength
  \!Lresiduallength=\!leaderlength
  \advance \!Lresiduallength by -\!Rresiduallength
  \ignorespaces}

\def\!hleaders{%
  \def\!Rule##1{\vrule height\linethickness width##1}%
  \def\!Skip##1{\hskip##1}%
  \leaders\hbox{\!RSlist}\hfill}

\def\!hpartialpattern#1{%
  \!dimenA=\!zpt \!dimenB=\!zpt
  \def\!Rule##1{#1{##1}\vrule height\linethickness width\!dimenD}%
  \def\!Skip##1{#1{##1}\hskip\!dimenD}%
  \!RSlist}

\def\!vleaders{%
  \def\!Rule##1{\hrule width\linethickness height##1}%
  \def\!Skip##1{\vskip##1}%
  \leaders\vbox{\!RSlist}\vfill}

\def\!vpartialpattern#1{%
  \!dimenA=\!zpt \!dimenB=\!zpt
  \def\!Rule##1{#1{##1}\hrule width\linethickness height\!dimenD}%
  \def\!Skip##1{#1{##1}\vskip\!dimenD}%
  \!RSlist}

\def\!Rtrunc#1{\!trunc{#1}>\!Rresiduallength}
\def\!Ltrunc#1{\!trunc{#1}<\!Lresiduallength}

\def\!trunc#1#2#3{%
  \!dimenA=\!dimenB
  \advance\!dimenB by #1%
  \!dimenD=\!dimenB  \ifdim\!dimenD#2#3\!dimenD=#3\fi
  \!dimenC=\!dimenA  \ifdim\!dimenC#2#3\!dimenC=#3\fi
  \advance \!dimenD by -\!dimenC}

\def\!start (#1,#2){%
  \!plotxorigin=\!xorigin  \advance \!plotxorigin by \!plotsymbolxshift
  \!plotyorigin=\!yorigin  \advance \!plotyorigin by \!plotsymbolyshift
  \!xS=\!M{#1}\!xunit \!yS=\!M{#2}\!yunit
  \!rotateaboutpivot\!xS\!yS
  \!copylist\!UDlist\to\!!UDlist
  \!getnextvalueof\!downlength\from\!!UDlist
  \!distacross=\!zpt
  \!intervalno=0 
  \global\totalarclength=\!zpt
  \ignorespaces}

\def\!ljoin (#1,#2){%
  \advance\!intervalno by 1
  \!xE=\!M{#1}\!xunit \!yE=\!M{#2}\!yunit
  \!rotateaboutpivot\!xE\!yE
  \!xdiff=\!xE \advance \!xdiff by -\!xS
  \!ydiff=\!yE \advance \!ydiff by -\!yS
  \!Pythag\!xdiff\!ydiff\!arclength
  \global\advance \totalarclength by \!arclength%
  \!drawlinearsegment
  \!xS=\!xE \!yS=\!yE
  \ignorespaces}

\def\!linearsolid{%
  \!npoints=\!arclength
  \!countA=\plotsymbolspacing
  \divide\!npoints by \!countA
  \ifnum \!npoints<1
    \!npoints=1
  \fi
  \divide\!xdiff by \!npoints
  \divide\!ydiff by \!npoints
  \!xpos=\!xS \!ypos=\!yS
  \loop\ifnum\!npoints>-1
    \!plotifinbounds
    \advance \!xpos by \!xdiff
    \advance \!ypos by \!ydiff
    \advance \!npoints by -1
  \repeat
  \ignorespaces}

\def\!lineardashed{%
  \ifdim\!distacross>\!arclength
    \advance \!distacross by -\!arclength  
  \else
    \loop\ifdim\!distacross<\!arclength
      \!divide\!distacross\!arclength\!dimenA
      \!removept\!dimenA\!t
      \!xpos=\!t\!xdiff \advance \!xpos by \!xS
      \!ypos=\!t\!ydiff \advance \!ypos by \!yS
      \!plotifinbounds
      \advance\!distacross by \plotsymbolspacing
      \!advancedashing
    \repeat
    \advance \!distacross by -\!arclength
  \fi
  \ignorespaces}

\def\!!advancedashing{%
  \advance\!downlength by -\plotsymbolspacing
  \ifdim \!downlength>\!zpt
  \else
    \advance\!distacross by \!downlength
    \!getnextvalueof\!uplength\from\!!UDlist
    \advance\!distacross by \!uplength
    \!getnextvalueof\!downlength\from\!!UDlist
  \fi}

\def\inboundscheckoff{%
  \def\!plotifinbounds{\!plot(\!xpos,\!ypos)}%
  \def\!initinboundscheck{\relax}\ignorespaces}

\inboundscheckoff

\def\!!plotifinbounds{%
  \ifdim \!xpos<\!checkleft
  \else
    \ifdim \!xpos>\!checkright
    \else
      \ifdim \!ypos<\!checkbot
      \else
         \ifdim \!ypos>\!checktop
         \else
           \!plot(\!xpos,\!ypos)
         \fi
      \fi
    \fi
  \fi}

\def\!!initinboundscheck{%
  \!checkleft=\!arealloc     \advance\!checkleft by \!xorigin
  \!checkright=\!arearloc    \advance\!checkright by \!xorigin
  \!checkbot=\!areabloc      \advance\!checkbot by \!yorigin
  \!checktop=\!areatloc      \advance\!checktop by \!yorigin}

\def\!logten#1#2{%
  \expandafter\!!logten#1\!nil
  \!removept\!dimenF#2%
  \ignorespaces}

\def\!!logten#1#2\!nil{%
  \if -#1%
    \!dimenF=\!zpt
    \def\!next{\ignorespaces}%
  \else
    \if +#1%
      \def\!next{\!!logten#2\!nil}%
    \else
      \if .#1%
        \def\!next{\!!logten0.#2\!nil}%
      \else
        \def\!next{\!!!logten#1#2..\!nil}%
      \fi
    \fi
  \fi
  \!next}

\def\!!!logten#1#2.#3.#4\!nil{%
  \!dimenF=1pt 
  \if 0#1%
    \!!logshift#3pt 
  \else 
    \!logshift#2/
    \!dimenE=#1.#2#3pt 
  \fi 
  \ifdim \!dimenE<\!rootten
    \multiply \!dimenE 10 
    \advance  \!dimenF -1pt
  \fi
  \!dimenG=\!dimenE
    \advance\!dimenG 10pt
  \advance\!dimenE -10pt 
  \multiply\!dimenE 10 
  \!divide\!dimenE\!dimenG\!dimenE
  \!removept\!dimenE\!t
  \!dimenG=\!t\!dimenE
  \!removept\!dimenG\!tt
  \!dimenH=\!tt\!tenAe
    \divide\!dimenH 100
  \advance\!dimenH \!tenAc
  \!dimenH=\!tt\!dimenH
    \divide\!dimenH 100
  \advance\!dimenH \!tenAa
  \!dimenH=\!t\!dimenH
    \divide\!dimenH 100 
  \advance\!dimenF \!dimenH}

\def\!logshift#1{%
  \if #1/%
    \def\!next{\ignorespaces}%
  \else
    \advance\!dimenF 1pt
    \def\!next{\!logshift}%
  \fi
  \!next}

 \def\!!logshift#1{%
   \advance\!dimenF -1pt
   \if 0#1%
     \def\!next{\!!logshift}%
   \else
     \if p#1%
       \!dimenF=1pt
       \def\!next{\!dimenE=1p}%
     \else
       \def\!next{\!dimenE=#1.}%
     \fi
   \fi
   \!next}

\def\beginpicture{%
  \setbox\!picbox=\hbox\bgroup%
  \!xleft=\maxdimen
  \!xright=-\maxdimen
  \!ybot=\maxdimen
  \!ytop=-\maxdimen}

\def\endpicture{%
  \ifdim\!xleft=\maxdimen
    \!xleft=\!zpt \!xright=\!zpt \!ybot=\!zpt \!ytop=\!zpt
  \fi
  \global\!Xleft=\!xleft \global\!Xright=\!xright
  \global\!Ybot=\!ybot \global\!Ytop=\!ytop
  \egroup%
  \ht\!picbox=\!Ytop  \dp\!picbox=-\!Ybot
  \ifdim\!Ybot>\!zpt
  \else
    \ifdim\!Ytop<\!zpt
      \!Ybot=\!Ytop
    \else
      \!Ybot=\!zpt
    \fi
  \fi
  \hbox{\kern-\!Xleft\lower\!Ybot\box\!picbox\kern\!Xright}}

\def\endpicturesave <#1,#2>{%
  \endpicture \global #1=\!Xleft \global #2=\!Ybot \ignorespaces}

\def\setcoordinatesystem{%
  \!ifnextchar{u}{\!getlengths }
    {\!getlengths units <\!xunit,\!yunit>}}
\def\!getlengths units <#1,#2>{%
  \!xunit=#1\relax
  \!yunit=#2\relax
  \!ifcoordmode
    \let\!SCnext=\!SCccheckforRP
  \else
    \let\!SCnext=\!SCdcheckforRP
  \fi
  \!SCnext}
\def\!SCccheckforRP{%
  \!ifnextchar{p}{\!cgetreference }
    {\!cgetreference point at {\!xref} {\!yref} }}
\def\!cgetreference point at #1 #2 {%
  \edef\!xref{#1}\edef\!yref{#2}%
  \!xorigin=\!xref\!xunit  \!yorigin=\!yref\!yunit
  \!initinboundscheck 
  \ignorespaces}
\def\!SCdcheckforRP{%
  \!ifnextchar{p}{\!dgetreference}%
    {\ignorespaces}}
\def\!dgetreference point at #1 #2 {%
  \!xorigin=#1\relax  \!yorigin=#2\relax
  \ignorespaces}

\long\def\put#1#2 at #3 #4 {%
  \!setputobject{#1}{#2}%
  \!xpos=\!M{#3}\!xunit  \!ypos=\!M{#4}\!yunit
  \!rotateaboutpivot\!xpos\!ypos%
  \advance\!xpos -\!xorigin  \advance\!xpos -\!xshift
  \advance\!ypos -\!yorigin  \advance\!ypos -\!yshift
  \kern\!xpos\raise\!ypos\box\!putobject\kern-\!xpos%
  \!doaccounting\ignorespaces}

\long\def\multiput #1#2 at {%
  \!setputobject{#1}{#2}%
  \!ifnextchar"{\!putfromfile}{\!multiput}}
\def\!putfromfile"#1"{%
  \expandafter\!multiput \input #1 /}
\def\!multiput{%
  \futurelet\!nextchar\!!multiput}
\def\!!multiput{%
  \if *\!nextchar
    \def\!nextput{\!alsoby}%
  \else
    \if /\!nextchar
      \def\!nextput{\!finishmultiput}%
    \else
      \def\!nextput{\!alsoat}%
    \fi
  \fi
  \!nextput}
\def\!finishmultiput/{%
  \setbox\!putobject=\hbox{}%
  \ignorespaces}

\def\!alsoat#1 #2 {%
  \!xpos=\!M{#1}\!xunit  \!ypos=\!M{#2}\!yunit
  \!rotateaboutpivot\!xpos\!ypos%
  \advance\!xpos -\!xorigin  \advance\!xpos -\!xshift
  \advance\!ypos -\!yorigin  \advance\!ypos -\!yshift
  \kern\!xpos\raise\!ypos\copy\!putobject\kern-\!xpos%
  \!doaccounting
  \!multiput}

\def\!alsoby*#1 #2 #3 {%
  \!dxpos=\!M{#2}\!xunit \!dypos=\!M{#3}\!yunit
  \!rotateonly\!dxpos\!dypos
  \!ntemp=#1%
  \!!loop\ifnum\!ntemp>0
    \advance\!xpos by \!dxpos  \advance\!ypos by \!dypos
    \kern\!xpos\raise\!ypos\copy\!putobject\kern-\!xpos%
    \advance\!ntemp by -1
  \repeat
  \!doaccounting
  \!multiput}

\def\accountingon{\def\!doaccounting{\!!doaccounting}\ignorespaces}

\accountingon
\def\!!doaccounting{%
  \!xtemp=\!xpos
  \!ytemp=\!ypos
  \ifdim\!xtemp<\!xleft
     \!xleft=\!xtemp
  \fi
  \advance\!xtemp by  \!wd
  \ifdim\!xright<\!xtemp
    \!xright=\!xtemp
  \fi
  \advance\!ytemp by -\!dp
  \ifdim\!ytemp<\!ybot
    \!ybot=\!ytemp
  \fi
  \advance\!ytemp by  \!dp
  \advance\!ytemp by  \!ht
  \ifdim\!ytemp>\!ytop
    \!ytop=\!ytemp
  \fi}

\long\def\!setputobject#1#2{%
  \setbox\!putobject=\hbox{#1}%
  \!ht=\ht\!putobject  \!dp=\dp\!putobject  \!wd=\wd\!putobject
  \wd\!putobject=\!zpt
  \!xshift=.5\!wd   \!yshift=.5\!ht   \advance\!yshift by -.5\!dp
  \edef\!putorientation{#2}%
  \expandafter\!SPOreadA\!putorientation[]\!nil%
  \expandafter\!SPOreadB\!putorientation<\!zpt,\!zpt>\!nil\ignorespaces}

\def\!SPOreadA#1[#2]#3\!nil{\!etfor\!orientation:=#2\do\!SPOreviseshift}

\def\!SPOreadB#1<#2,#3>#4\!nil{\advance\!xshift by -#2\advance\!yshift by -#3}

\def\!SPOreviseshift{%
  \if l\!orientation
    \!xshift=\!zpt
  \else
    \if r\!orientation
      \!xshift=\!wd
    \else
      \if b\!orientation
        \!yshift=-\!dp
      \else
        \if B\!orientation
          \!yshift=\!zpt
        \else
          \if t\!orientation
            \!yshift=\!ht
          \fi
        \fi
      \fi
    \fi
  \fi}

\long\def\!dimenput#1#2(#3,#4){%
  \!setputobject{#1}{#2}%
  \!xpos=#3\advance\!xpos by -\!xshift
  \!ypos=#4\advance\!ypos by -\!yshift
  \kern\!xpos\raise\!ypos\box\!putobject\kern-\!xpos%
  \!doaccounting\ignorespaces}

\def\!setdimenmode{%
  \let\!M=\!M!!\ignorespaces}
\def\!setcoordmode{%
  \let\!M=\!M!\ignorespaces}
\def\!ifcoordmode{%
  \ifx \!M \!M!}
\def\!ifdimenmode{%
  \ifx \!M \!M!!}
\def\!M!#1#2{#1#2}
\def\!M!!#1#2{#1}
\!setcoordmode
\let\setdimensionmode=\!setdimenmode
\let\setcoordinatemode=\!setcoordmode

\def\!stack[#1]{%
  \let\!lglue=\hfill \let\!rglue=\hfill
  \expandafter\let\csname !#1glue\endcsname=\relax
  \!ifnextchar<{\!!stack}{\!!stack<\stackleading>}}
\def\!!stack<#1>#2{%
  \vbox{\def\!valueslist{}\!ecfor\!value:=#2\do{%
    \expandafter\!rightappend\!value\withCS{\\}\to\!valueslist}%
    \!lop\!valueslist\to\!value
    \let\\=\cr\lineskiplimit=\maxdimen\lineskip=#1%
    \baselineskip=-1000pt\halign{\!lglue##\!rglue\cr \!value\!valueslist\cr}}%
  \ignorespaces}

\def\!lines[#1]#2{%
  \let\!lglue=\hfill \let\!rglue=\hfill
  \expandafter\let\csname !#1glue\endcsname=\relax
  \vbox{\halign{\!lglue##\!rglue\cr #2\crcr}}%
  \ignorespaces}

\def\!Lines[#1]#2{%
  \let\!lglue=\hfill \let\!rglue=\hfill
  \expandafter\let\csname !#1glue\endcsname=\relax
  \vtop{\halign{\!lglue##\!rglue\cr #2\crcr}}%
  \ignorespaces}

\def\setplotsymbol(#1#2){%
  \!setputobject{#1}{#2}
  \setbox\!plotsymbol=\box\!putobject%
  \!plotsymbolxshift=\!xshift
  \!plotsymbolyshift=\!yshift
  \ignorespaces}

\setplotsymbol({\fiverm .})

\def\!!plot(#1,#2){%
  \!dimenA=-\!plotxorigin \advance \!dimenA by #1
  \!dimenB=-\!plotyorigin \advance \!dimenB by #2
  \kern\!dimenA\raise\!dimenB\copy\!plotsymbol\kern-\!dimenA%
  \ignorespaces}

\def\!!!plot(#1,#2){%
  \!dimenA=-\!plotxorigin \advance \!dimenA by #1
  \!dimenB=-\!plotyorigin \advance \!dimenB by #2
  \kern\!dimenA\raise\!dimenB\copy\!plotsymbol\kern-\!dimenA%
  \!countE=\!dimenA
  \!countF=\!dimenB
  \immediate\write\!replotfile{\the\!countE,\the\!countF.}%
  \ignorespaces}

\def\savelinesandcurves on "#1" {%
  \immediate\closeout\!replotfile
  \immediate\openout\!replotfile=#1%
  \let\!plot=\!!!plot}

\def\dontsavelinesandcurves {%
  \let\!plot=\!!plot}
\dontsavelinesandcurves

{\catcode`\%=11\xdef\!Commentsignal{
\def\writesavefile#1 {%
  \immediate\write\!replotfile{\!Commentsignal #1}%
  \ignorespaces}

\def\replot"#1" {%
  \expandafter\!replot\input #1 /}
\def\!replot#1,#2. {%
  \!dimenA=#1sp
  \kern\!dimenA\raise#2sp\copy\!plotsymbol\kern-\!dimenA
  \futurelet\!nextchar\!!replot}
\def\!!replot{%
  \if /\!nextchar
    \def\!next{\!finish}%
  \else
    \def\!next{\!replot}%
  \fi
  \!next}

\def\!Pythag#1#2#3{%
  \!dimenE=#1\relax
  \ifdim\!dimenE<\!zpt
    \!dimenE=-\!dimenE
  \fi
  \!dimenF=#2\relax
  \ifdim\!dimenF<\!zpt
    \!dimenF=-\!dimenF
  \fi
  \advance \!dimenF by \!dimenE
  \ifdim\!dimenF=\!zpt
    \!dimenG=\!zpt
  \else
    \!divide{8\!dimenE}\!dimenF\!dimenE
    \advance\!dimenE by -4pt
      \!dimenE=2\!dimenE
    \!removept\!dimenE\!!t
    \!dimenE=\!!t\!dimenE
    \advance\!dimenE by 64pt
    \divide \!dimenE by 2
    \!dimenH=7pt
    \!!Pythag\!!Pythag\!!Pythag
    \!removept\!dimenH\!!t
    \!dimenG=\!!t\!dimenF
    \divide\!dimenG by 8
  \fi
  #3=\!dimenG
  \ignorespaces}

\def\!!Pythag{
  \!divide\!dimenE\!dimenH\!dimenI
  \advance\!dimenH by \!dimenI
    \divide\!dimenH by 2}

\def\placehypotenuse for <#1> and <#2> in <#3> {%
  \!Pythag{#1}{#2}{#3}}

\def\!qjoin (#1,#2) (#3,#4){%
  \advance\!intervalno by 1
  \!ifcoordmode
    \edef\!xmidpt{#1}\edef\!ymidpt{#2}%
  \else
    \!dimenA=#1\relax \edef\!xmidpt{\the\!dimenA}%
    \!dimenA=#2\relax \edef\!ymidpt{\the\!dimenA}%
  \fi
  \!xM=\!M{#1}\!xunit  \!yM=\!M{#2}\!yunit   \!rotateaboutpivot\!xM\!yM
  \!xE=\!M{#3}\!xunit  \!yE=\!M{#4}\!yunit   \!rotateaboutpivot\!xE\!yE
  \!dimenA=\!xM  \advance \!dimenA by -\!xS
  \!dimenB=\!xE  \advance \!dimenB by -\!xM
  \!xB=3\!dimenA \advance \!xB by -\!dimenB
  \!xC=2\!dimenB \advance \!xC by -2\!dimenA
  \!dimenA=\!yM  \advance \!dimenA by -\!yS%
  \!dimenB=\!yE  \advance \!dimenB by -\!yM%
  \!yB=3\!dimenA \advance \!yB by -\!dimenB%
  \!yC=2\!dimenB \advance \!yC by -2\!dimenA%
  \!xprime=\!xB  \!yprime=\!yB
  \!dxprime=.5\!xC  \!dyprime=.5\!yC
  \!getf \!midarclength=\!dimenA
  \!getf \advance \!midarclength by 4\!dimenA
  \!getf \advance \!midarclength by \!dimenA
  \divide \!midarclength by 12
  \!arclength=\!dimenA
  \!getf \advance \!arclength by 4\!dimenA
  \!getf \advance \!arclength by \!dimenA
  \divide \!arclength by 12
  \advance \!arclength by \!midarclength
  \global\advance \totalarclength by \!arclength
  \ifdim\!distacross>\!arclength
    \advance \!distacross by -\!arclength
  \else
    \!initinverseinterp
    \loop\ifdim\!distacross<\!arclength
      \!inverseinterp
      \!xpos=\!t\!xC \advance\!xpos by \!xB
        \!xpos=\!t\!xpos \advance \!xpos by \!xS
      \!ypos=\!t\!yC \advance\!ypos by \!yB
        \!ypos=\!t\!ypos \advance \!ypos by \!yS
      \!plotifinbounds
      \advance\!distacross \plotsymbolspacing
      \!advancedashing
    \repeat
    \advance \!distacross by -\!arclength
  \fi
  \!xS=\!xE
  \!yS=\!yE
  \ignorespaces}

\def\!getf{\!Pythag\!xprime\!yprime\!dimenA%
  \advance\!xprime by \!dxprime
  \advance\!yprime by \!dyprime}

\def\!initinverseinterp{%
  \ifdim\!arclength>\!zpt
    \!divide{8\!midarclength}\!arclength\!dimenE
    \ifdim\!dimenE<\!wmin \!setinverselinear
    \else
      \ifdim\!dimenE>\!wmax \!setinverselinear
      \else
        \def\!inverseinterp{\!inversequad}\ignorespaces
         \!removept\!dimenE\!Ew
         \!dimenF=-\!Ew\!dimenE
         \advance\!dimenF by 32pt
         \!dimenG=8pt
         \advance\!dimenG by -\!dimenE
         \!dimenG=\!Ew\!dimenG
         \!divide\!dimenF\!dimenG\!beta
         \!gamma=1pt
         \advance \!gamma by -\!beta
      \fi
    \fi
  \fi
  \ignorespaces}

\def\!inversequad{%
  \!divide\!distacross\!arclength\!dimenG
  \!removept\!dimenG\!v
  \!dimenG=\!v\!gamma
  \advance\!dimenG by \!beta
  \!dimenG=\!v\!dimenG
  \!removept\!dimenG\!t}

\def\!setinverselinear{%
  \def\!inverseinterp{\!inverselinear}%
  \divide\!dimenE by 8 \!removept\!dimenE\!t
  \!countC=\!intervalno \multiply \!countC 2
  \!countB=\!countC     \advance \!countB -1
  \!countA=\!countB     \advance \!countA -1
  \wlog{\the\!countB th point (\!xmidpt,\!ymidpt) being plotted
    doesn't lie in the}%
  \wlog{ middle third of the arc between the \the\!countA th
    and \the\!countC th points:}%
  \wlog{ [arc length \the\!countA\space to \the\!countB]/[arc length
    \the \!countA\space to \the\!countC]=\!t.}%
  \ignorespaces}

\def\!inverselinear{%
  \!divide\!distacross\!arclength\!dimenG
  \!removept\!dimenG\!t}

\def\startrotation{%
  \let\!rotateaboutpivot=\!!rotateaboutpivot
  \let\!rotateonly=\!!rotateonly
  \!ifnextchar{b}{\!getsincos }%
    {\!getsincos by {\!cosrotationangle} {\!sinrotationangle} }}
\def\!getsincos by #1 #2 {%
  \edef\!cosrotationangle{#1}%
  \edef\!sinrotationangle{#2}%
  \!ifcoordmode
    \let\!ROnext=\!ccheckforpivot
  \else
    \let\!ROnext=\!dcheckforpivot
  \fi
  \!ROnext}
\def\!ccheckforpivot{%
  \!ifnextchar{a}{\!cgetpivot}%
    {\!cgetpivot about {\!xpivotcoord} {\!ypivotcoord} }}
\def\!cgetpivot about #1 #2 {%
  \edef\!xpivotcoord{#1}%
  \edef\!ypivotcoord{#2}%
  \!xpivot=#1\!xunit  \!ypivot=#2\!yunit
  \ignorespaces}
\def\!dcheckforpivot{%
  \!ifnextchar{a}{\!dgetpivot}{\ignorespaces}}
\def\!dgetpivot about #1 #2 {%
  \!xpivot=#1\relax  \!ypivot=#2\relax
  \ignorespaces}

\def\stoprotation{%
  \let\!rotateaboutpivot=\!!!rotateaboutpivot
  \let\!rotateonly=\!!!rotateonly
  \ignorespaces}

\def\!!rotateaboutpivot#1#2{%
  \!dimenA=#1\relax  \advance\!dimenA -\!xpivot
  \!dimenB=#2\relax  \advance\!dimenB -\!ypivot
  \!dimenC=\!cosrotationangle\!dimenA
    \advance \!dimenC -\!sinrotationangle\!dimenB
  \!dimenD=\!cosrotationangle\!dimenB
    \advance \!dimenD  \!sinrotationangle\!dimenA
  \advance\!dimenC \!xpivot  \advance\!dimenD \!ypivot
  #1=\!dimenC  #2=\!dimenD
  \ignorespaces}

\def\!!rotateonly#1#2{%
  \!dimenA=#1\relax  \!dimenB=#2\relax
  \!dimenC=\!cosrotationangle\!dimenA
    \advance \!dimenC -\!rotsign\!sinrotationangle\!dimenB
  \!dimenD=\!cosrotationangle\!dimenB
    \advance \!dimenD  \!rotsign\!sinrotationangle\!dimenA
  #1=\!dimenC  #2=\!dimenD
  \ignorespaces}
\def\!rotsign{}
\def\!!!rotateaboutpivot#1#2{\relax}
\def\!!!rotateonly#1#2{\relax}
\stoprotation

\def\!reverserotateonly#1#2{%
  \def\!rotsign{-}%
  \!rotateonly{#1}{#2}%
  \def\!rotsign{}%
  \ignorespaces}

\def\!getspan span <#1>{%
  \!dshade=#1\relax
  \!ifcoordmode
    \let\!GRnext=\!GRccheckforAP
  \else
    \let\!GRnext=\!GRdcheckforAP
  \fi
  \!GRnext}
\def\!GRccheckforAP{%
  \!ifnextchar{p}{\!cgetanchor }
    {\!cgetanchor point at {\!xshadesave} {\!yshadesave} }}
\def\!cgetanchor point at #1 #2 {%
  \edef\!xshadesave{#1}\edef\!yshadesave{#2}%
  \!xshade=\!xshadesave\!xunit  \!yshade=\!yshadesave\!yunit
  \ignorespaces}
\def\!GRdcheckforAP{%
  \!ifnextchar{p}{\!dgetanchor}%
    {\ignorespaces}}
\def\!dgetanchor point at #1 #2 {%
  \!xshade=#1\relax  \!yshade=#2\relax
  \ignorespaces}

\def\setshadesymbol{%
  \!ifnextchar<{\!setshadesymbol}{\!setshadesymbol<,,,> }}

\def\!setshadesymbol <#1,#2,#3,#4> (#5#6){%
  \!setputobject{#5}{#6}%
  \setbox\!shadesymbol=\box\!putobject%
  \!shadesymbolxshift=\!xshift \!shadesymbolyshift=\!yshift
  \!dimenA=\!xshift \advance\!dimenA \!smidge
  \!override\!dimenA{#1}\!lshrinkage%
  \!dimenA=\!wd \advance \!dimenA -\!xshift
    \advance\!dimenA \!smidge
    \!override\!dimenA{#2}\!rshrinkage
  \!dimenA=\!dp \advance \!dimenA \!yshift
    \advance\!dimenA \!smidge
    \!override\!dimenA{#3}\!bshrinkage
  \!dimenA=\!ht \advance \!dimenA -\!yshift
    \advance\!dimenA \!smidge
    \!override\!dimenA{#4}\!tshrinkage
  \ignorespaces}
\def\!smidge{-.2pt}%

\def\!override#1#2#3{%
  \edef\!!override{#2}%
  \ifx \!!override\empty
    #3=#1\relax
  \else
    \if z\!!override
      #3=\!zpt
    \else
      \ifx \!!override\!blankz
        #3=\!zpt
      \else
        #3=#2\relax
      \fi
    \fi
  \fi
  \ignorespaces}
\def\!blankz{ z}

\setshadesymbol ({\fiverm .})

\def\!startvshade#1(#2,#3,#4){%
  \let\!!xunit=\!xunit%
  \let\!!yunit=\!yunit%
  \let\!!xshade=\!xshade%
  \let\!!yshade=\!yshade%
  \def\!getshrinkages{\!vgetshrinkages}%
  \let\!setshadelocation=\!vsetshadelocation%
  \!xS=\!M{#2}\!!xunit
  \!ybS=\!M{#3}\!!yunit
  \!ytS=\!M{#4}\!!yunit
  \!shadexorigin=\!xorigin  \advance \!shadexorigin \!shadesymbolxshift
  \!shadeyorigin=\!yorigin  \advance \!shadeyorigin \!shadesymbolyshift
  \ignorespaces}

\def\!starthshade#1(#2,#3,#4){%
  \let\!!xunit=\!yunit%
  \let\!!yunit=\!xunit%
  \let\!!xshade=\!yshade%
  \let\!!yshade=\!xshade%
  \def\!getshrinkages{\!hgetshrinkages}%
  \let\!setshadelocation=\!hsetshadelocation%
  \!xS=\!M{#2}\!!xunit
  \!ybS=\!M{#3}\!!yunit
  \!ytS=\!M{#4}\!!yunit
  \!shadexorigin=\!xorigin  \advance \!shadexorigin \!shadesymbolxshift
  \!shadeyorigin=\!yorigin  \advance \!shadeyorigin \!shadesymbolyshift
  \ignorespaces}

\def\!lattice#1#2#3#4#5{%
  \!dimenA=#1
  \!dimenB=#2
  \!countB=\!dimenB
  \!dimenC=#3
  \advance\!dimenC -\!dimenA
  \!countA=\!dimenC
  \divide\!countA \!countB
  \ifdim\!dimenC>\!zpt
    \!dimenD=\!countA\!dimenB
    \ifdim\!dimenD<\!dimenC
      \advance\!countA 1 
    \fi
  \fi
  \!dimenC=\!countA\!dimenB
    \advance\!dimenC \!dimenA
  #4=\!countA
  #5=\!dimenC
  \ignorespaces}

\def\!qshade#1(#2,#3,#4)#5(#6,#7,#8){%
  \!xM=\!M{#2}\!!xunit
  \!ybM=\!M{#3}\!!yunit
  \!ytM=\!M{#4}\!!yunit
  \!xE=\!M{#6}\!!xunit
  \!ybE=\!M{#7}\!!yunit
  \!ytE=\!M{#8}\!!yunit
  \!getcoeffs\!xS\!ybS\!xM\!ybM\!xE\!ybE\!ybB\!ybC
  \!getcoeffs\!xS\!ytS\!xM\!ytM\!xE\!ytE\!ytB\!ytC
  \def\!getylimits{\!qgetylimits}%
  \!shade{#1}\ignorespaces}

\def\!lshade#1(#2,#3,#4){%
  \!xE=\!M{#2}\!!xunit
  \!ybE=\!M{#3}\!!yunit
  \!ytE=\!M{#4}\!!yunit
  \!dimenE=\!xE  \advance \!dimenE -\!xS
  \!dimenC=\!ytE \advance \!dimenC -\!ytS
  \!divide\!dimenC\!dimenE\!ytB
  \!dimenC=\!ybE \advance \!dimenC -\!ybS
  \!divide\!dimenC\!dimenE\!ybB
  \def\!getylimits{\!lgetylimits}%
  \!shade{#1}\ignorespaces}

\def\!getcoeffs#1#2#3#4#5#6#7#8{%
  \!dimenC=#4\advance \!dimenC -#2
  \!dimenE=#3\advance \!dimenE -#1
  \!divide\!dimenC\!dimenE\!dimenF
  \!dimenC=#6\advance \!dimenC -#4
  \!dimenH=#5\advance \!dimenH -#3
  \!divide\!dimenC\!dimenH\!dimenG
  \advance\!dimenG -\!dimenF
  \advance \!dimenH \!dimenE
  \!divide\!dimenG\!dimenH#8
  \!removept#8\!t
  #7=-\!t\!dimenE
  \advance #7\!dimenF
  \ignorespaces}

\def\!shade#1{%
  \!getshrinkages#1<,,,>\!nil
  \advance \!dimenE \!xS
  \!lattice\!!xshade\!dshade\!dimenE
    \!parity\!xpos
  \!dimenF=-\!dimenF
    \advance\!dimenF \!xE
  \!loop\!not{\ifdim\!xpos>\!dimenF}
    \!shadecolumn%
    \advance\!xpos \!dshade
    \advance\!parity 1
  \repeat
  \!xS=\!xE
  \!ybS=\!ybE
  \!ytS=\!ytE
  \ignorespaces}

\def\!vgetshrinkages#1<#2,#3,#4,#5>#6\!nil{%
  \!override\!lshrinkage{#2}\!dimenE
  \!override\!rshrinkage{#3}\!dimenF
  \!override\!bshrinkage{#4}\!dimenG
  \!override\!tshrinkage{#5}\!dimenH
  \ignorespaces}
\def\!hgetshrinkages#1<#2,#3,#4,#5>#6\!nil{%
  \!override\!lshrinkage{#2}\!dimenG
  \!override\!rshrinkage{#3}\!dimenH
  \!override\!bshrinkage{#4}\!dimenE
  \!override\!tshrinkage{#5}\!dimenF
  \ignorespaces}

\def\!shadecolumn{%
  \!dxpos=\!xpos
  \advance\!dxpos -\!xS
  \!removept\!dxpos\!dx
  \!getylimits
  \advance\!ytpos -\!dimenH
  \advance\!ybpos \!dimenG
  \!yloc=\!!yshade
  \ifodd\!parity
     \advance\!yloc \!dshade
  \fi
  \!lattice\!yloc{2\!dshade}\!ybpos%
    \!countA\!ypos
  \!dimenA=-\!shadexorigin \advance \!dimenA \!xpos
  \loop\!not{\ifdim\!ypos>\!ytpos}
    \!setshadelocation
    \!rotateaboutpivot\!xloc\!yloc%
    \!dimenA=-\!shadexorigin \advance \!dimenA \!xloc
    \!dimenB=-\!shadeyorigin \advance \!dimenB \!yloc
    \kern\!dimenA \raise\!dimenB\copy\!shadesymbol \kern-\!dimenA
    \advance\!ypos 2\!dshade
  \repeat
  \ignorespaces}

\def\!qgetylimits{%
  \!dimenA=\!dx\!ytC
  \advance\!dimenA \!ytB
  \!ytpos=\!dx\!dimenA
  \advance\!ytpos \!ytS
  \!dimenA=\!dx\!ybC
  \advance\!dimenA \!ybB
  \!ybpos=\!dx\!dimenA
  \advance\!ybpos \!ybS}

\def\!lgetylimits{%
  \!ytpos=\!dx\!ytB
  \advance\!ytpos \!ytS
  \!ybpos=\!dx\!ybB
  \advance\!ybpos \!ybS}

\def\!vsetshadelocation{
  \!xloc=\!xpos
  \!yloc=\!ypos}
\def\!hsetshadelocation{
  \!xloc=\!ypos
  \!yloc=\!xpos}

\def\!axisticks {%
  \def\!nextkeyword##1 {%
    \expandafter\ifx\csname !ticks##1\endcsname \relax
      \def\!next{\!fixkeyword{##1}}%
    \else
      \def\!next{\csname !ticks##1\endcsname}%
    \fi
    \!next}%
  \!axissetup
    \def\!axissetup{\relax}%
  \edef\!ticksinoutsign{\!ticksinoutSign}%
  \!ticklength=\longticklength
  \!tickwidth=\linethickness
  \!gridlinestatus
  \!setticktransform
  \!maketick
  \!tickcase=0
  \def\!LTlist{}%
  \!nextkeyword}

\def\ticksout{%
  \def\!ticksinoutSign{+}}

\ticksout

\def\nogridlines{%
  \def\!gridlinestatus{\!gridlinestoofalse}}
\nogridlines

\def\loggedticks{%
  \def\!setticktransform{\let\!ticktransform=\!logten}}
\def\unloggedticks{%
  \def\!setticktransform{\let\!ticktransform=\!donothing}}
\def\!donothing#1#2{\def#2{#1}}
\unloggedticks

\expandafter\def\csname !ticks/\endcsname{%
  \!not {\ifx \!LTlist\empty}
    \!placetickvalues
  \fi
  \def\!tickvalueslist{}%
  \def\!LTlist{}%
  \expandafter\csname !axis/\endcsname}

\def\!maketick{%
  \setbox\!boxA=\hbox{%
    \beginpicture
      \!setdimenmode
      \setcoordinatesystem point at {\!zpt} {\!zpt}
      \linethickness=\!tickwidth
      \ifdim\!ticklength>\!zpt
        \putrule from {\!zpt} {\!zpt} to
          {\!ticksinoutsign\!tickxsign\!ticklength}
          {\!ticksinoutsign\!tickysign\!ticklength}
      \fi
      \if!gridlinestoo
        \putrule from {\!zpt} {\!zpt} to
          {-\!tickxsign\!xaxislength} {-\!tickysign\!yaxislength}
      \fi
    \endpicturesave <\!Xsave,\!Ysave>}%
    \wd\!boxA=\!zpt}

\def\!ticksin{%
  \def\!ticksinoutsign{-}%
  \!maketick
  \!nextkeyword}

\def\!ticksout{%
  \def\!ticksinoutsign{+}%
  \!maketick
  \!nextkeyword}

\def\!tickslength<#1> {%
  \!ticklength=#1\relax
  \!maketick
  \!nextkeyword}

\def\!tickslong{%
  \!tickslength<\longticklength> }

\def\!ticksshort{%
  \!tickslength<\shortticklength> }

\def\!tickswidth<#1> {%
  \!tickwidth=#1\relax
  \!maketick
  \!nextkeyword}

\def\!ticksandacross{%
  \!gridlinestootrue
  \!maketick
  \!nextkeyword}

\def\!ticksbutnotacross{%
  \!gridlinestoofalse
  \!maketick
  \!nextkeyword}

\def\!tickslogged{%
  \let\!ticktransform=\!logten
  \!nextkeyword}

\def\!ticksunlogged{%
  \let\!ticktransform=\!donothing
  \!nextkeyword}

\def\!ticksunlabeled{%
  \!tickcase=0
  \!nextkeyword}

\def\!ticksnumbered{%
  \!tickcase=1
  \!nextkeyword}

\def\!tickswithvalues#1/ {%
  \edef\!tickvalueslist{#1! /}%
  \!tickcase=2
  \!nextkeyword}

\def\!ticksquantity#1 {%
  \ifnum #1>1
    \!updatetickoffset
    \!countA=#1\relax
    \advance \!countA -1
    \!ticklocationincr=\!axisLength
      \divide \!ticklocationincr \!countA
    \!ticklocation=\!axisstart
    \loop \!not{\ifdim \!ticklocation>\!axisend}
      \!placetick\!ticklocation
      \ifcase\!tickcase
          \relax 
        \or
          \relax 
        \or
          \expandafter\!gettickvaluefrom\!tickvalueslist
          \edef\!tickfield{{\the\!ticklocation}{\!value}}%
          \expandafter\!listaddon\expandafter{\!tickfield}\!LTlist%
      \fi
      \advance \!ticklocation \!ticklocationincr
    \repeat
  \fi
  \!nextkeyword}

\def\!ticksat#1 {%
  \!updatetickoffset
  \edef\!Loc{#1}%
  \if /\!Loc
    \def\next{\!nextkeyword}%
  \else
    \!ticksincommon
    \def\next{\!ticksat}%
  \fi
  \next}

\def\!ticksfrom#1 to #2 by #3 {%
  \!updatetickoffset
  \edef\!arg{#3}%
  \expandafter\!separate\!arg\!nil
  \!scalefactor=1
  \expandafter\!countfigures\!arg/
  \edef\!arg{#1}%
  \!scaleup\!arg by\!scalefactor to\!countE
  \edef\!arg{#2}%
  \!scaleup\!arg by\!scalefactor to\!countF
  \edef\!arg{#3}%
  \!scaleup\!arg by\!scalefactor to\!countG
  \loop \!not{\ifnum\!countE>\!countF}
    \ifnum\!scalefactor=1
      \edef\!Loc{\the\!countE}%
    \else
      \!scaledown\!countE by\!scalefactor to\!Loc
    \fi
    \!ticksincommon
    \advance \!countE \!countG
  \repeat
  \!nextkeyword}

\def\!updatetickoffset{%
  \!dimenA=\!ticksinoutsign\!ticklength
  \ifdim \!dimenA>\!offset
    \!offset=\!dimenA
  \fi}

\def\!placetick#1{%
  \if!xswitch
    \!xpos=#1\relax
    \!ypos=\!axisylevel
  \else
    \!xpos=\!axisxlevel
    \!ypos=#1\relax
  \fi
  \advance\!xpos \!Xsave
  \advance\!ypos \!Ysave
  \kern\!xpos\raise\!ypos\copy\!boxA\kern-\!xpos
  \ignorespaces}

\def\!gettickvaluefrom#1 #2 /{%
  \edef\!value{#1}%
  \edef\!tickvalueslist{#2 /}%
  \ifx \!tickvalueslist\!endtickvaluelist
    \!tickcase=0
  \fi}
\def\!endtickvaluelist{! /}

\def\!ticksincommon{%
  \!ticktransform\!Loc\!t
  \!ticklocation=\!t\!!unit
  \advance\!ticklocation -\!!origin
  \!placetick\!ticklocation
  \ifcase\!tickcase
    \relax 
  \or 
    \ifdim\!ticklocation<-\!!origin
      \edef\!Loc{$\!Loc$}%
    \fi
    \edef\!tickfield{{\the\!ticklocation}{\!Loc}}%
    \expandafter\!listaddon\expandafter{\!tickfield}\!LTlist%
  \or 
    \expandafter\!gettickvaluefrom\!tickvalueslist
    \edef\!tickfield{{\the\!ticklocation}{\!value}}%
    \expandafter\!listaddon\expandafter{\!tickfield}\!LTlist%
  \fi}

\def\!separate#1\!nil{%
  \!ifnextchar{-}{\!!separate}{\!!!separate}#1\!nil}
\def\!!separate-#1\!nil{%
  \def\!sign{-}%
  \!!!!separate#1..\!nil}
\def\!!!separate#1\!nil{%
  \def\!sign{+}%
  \!!!!separate#1..\!nil}
\def\!!!!separate#1.#2.#3\!nil{%
  \def\!arg{#1}%
  \ifx\!arg\!empty
    \!countA=0
  \else
    \!countA=\!arg
  \fi
  \def\!arg{#2}%
  \ifx\!arg\!empty
    \!countB=0
  \else
    \!countB=\!arg
  \fi}

\def\!countfigures#1{%
  \if #1/%
    \def\!next{\ignorespaces}%
  \else
    \multiply\!scalefactor 10
    \def\!next{\!countfigures}%
  \fi
  \!next}

\def\!scaleup#1by#2to#3{%
  \expandafter\!separate#1\!nil
  \multiply\!countA #2\relax
  \advance\!countA \!countB
  \if -\!sign
    \!countA=-\!countA
  \fi
  #3=\!countA
  \ignorespaces}

\def\!scaledown#1by#2to#3{%
  \!countA=#1\relax
  \ifnum \!countA<0 
    \def\!sign{-}
    \!countA=-\!countA
  \else
    \def\!sign{}%
  \fi
  \!countB=\!countA
  \divide\!countB #2\relax
  \!countC=\!countB
    \multiply\!countC #2\relax
  \advance \!countA -\!countC
  \edef#3{\!sign\the\!countB.}
  \!countC=\!countA 
  \ifnum\!countC=0 
    \!countC=1
  \fi
  \multiply\!countC 10
  \!loop \ifnum #2>\!countC
    \edef#3{#3\!zero}%
    \multiply\!countC 10
  \repeat
  \edef#3{#3\the\!countA}
  \ignorespaces}

\def\!placetickvalues{%
  \advance\!offset \tickstovaluesleading
  \if!xswitch
    \setbox\!boxA=\hbox{%
      \def\\##1##2{%
        \!dimenput {##2} [B] (##1,\!axisylevel)}%
      \beginpicture
        \!LTlist
      \endpicturesave <\!Xsave,\!Ysave>}%
    \!dimenA=\!axisylevel
      \advance\!dimenA -\!Ysave
      \advance\!dimenA \!tickysign\!offset
      \if -\!tickysign
        \advance\!dimenA -\ht\!boxA
      \else
        \advance\!dimenA  \dp\!boxA
      \fi
    \advance\!offset \ht\!boxA
      \advance\!offset \dp\!boxA
    \!dimenput {\box\!boxA} [Bl] <\!Xsave,\!Ysave> (\!zpt,\!dimenA)
  \else
    \setbox\!boxA=\hbox{%
      \def\\##1##2{%
        \!dimenput {##2} [r] (\!axisxlevel,##1)}%
      \beginpicture
        \!LTlist
      \endpicturesave <\!Xsave,\!Ysave>}%
    \!dimenA=\!axisxlevel
      \advance\!dimenA -\!Xsave
      \advance\!dimenA \!tickxsign\!offset
      \if -\!tickxsign
        \advance\!dimenA -\wd\!boxA
      \fi
    \advance\!offset \wd\!boxA
    \!dimenput {\box\!boxA} [Bl] <\!Xsave,\!Ysave> (\!dimenA,\!zpt)
  \fi}

\normalgraphs
\catcode`!=12 

\catcode`@=11 \catcode`!=11

\let\!pictexendpicture=\endpicture
\let\!pictexframe=\frame
\let\!pictexlinethickness=\linethickness
\let\!pictexmultiput=\multiput
\let\!pictexput=\put

\def\beginpicture{%
  \setbox\!picbox=\hbox\bgroup%
  \let\endpicture=\!pictexendpicture
  \let\frame=\!pictexframe
  \let\linethickness=\!pictexlinethickness
  \let\multiput=\!pictexmultiput
  \let\put=\!pictexput
  \let\input=\@@input   
  \!xleft=\maxdimen
  \!xright=-\maxdimen
  \!ybot=\maxdimen
  \!ytop=-\maxdimen}

\let\frame=\!latexframe

\let\pictexframe=\!pictexframe

\let\linethickness=\!latexlinethickness
\let\pictexlinethickness=\!pictexlinethickness

\let\\=\@normalcr
\catcode`@=12 \catcode`!=12

\begin{document}

\title{Scalar-Flat K\"{a}hler Surfaces of All Genera}
\author{\parbox{3in}{\center Jongsu Kim\\
ICTP, Trieste}\vspace{.115in}
\\
\parbox{3in}{\center Claude LeBrun\thanks{Supported in part
by NSF grant DMS 92-04093.}\\SUNY, Stony
 Brook}
\vspace{.155in}
\\
 and\\
 \parbox{3in}{\center  Massimiliano Pontecorvo\\
Universit\`a di Napoli
}}
\maketitle
\begin{abstract}
Let $(M,J)$ be a compact complex 2-manifold which
 which admits a K\"ahler metric for which
the integral of the scalar curvature is
non-negative. Also
suppose that $M$ does { not} admit a
Ricci-flat K\"ahler metric. Then
 if $M$ is blown up at sufficiently many points, the
resulting complex surface $(\tilde{M}, \tilde{J})$   admits
K\"ahler metrics with scalar curvature identically equal to zero.
This proves Conjecture 1 of \cite{LS1}.
 \end{abstract}

\vfill    \pagebreak

\section{Introduction}
 The problem of determining which compact
complex manifolds  $(M,J)$ admit
K\"ahler metrics $g$ with constant scalar curvature was first
formulated and studied by
Calabi \cite{cal} in the late 1950's.
Posed in this generality, Calabi's
question is one to which the answer
still eludes us; but there
are two unrelated fronts
on which notable progress  has  been made.

The  most dramatic   progress has
 been made in relation to the case  in which one
also requires that the K\"ahler class
be a multiple of the manifold's first
Chern class--- that is, in relation to
the question of
  which compact
complex manifolds  admit
  K\"ahler-{\em Einstein} metrics.
   This may be reformulated as an
existence problem for
 solutions of the
complex Monge-Amp\`ere equation,
and, working from this point of view,
 Aubin \cite{aub}, Yau \cite{yau}, and others
\cite{ty,siu,nadel,tian}
have given us a fairly complete solution of the problem.

Progress, albeit of a more modest kind,
has also  been made   regarding the
case in which the  scalar curvature
of $g$ is required
to  equal
zero and $\mbox{dim}_{\Bbb C}M=2$;
solutions $(M,J,g)$ of this problem are the
 {\em scalar-flat K\"ahler surfaces} of the title.
The essential reason why this case is
more tractable than others is that
the underlying oriented Riemannian manifold
$(M,g)$ of any scalar-flat K\"ahler surface is
automatically \cite{gau} {\em anti-self-dual},
allowing one to invoke the Penrose twistor
correspondence \cite{P,AHS} and giving rise to phenomena
familiar from the theory of
totally integrable systems. In particular,
all scalar-flat K\"ahler surfaces with semi-free
Killing fields
 can be written down
 explicitly \cite{L,Ltan}, and deforming these \cite{LS1} leads
to a reasonably
complete picture  when the fundamental
group is large.

A fundamental limitation of the above approach is
 that any
simply-connected
scalar-flat K\"ahler surface  with a semi-free Killing field
automatically must have large $\pi_1$, and solutions with
small fundamental group are thus {\em a priori}
 inaccessible by this method.
The present article, however,  will
 finally prove the existence of
 (non-Ricci-flat) scalar-flat K\"ahler surfaces
with $\pi_1=0$ and ${\Bbb Z}\oplus {\Bbb Z}$.
Our key trick is  a Kummer-type
construction which allows us to produce
 new solutions by
smoothing the orbifold singularities of  ${\Bbb Z}_2$-quotients
of old solutions. Successfully carrying this out  involves
melding two previous
extensions \cite{KP,LS2} of the
Donaldson-Friedman   method \cite{DF} of constructing
anti-self-dual metrics on connected sums. After
comparing the solutions
obtained in this way with the biholomorphism types allowed by
surface classification \cite{yau0}, we deduce the following:

\begin{main}\label{maj}
 Let $(M,J)$ be a compact complex 2-manifold which admits
a K\"ahler metric for which the integral of the
scalar curvature is non-negative. Then precisely
one of the following holds:
\begin{itemize}
\item $(M,J)$ admits a Ricci-flat K\"ahler metric; or
\item any blow-up of  $(M,J)$ has blow-ups $(\tilde{M},\tilde{J})$
which admit
scalar-flat K\"ahler metrics.
\end{itemize}
\end{main}
This proves Conjecture 1 of \cite{LS1}.

A corollary of this is the following:

\begin{main}\label{min}
 Let $(M,J)$ be a compact complex 2-manifold which admits
a K\"ahler metric for which the integral of the
scalar curvature is positive. Then
 any blow-up of  $(M,J)$ has blow-ups $(\tilde{M},\tilde{J})$
which admit
 K\"ahler metrics of constant positive scalar curvature.
\end{main}

One might  note the formal similarity
 between Theorem \ref{maj} and a recent result of
Taubes \cite{tau}
which asserts
that one can find anti-self-dual
metrics on  the connected sum of any smooth oriented 4-manifold $M$
 with enough copies of $\overline{\Bbb CP}_2$. However,  Taubes' proof
is direct,
whereas  ours falls back on classification theory.
One would hope that different proof
of Theorem \ref{maj}, proceeding  along Taubes' lines, might
shed more light on Calabi's  general problem.

\section{The Quotient Theorem} \label{quoth}

In this section, we state the central technical result of this
article, and then set up the framework in which it
will be proved. This result, from which our main results
will be deduced,   is the following:

\begin{thm}[Quotient Theorem] \label{key}
Let $(N,J_N,g_N)$ be a non-minimal
compact complex surface with scalar-flat K\"ahler metric,
and let $\Phi : N\to N$, $\Phi^2={\bf 1}$,
 be a holomorphic isometry
with only isolated fixed points.
Let $(M, J_M)$ be obtained from $N/\Phi$ by replacing each
singular point with
a ${\Bbb CP}_1$ of self-intersection $-2$.
Then there exist scalar-flat K\"ahler metrics $g_M$ on
$(M, J_M)$.
\end{thm}

Here a compact complex surface $N$ is called {\em non-minimal}
if is obtained from  another surface by blowing up;
this is equivalent to saying
that $N$ contains a ${\Bbb CP}_1$ of self-intersection $-1$.
Our description of $M$  amounts to saying that
if $N$ is blown up at the fixed points of $\Phi$, the resulting
complex surface $\tilde{N}$ is a branched double cover of
$M$, with the newly-introduced exceptional divisors as ramification locus.
More abstractly, $M$ is the minimal resolution of the
ordinary double-point singularities
of the variety  $N/\Phi$.

We will prove this theorem by using the theory of
twistor spaces \cite{AHS,P}. For our
purposes, a {\em twistor space}  means a compact complex
3-manifold $Z$ equipped with a free anti-holomorphic
involution $\sigma : Z\to Z$ and a foliation by
$\sigma$-invariant
rational curves ${\Bbb CP}_1\subset Z$ with normal
bundle ${\cal  O}(1)\oplus
{\cal O}(1)$.
Let  $X$ denote the
 leaf space of this foliation
by the so-called {\em real twistor lines}, and let  $\wp:\tilde{Z}\to X$
denote the quotient map. There is then a canonical
anti-self-dual conformal
metric $[g]$ on $X$, characterized by
the requirement that the image of every holomorphic tangent
space  $T^{1,0}_z\tilde{Z}$ should be a
$g$-isotropic  subspace of ${\Bbb C}\otimes TX$.
Conversely, every anti-self-dual manifold   arises in this
way, and does so in an essentially unique manner.

If a twistor space $Z$ contains a compact complex surface $D$ which
is disjoint from its conjugate $\bar{D}:=\sigma (D)$ and
has homological intersection number 1 with a twistor line,
then   $\wp|_D:D\to X$
is a diffeomorphism, and
 $[g]$ pulls back from $X$ to yield a conformal class of
Hermitian metrics  on $D$. Because $D$ is compact, the anti-self-duality
of $[g]$ implies \cite{boyer,pont} that this conformal class is locally
represented on $D$ by scalar-flat K\"ahler metrics. If, moreover,
$b_1(D)$ is even, there is a  globally-defined
scalar-flat K\"ahler metric $g\in [g]$, and this
global representative is uniquely determined once its
total volume is specified; conversely, every scalar-flat
K\"ahler surface arises from this construction, and does
so in an essentially unique manner.  Thus,
in order to prove Theorem \ref{key}, it
 suffices to produce a twistor space
$Z$   containing a copy of
the complex surface $(M,J_M)$ which is disjoint
from its conjugate and intersects some twistor line
transversely in  one point.
We will do just this by  refining the methods of
\cite{LS2}, where   anti-self-dual metrics were constructed on
the underlying smooth manifold $M$.

To begin this construction, let $Z_N$ denote the twistor
space of $(N,g_N)$, and let $L_1, \ldots , L_k$ be the
twistor lines of the fixed points
of $\Phi$. Let $\tilde{Z}_N$ be the blow-up of
$Z_N$ along these lines, and let $Q_1, \ldots,
Q_k$ be the
exceptional divisors in $\tilde{Z}_N$
 corresponding  to $L_1, \ldots , L_k$; thus
$Q_j\cong {\Bbb CP}_1\times {\Bbb CP}_1$, $j=1, \ldots k$,
and each of these 2-quadrics has normal bundle ${\cal O}(1,-1)\to
{\Bbb CP}_1\times {\Bbb CP}_1$.
Since the derivative of $\Phi$ at its $k$ isolated
fixed points must be $-1$, the induced biholomorphism
$\hat{\Phi}:\tilde{Z}_N\to \tilde{Z}_N$ fixes each
$Q_j$ and
acts on its normal bundle by $-1$.
The quotient $Z_-:=\tilde{Z}_N/\hat{\Phi}$
can thus be given the structure of a compact complex
manifold in  a unique way that the quotient map
 $\tilde{Z}_N\to Z_-$ becomes a branched
covering, with $Q:=\bigcup Q_j$ as ramification locus.
Let  $Q_{j-}$ denote the image of this $Q_j$ in $Z_-$,
which is an imbedded  quadric with normal bundle ${\cal O}(2,-2)$,
and let $Q_-=\bigcup Q_{j-}$.

The twistor space $Z_N$ contains a hypersurface $D_N$
corresponding to the complex structure $J_N$, as well as
a disjoint hypersurface   $\bar{D}_N:=\sigma (D_N)$ corresponding
to the conjugate complex structure $-J_N$; indeed,
$D_N$ and $\bar{D}_N$ are  respectively isomorphic to
$(N,\pm J_N)$ as complex surfaces. As the action of $\hat{\Phi}$
sends each such surface to itself,
there are disjoint hypersurfaces $D_-$ and $\bar{D}_-$ in $Z_-$
obtained by first taking the  proper transforms
in $\tilde{Z}_N$
of $D_N$ and $\bar{D}_N$
and then projecting these hypersurfaces to
$Z_-$. Notice that $D_-$ is exactly a copy of
$(M,J_M)$, whereas  $\bar{D}_-$ is
a copy of  $(M,-J_M)$.
Set $\ell_{j-}:=D_-\cap Q_j$, $\bar{\ell}_{j-}:=\bar{D}_-\cap Q_j$,
$\ell_-:=\bigcup \ell_{j-}$, and
$\bar{\ell}_-:=\bigcup\bar{\ell}_{j-}$.

Our next step is to let $Z_+$ consist of $k$ disjoint copies of the
complex 3-fold
$\tilde{Z}_{EH}$ obtained from the orbifold
twistor space of the conformally compactified Eguchi-Hanson
metric by blowing up the twistor line of infinity.
To describe    $\tilde{Z}_{EH}$ explicitly \cite{hit,LS2},
start with the
${\Bbb CP}_3$-bundle  $\pi: {\cal B}\to {\Bbb CP}_1$ defined by
$${\cal B}={\Bbb P}({\cal O}(2)^{\oplus  3}\oplus {\cal O}) ~ ,$$
where our conventions are that ${\Bbb P}(E):= (E-0)/{\Bbb C}^{\times}$.
Let ${\cal O}(1,0):=\pi^*{\cal O}(1)$, and let ${\cal O}(0,-1)$ be the
universal bundle, whose principal ${\Bbb C}^{\times}$-bundle is
$[({\cal O}(2)^{\oplus  3}\oplus {\cal O})-0]\to {\cal B}$. The
``homogeneous coordinates'' of  ${\cal O}(2)^{\oplus  3}\oplus {\cal O}$
are canonical sections $x,y,z\in \Gamma {\cal O} (2,1)$ and
$t\in \Gamma{\cal O} (0,1)$. Let $a\in \Gamma ({\Bbb CP}_1, {\cal O}(2))$
be a non-trivial section which is invariant under the anti-holomorphic
involution $\conj$ of ${\cal O}(2)=T{\Bbb CP}_1$ induced by the
antipodal map of $S^2={\Bbb CP}_1$; and let the 2 distinct
zeroes of $a$ be  called  the north and south poles. Our
blown-up twistor space is then obtained from the hypersurface
$$xy=z^2-t^2a^2$$
in ${\cal B}$ by replacing the
  two singular points $x=y=z=0$ with ${\Bbb CP}_1$'s.
 The quadric $t=0$
corresponds to the blow-up of the
orbifold twistor line ``at infinity,'' whereas
 the real structure is given by
$$[x:y:z:t]\to [\conj (y):\conj (x):\conj (z): \bar{t} ~ ] .$$
 If we let $D_{EH}\subset \tilde{Z}_{EH}$ denote
the Hirzebruch surface
over the south pole and let $\bar{D}_{EH}$
denote the Hirzebruch surface over the north pole,
then we may define
$D_+\subset Z_+$ and $\bar{D}_+\subset Z_+$ to
consist of
$k$ disjoint copies of $D_{EH}$ and  $\bar{D}_{EH}$, respectively.
Let us use $Q_{j+}$ to denote the appropriate copy of the
$t=0$ quadric in $\tilde{Z}_{EH}$ and $Q_+$ to denote
$\bigcup Q_{j+}$. Set $\ell_{j+}=Q_{j+}\cap D_+$,
$\bar{\ell}_{j+}=Q_{j+}\cap \bar{D}_+$, $\ell_+=\bigcup \ell_{j+}$,
and $\bar{\ell}_+=\bigcup \bar{\ell}_{j+}$.

Now let $Z_0=Z_-\cup_{Q_j} Z_+$ be obtained from
the disjoint union $Z_-\sqcup Z_+$ by biholomorphically
identifying  $Q_{j-}\subset Z_-$ with $Q_{j+}\subset Z_+$
in such a way that the real structures agree
and such that  $\ell_{j-}$ is identified with $\ell_{j+}$.
 (This actually
specifies the gluing procedure uniquely, modulo real automorphisms
of $Z_+$.)
Set $D_0=D_-\cup D_+$ and $\bar{D}_0=\bar{D}_-\cup \bar{D}_+$.
The result is that  $Z_0$, $D_0$ and  $\bar{D}_0 $ are
complex spaces with normal crossing singularities,
\label{sing} and there is an induced real structure
$\sigma: Z_0\to Z_0$ which interchanges $D_0$ and  $\bar{D}_0$.
By a slight abuse of notation, we
will denote the image of $Q_{j\pm}$ in $Z_0$ by $Q_j$,
whereas the images of $\ell_{j\pm}$, $\bar{\ell}_{j\pm}$,
$\ell_{\pm}$, $\bar{\ell}_{\pm}$, and $Q_{\pm}$ will
respectively be denoted by $\ell_{j}$, $\bar{\ell}_{j}$,
$\ell$, $\bar{\ell}$, and $Q$.

Our method of proving Theorem \ref{key}  will now
be as
follows: {\em we will construct  twistor spaces
$Z_t$ containing hypersurfaces $D_t\cong (M,J_M)$  by
simultaneously smoothing the singularities of
$Z_0$ and $D_0\subset Z_0$.}
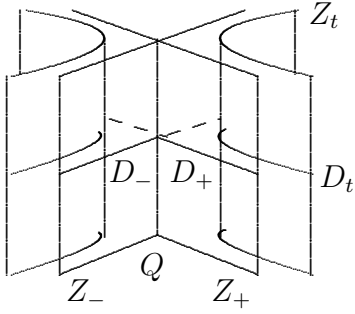
\begin{center}
\mbox{
\beginpicture
\setplotarea x from 0 to 200, y from -20 to 100
\put {$Z_+$} [B1] at 125  -10
\put {$Z_-$} [B1] at  70  -10
\put {$Z_t$} [B1] at 160 95
\put {$Q$} [B1] at    95 0
\put {$D_+$} [B1] at  110 35
\put {$D_-$} [B1] at   87 35
\put {$D_t$} [B1] at   165 33
{\setlinear
\plot 65 100 135 75 /
\plot 97 52 135 38 /
\plot 60 75  130 100 /
\plot 60 38  97 52 /
\plot 135 75 135 0 /
\plot 60 75 60 0 /
\plot 60 0 97 15 /
\plot 97 15 135 0 /
\plot 97 89 97 15 /
\plot 121 89 121 14 /
\plot 155 75 155 0 /
\plot 149 100 149 76 /
\plot 40 75 40 0 /
\plot 77 89 77 14 /
\plot 45 100 45 76 /
}
{\setlinear
\setdashes
\plot 75 59 97 52 /
\plot 97 53  117 59 /
}
{\setquadratic
\plot  148  75   113  89   143  100 /
\plot  148  0    117 11        115  17 /
\plot 34 75 69 89 39 100 /
\plot 34 38     62  48         67 54  /
\plot  148  38    117 48        115  55 /
\plot 34 0     62 10         67 16  /
}
{\setquadratic
\setdashes
}
\endpicture
}
\end{center}
 Because the  twistor spaces
$Z_t$ must admit real structures, we will of course also need
to smooth the singularities of $\bar{D}_0$, too,
and it is natural to also include this stipulation from the
outset.
In \S \ref{sta}, Theorem \ref{key} will now be proved
constructing  relative smoothings of precisely this type.

\section{Proof of the Quotient Theorem}\label{sta}

Continuing the discussion of  \S \ref{quoth},
let $\ddbar\subset Z_0$ denote the disjoint union
of $D_0$ and $\bar{D}_0=\sigma(D_0)$, and let
$f:\ddbar\hookrightarrow Z_0$ be the
tautological
 holomorphic imbedding. In this section we shall prove
 Theorem  \ref{key} by  studying the
deformation theory of the pair of the pair $(Z_0, \ddbar)$.
Our approach  depends crucially  upon the results of
Ran \cite{R1,R2}.

Let us first warm up
by discussing deformations of the singular surface
$D_0$, noting all along that such a discussion
 will, by conjugation, automatically
also implicitly completely
describe the deformation theory of $\bar{D}_0$, and hence that of
the disjoint union $\ddbar=D_0\sqcup \bar{D}_0$, too.
Now since
$D_0=D_-\cup_l D_+$ is obtained from $D_-\sqcup D_+$ by identifying
$(-2)$-curves $\ell_{j-}\subset D_-$ with  $(+2)$-curves
 $\ell_{j+}\subset D_+$,
$D_0$ is a singular complex surface with normal crossing
singularities along
$\ell\cong \ell_-\cong \ell_+$, and  satisfies the
so-called {\it d-semistable condition}---  the
two normal bundles $\nu_{\ell_{\pm}, D_{\pm}}$ of
the singular hypersurface
are dual to each other.

Because $D_+$ has $k$ connected components, each of which
is isomorphic to the second
Hirzebruch surface ${\Bbb P} (\CO\oplus\CO(2))$, and
as each  connected component of $\ell_+$ corresponds to the
zero section of $\CO(2)\subset {\Bbb P} (\CO\oplus\CO(2))$,
it follows that
$$
h^j(\Theta_{D_+}\otimes\i_{\ell_+})=0~~\forall j,$$
$$
h^j(\Theta_{D_+})= \left\{
\begin{array}{ll}6k&j=0\\0&j\neq 0,\end{array}\right.  $$
and
$$
h^j(\Theta_{D_+,\ell_+})=h^j(\Theta_{\ell_+}) =
 h^j(\nu_{\ell_+, D_+})=\left\{
\begin{array}{ll}3k&j=0\\0&j\neq 0\end{array}\right. ,
$$
where $\Theta_{Y,X}$ denotes the sheaf of
holomorphic vector fields on $Y$ which represent $0$ in the normal bundle
$\nu_{X,Y}$ of the complex submanifold $X\subset Y$. Using these facts,
we may now prove the following:

\begin{lem} Let  $\tau^0_{D_0}$
denote the sheaf of derivations of $\CO_{D_0}$.
Then
$$
H^j(\tau^0_{D_0}) \cong   \left\{ \begin{array}{ll}
H^0(\Theta_{D_-})     & j=0 \\
H^1(\Theta_{D_-,\ell_-}) & j=1 \\
0                     & j\geq 2. \end{array}  \right.
$$
\end{lem}

\begin{proof} Consider the  normalization exact sequence
$$ 0 \longrightarrow \tau^0_\d \longrightarrow
   q_* \Theta_{D_-\sqcup D_+,
\ell_- \sqcup \ell_+} \longrightarrow \iota_* \Theta_\ell
   \longrightarrow  0 ~, $$
where $q:D_-\sqcup D_+\to D_0$ is the quotient map and
$\iota:\ell\to D_0$ is the inclusion.
Since the restriction map
$H^0(\Theta_{D_+,\ell_+})\to H^0(\Theta_\ell)$
is an isomorphism, the associated long exact sequence tells us that
$H^0(\tau^0_{D_0})= H^0(\Theta_{D_-,\ell_-})= H^0(\Theta_{D_-})$,
whereas $H^1(\tau^0_{D_0})= H^1(\Theta_{D_-,\ell_-}).$

Since $H^2(\Theta_{D_+,\ell_+})=0$,  we  also read off  that
$H^2(\tau^0_{D_0}) =H^2(\Theta_{D_-,\ell_-})$,
and it  only remains for us
to show that the latter cohomology group vanishes.
But
$D_-$ is a ruled surface, and the generic
$\bcp_1$  fiber  of $D_-$ is disjoint from
$\ell_-$; indeed,  $N$ is
ruled because \cite{yau0} it is non-minimal and
admits a scalar-flat K\"ahler metric,
 $\tilde{N}$ is obtained from $N$ by blowing up,
and $\tilde{N}$ is a branched cover of $D_-\cong M$,
with branch locus $\ell_-$ the union of the exceptional
divisors introduced by the blow-up $\tilde{N}\to N$.
 Moreover, $\Theta_{D_-,\ell_-}$ is a locally
free sheaf, and Serre duality therefore says that
$H^2(D_-,\Theta_{D_-,\ell_-})=
[H^0(D_-,Hom (\Theta_{D_-,\ell_-},\Omega^2))]^*$.
The restriction of $Hom (\Theta_{D_-,\ell_-},\Omega^2_{D_-})$
to a generic $\bcp_1$  fiber  of $D_-$
is thus isomorphic to  $\CO (-2)\oplus \CO (-4)$,
and any global section of this sheaf therefore vanishes.
Hence $H^2(\tau^0_{D_0}) =H^2(D_-,\Theta_{D_-,\ell_-})=0$, as claimed.
\end{proof}

 This implies that the  deformation theory of
$D_0$ is  unobstructed:

\begin{propn} The complex space $\d$ admits a versal deformation
$$ {\cal D} \stackrel{\varpi}{\longrightarrow}
H^1(D_-,\Theta_{D_-,\ell_-})\times {\Bbb C}^k ~,$$
 with fibers $D_t:=\varpi^{-1}(t_1,t_2)$,
 $t_1\in H^1(D_-,\Theta_{D_-,\ell_-})$,  $t_2\in  \BC^k$, satisfying
\begin{enumerate}
\item $D_t$ is smooth iff
  $t_2\in (\BC^{\times})^k$;
\item $D_t\cong D_-=M$ when $t_1=0$ and $t_2\in (\BC^{\times})^k$; and
\item all  small deformations of $D_-= M$
occur as smooth fibers
 $D_t$ of $\varpi$.
\end{enumerate}
\label{surf}
\end{propn}
\begin{proof} The deformation theory \cite{F,DF}
of $\d$ is governed by the vector spaces
 $T^j_\d={\rm Ext}^j(\Omega^1_\d, \CO_\d)$.
These  may be computed by means of the Ext spectral
sequence \cite{go}
$$E_2^{p,q}=H^p(\d , \tau^j_\d)\Longrightarrow T^{p+q}_\d,$$ where
 $\tau^j_\d=Ext^j(\Omega^1, \CO)$.
Because $D_0$ is a locally complete intersection,
$\tau^j_\d=0$ for $j\geq 2$, and the spectral sequence
therefore degenerates into the exact sequences
$$H^{j-2}(\tau^1_\d)\to H^j(\tau^0_\d)\to T^j_\d \to H^{j-1}(\tau^1_\d)
\to H^{j+1}(\tau^0_\d).$$
Meanwhile, the d-semi-stable condition
 tells us tells us that
$\tau^0_\d\cong \CO_{\ell}$, so
$$H^j(\tau^1_\d)\cong\left\{ \begin{array}{ll}
 \BC ^k&j=0\\0&j\neq 0.\end{array}\right. $$
The lemma now tells us that
 $T^2_\d=0$, and that there is an  exact sequence
$$ 0\to H^1(D_-,\Theta_{D_-,\ell_-})\to T^1_\d \to \BC ^k \to 0~.$$

Since  $T^2_\d=0$, the deformation theory  \cite{F}
of $\d$ is unobstructed, and
there is a versal family $\cal D$ over a neighborhood of $0\in T^1_\d$.
This family has the property that
any 1-dimensional subfamily of this family  smooths
 the normal crossing at $\ell_j$ if
the image of its derivative in
$H^0(\tau^1_\ell)\cong H^0(\ell, \CO)=\BC^k$ is non-zero on $\ell_j$.
Moreover, an effectively parametrized family
with central fiber $\d$
is a subfamily of $\cal D$.

Let us now consider the explicit smoothing $\cal M$ of $\d$ gotten by
blowing up $M\times \BC^k$ along the smooth submanifold
$(\ell_1\times \{a_1=0\})\cup \cdots \cup(\ell_k\times \{a_k=0\}).$
Since $\ell_j\subset M$ has self-intersection $-2$,
the central fiber of ${\cal M} \to \BC^k$ is obtained from
$D_-=M$ by attaching a Hirzebruch surface ${\Bbb P}(\CO (-2)\oplus \CO)$
to $M$ at each $\ell_j$; in other words,  the central fiber is
isomorphic to $\d$.
Since  $\cal M$ is an effectively parametrized family
with central fiber $\d$, it must be contained in $\cal D$ by versality,
and because the fiber of $\cal M$ over any
$t_2\in  (\BC^{\times})^k$ is smooth, the image of
$\BC^k$ in $T^1_\d$ is transverse to the kernel
of $T^1_\d \to H^0(\tau^1_\ell)\cong \BC^k$.
We may now choose new coordinates on $T^1_\d$,
identifying it with
 $H^1(D_-,\Theta_{D_-,\ell_-})\times {\Bbb C}^k$
in such a way that the above explicit family
corresponds to the ${\Bbb C}^k$ subspace,
and so that the $H^1(D_-,\Theta_{D_-,\ell_-})$ subspace
corresponds to the original kernel of the
projection.

If ${\bf M}\to \Delta$ is any deformation of
$M$ over a neighborhood
of $0\in {\Bbb C}$, its blow-up $\tilde{\bf M}$
at the submanifold $\ell_-$ of the central fiber is
a   smoothing of $\d$. By  versality,
 $\tilde{\bf M}$ must therefore be a subfamily of ${\cal D}$,
amd every small deformation of $M$ thus occurs as a smooth
fiber $D_t$.
\end{proof}

\begin{remark} At the tangent space level, the
above proof shows that the exact sequence
$$0\to H^1(\tau^0_\d)\to  T^1_\d  \to H^0(\tau^1_\d)\to 0$$
has a geometrically preferred splitting. This will  later prove
useful.
\end{remark}

\bigskip

Having discussed the deformation theory of
$\ddbar$, the next  step is obviously to discuss that of
$Z_0$; but this, in fact, has already been studied in
\cite{LS2}. As before, the vector spaces $T^q_{Z_0} =
{\rm Ext}^q(\Omega^1_{Z_0}, \CO_{Z_0})$ that control the deformation theory
fit into an exact sequence
$$\begin{array}{cccc}
0 \rightarrow  H^1(Z_0,\tau^0_{Z_0}) \rightarrow T^1_{Z_0} \rightarrow &\!\!
H^0(Q,\CO) \!\! &\rightarrow H^2(Z_0,\tau^0_{Z_0}) \rightarrow T^2_{Z_0}
\rightarrow &\!\! H^1(Q,\CO) \\
& \|&& \| \\
& {\Bbb C}^k&& 0 \\
\end{array}
$$
and the cohomology of
$\tau^0_{Z_0} = {\cal H}om(\Omega^1_{Z_0} ,\CO_{Z_0})$
can in turn be computed via the long exact sequence
$$
H^{j-1}(Q, \Theta_Q)\to
H^j(Z_0,\tau^0_{Z_0})\to H^j(Z_-,\Theta_{Z_-,Q})\oplus
H^j(Z_+,\Theta_{Z_+,Q})\to H^j(Q, \Theta_Q)
$$
Using the explicit form of  $Z_+$,
one then may check that $H^2(Z_+,\Theta_{Z_+,Q})=0$.
On the other hand, $Z_-$ is has a branched cover which is a blow-up
of the twistor space $Z_N$, and this implies that
$H^2(Z_-,\Theta_{Z_-,Q})=[H^2(Z_N,\Theta_{Z_-,Q})]_{\Phi}$,
where the subscript indicates the $+1$-eigenspace
of the automorphism induced by $\Phi$.
One may thus conclude that $T^2_{Z_0}=0$, and that the
smoothing theory of $Z_0$ is therefore unobstructed, once one
knows the following result, which
 was stated in
\cite{LS1}:
\begin{thm}
Suppose that $N$ is a non-minimal compact complex surface
with scalar-flat K\"ahler metric $g_N$. Then its
twistor space satisfies $H^2(Z_N, \Theta)=H^2(Z_N, \Theta\otimes
\kappa^{-1/2})=H^2(Z_N, \Theta_{Z,D\bar{D}})=0$, where
$D$ and $\bar{D}$ are canonical divisors associated with
$J_N$ and $-J_N$.
\end{thm}

\noindent  Unfortunately, while the  proof given in \cite{LS1}
suffices for all  cases needed for our
applications, it  overlooks the case of
surfaces with
 non-semi-free
${\Bbb C}^{\times}$-actions,
and we have therefore chosen to
include a completed proof in the
present article.
As this proof is rather long, however,
 and  involves ideas
quite unrelated to the thrust of the present discussion,
it has  been
 relegated to an appendix (\S \ref{van}).

\medskip

Finally, we
 turn to the deformation theory of the pair $(Z_0, \ddbar)$,
which
is governed by the derived functors $T^j_f$ of Ran \cite{R1} \cite{R2},
in the sense that
$T^1_f$ corresponds to infinitesimal deformations of the imbedding $f$
and obstructions lie in $T^2_f$.
These vector spaces may be computed by means of  a long exact sequence

\begin{eqnarray}
0 & \longrightarrow  & T^0_f \longrightarrow
         T^0_{D_0\bar{D}_0} \oplus T^0_{Z_0}
    \longrightarrow
         \ext^0_f(\Omega^1_{Z_0},\CO_{D_0\bar{D}_0})
\nonumber\\
  & \longrightarrow  & T^1_f
{\longrightarrow}
 T^1_{D_0\bar{D}_0} \oplus T^1_{Z_0}
    {\longrightarrow}
 \ext^1_f(\Omega^1_{Z_0},\CO_{D_0\bar{D}_0})
 \label{middle} \\
  & \longrightarrow  & T^2_f \longrightarrow
         T^2_{D_0\bar{D}_0} \oplus T^2_{Z_0}
    \longrightarrow \cdots \nonumber
\end{eqnarray}
Here  $T^j_\ddbar :=
       \ext^j(\Omega^1_\d,\CO_\d) \oplus \ext^j(\Omega^1_{\dbar},\CO_{\dbar})$
and $T^j_{Z_0} := \ext^j(\Omega^1_{Z_0},\CO_{Z_0})$ are the
usual global Ext groups,
whereas $\ext^j_f(\Omega^1_{Z_0},\CO_\ddbar)$ are the derived functors of
$$\Hom _f(\Omega^1_{Z_0},\CO_\ddbar) :=
 \Hom _{\CO_{\ddbar}}(f^*\Omega^1_{Z_0},\CO_\ddbar) \cong
 \Hom _{\CO_{Z_0}}(\Omega^1_{Z_0},f_*\CO_\ddbar)$$
in either variable.

The local form of the singularities of
the pair $(Z_0, \ddbar)$ is just the same
as in \cite[\S 3]{KP}. Namely, around the singular locus we can take
coordinates $(w_1,\ldots,w_4)\in\BC^4$ so that $Z_0$ is given by
$\{w_1w_2=0\}$ and $\d$ (or $\dbar$) is the hypersurface $\{w_4=0\}.$
Because of this we have a spectral sequence
$$E_2^{p,q}=H^p(\lext^p(\Omega^1_{{Z_0}_{|D_0\bar{D}_0}},\CO_{D_0\bar{D}_0}))
  \Rightarrow \ext^{p+q}_f(\Omega^1_{Z_0},\CO_{D_0\bar{D}_0})$$
where the local ${\cal E}xt$
sheaves can be computed  \cite[3.10]{KP} to be
$$
 \lext^r(\Omega^1_{Z_0}|_{D_0\bar{D}_0}, \CO_{D_0\bar{D}_0})
 \cong  \left\{
 \begin{array}{ll}
  \taur,              & r=0     \\
  \tau^1_{D_0\bar{D}_0} \cong \CO_\llbar , & r=1     \\
  0,                    & r\geq 2.
 \end{array}\right.
$$

We now claim that the relative deformation theory
under consideration is unobstructed as a consequence of Theorem 2.
As in \cite[\S 4]{KP},
we start by considering the  exact sequence

$$ 0 \longrightarrow \tau^0_{Z_0,D_0\bar{D}_0}
   \longrightarrow \tau^0_{D_0\bar{D}_0}\oplus\tau^0_{Z_0} \longrightarrow
   \taurr \longrightarrow 0.   $$

\begin{lem}
Suppose that $N$ satisfies the hypothesis of the Quotient Theorem. Then
$H^2(\tau^0_{Z_0,D_0\bar{D}_0})=H^2(\taurr)=0.$
\end{lem}

\begin{proof} Since we have already observed that
$H^2(\tau^0_{D_0\bar{D}_0})\oplus
H^2(\tau^0_{Z_0})=0$, it is enough to show that
$H^j(\tau^0_{Z_0,D_0\bar{D}_0})$ for $j=2,3.$
To this end, consider the exact normalization  sequence
$$0 \longrightarrow \tau^0_{\z,D_0\bar{D}_0} \longrightarrow
    \Theta_{Z_+,D_+\bar{D}_+Q}
\oplus \Theta_{Z_-,D_-\bar{D}_-Q_-}
\longrightarrow
    \Theta_{Q,\llbar} \longrightarrow 0 ~~ , $$
  where  $\ell=Q_- \cap D_- =Q_- \cap D_+$ and
$\bar{\ell} = Q_- \cap \bar{D}_- =Q \cap \bar{D}_+.$
Because
$$ 0 \longrightarrow \Theta_{Q,\llbar} \longrightarrow \Theta_{Q}
    \longrightarrow \nu_\ell\oplus \nu_{\bar{\ell}} \longrightarrow 0 $$
is exact and
$H^0(\Theta_Q) \rightarrow H^0(\nu_\ell \oplus \nu_{\bar{\ell}}) $
is surjective,
it follows  that
$$ H^j(\tau^0_{\z,D_0\bar{D}_0}) =
   H^j(\Theta_{Z_+,
D_+\bar{D}_+Q})
\oplus H^j(\Theta_{Z_-,
D_-\bar{D}_-Q_-}),\quad j=2,3 .$$
But  rational curves of normal bundle
$\CO (1)\oplus \CO (1)$ sweep out an open subset of $Z_{\pm}$,
so  Serre duality tells us  that
 $H^3(\Theta_{Z_\pm,D_\pm\bar{D}_\pm Q_\pm})=0$.
It thus only remains  to show that
$H^2(\Theta_{Z_\pm,D_\pm\bar{D}_\pm Q_\pm})=0.$

We start by considering $Z_-$, which,  by construction,  fits into
a diagram

\setlength{\unitlength}{1ex}
\begin{center}\begin{picture}(20,17)(0,3)
\put(10,17){\makebox(0,0){$\tilde{Z}_N$}}
\put(2,5){\makebox(0,0){$Z_-$}}
\put(19,5){\makebox(0,0){${Z}_N$}}
\put(15.5,12){\makebox(0,0){$\beta$}}
\put(5,12){\makebox(0,0){$\alpha$}}
\put(11,15.5){\vector(2,-3){6}}
\put(9,15.5){\vector(-2,-3){6}}
\end{picture}\end{center}
in which $\beta$  blows  a disjoint union $Q$ of quadrics
down
into the union $L_{\Phi}$ of $\Phi$-invariant
twistor lines in $Z_N$, and in which
$\alpha$ is a 2-fold branched covering map
with branch locus $Q$. Let  $\cal L$ be the
divisor square-root of  $Q\subset {Z}_N$
which is associated with
this branched cover. Then
\be  \alpha^j_{\ast}
   \Theta_{\tilde Z _{N},\tilde{D}_N\tilde{\bar{D}}_NQ} =
\left\{ \begin{array}{ll} \Theta_{Z_{-},D_- \bar{D}_-Q_-} \oplus
(\Theta_{Z_{-},D_- \bar{D}_-Q_-} \otimes {\cal L}) & j=0 \\
0  & j\neq 0,\end{array}\right. \label{arf}\ee
whereas
\be \label{barf} \beta^j_{\ast}
\Theta_{\tilde{Z}_N,\tilde{D}_{N}\tilde{\bar{D}}_{N}Q} =
\left\{ \begin{array}{ll}  \Theta_{Z_{N},D_{N}\bar{D}_{N}, L_{\Phi}}
&j=0\\0  &j\neq 0.\end{array}\right.\ee
In combination with the short exact sequence
$$ 0 \longrightarrow \Theta_{Z_{N},D_{N}\bar{D}_{N},L_{\Phi}} \longrightarrow
 \Theta_{Z_{N},D_{N}\bar{D}_{N}} \longrightarrow
\nu_{L_{\Phi},Z_{N}} \longrightarrow 0$$
and the observation that
$\nu_{L_{\Phi},Z_{N}}\cong\CO(1) \oplus \CO(1)$ on each
$\bcp_1$ component, (\ref{barf}) tells us that
$$H^2({\ztil_{N}},\Theta_{\tilde{Z_N},
\tilde{D}_{N}\tilde{\bar{D}}_{N}Q})=
H^2(Z_{N},\Theta_{Z_{N},D_{N}\bar{D}_{N},L_{\Phi}})=
H^2(Z_{N},\Theta_{Z_{N},D_{N}\bar{D}_{N}}) =0.$$
But the Leray spectral sequence of (\ref{arf}) says us
that
$$ [ H^2( {\ztil_{N}},
    \Theta_{\tilde{Z_N},\tilde{D}_{N}\tilde{\bar{D}}_{N}Q}  )]_{\Phi}
\cong H^2( Z_{-},\Theta_{Z_{-},D_- \bar{D}_-Q_-}),$$
so that $H^2( Z_{-},\Theta_{Z_{-},D_- \bar{D}_-Q})=0$,
as  claimed.

\smallskip
To finish the proof we have to show
$ H^2(Z_+, \Theta_{Z_+,D_+\bar{D}_+Q} )=0$.
But \cite[Lemma 2]{LS2} says that
$H^2(Z_+, \Theta_{Z_{+}, Q_{+}} )=0.$
 We now invoke the exact sequence
$$ 0 \longrightarrow \Theta_{Z_{+}, D_{+}\bar{D}_{+}Q_{+}}
     \longrightarrow \Theta_{Z_+, Q} \longrightarrow
     \nu_{D_+\bar{D}_+ , Z_+} \longrightarrow 0.
$$
But because
$\nu_{D_+\bar{D}_+,Z_+}$ is trivial,
it follows that  $H^1(\nu_{D_+\bar{D}_+,Z_+})
=H^1({D_+\bar{D}_+}, \CO) =0,$
since each component of $D_+\bar{D}_+$ is a simply-connected
surface.
This implies  that
$H^2(Z_+, \Theta_{Z_+,D_+\bar{D}_+Q} ) =
  H^2(Z_+, \Theta_{Z_+,Q} )=0$, as desired.
\end{proof}

As a consequence, we have a  commutative diagram
$$\begin{array}{ccccccl}
0       &&    0         &&          0       &&     \\
\downarrow &          &\downarrow &        &\downarrow &    &\\

H^1(\tau^0_{\z,\ddbar})
&{\to}&
H^1(\tau^0_{\ddbar}) \oplus H^1(\tau^0_\z)
 &{\to}&
 H^1(\taurr)  &{\to}&0\\

\downarrow &
 &\downarrow &&   \downarrow  &&   \\

  T^1_f  &{\to}&
 T^1_{\ddbar} \oplus T^1_\z
 &{\to}&
  {\rm Ext}^1_{\ddbar}( \Omega^1_{\z}, {\cal O})
 &\to & T^2_f \to 0 \\

\downarrow &  &\downarrow &&   \downarrow &&   \\

\BC^k &{\to}&
  H^0(\tau^1_\ddbar)  \oplus H^0(\tau^1_\z)
  &{\to}&
  H^0(\tau^1_{\z |\ddbar})  &\to & 0 \\

\downarrow &          &\downarrow &        &\downarrow &    &   \\
  0  &&       0      &&     0     &&
\end{array} $$
with exact rows and columns, and where the middle
row is (\ref{middle}). In particular,
 $T^2_f=0$, so there exists
a versal deformation
$$\begin{array}{rlcc}
{\cal D}\sqcup \bar{\cal D}
\hookrightarrow &{\cal Z}& &\\
 _{\varpi_D} \searrow & \downarrow_{\varpi_Z}&&\\
&U&\stackrel{\mbox{\tiny open}}{\subset} & T^1_f
\end{array}
$$
such that the fibers of $\varpi_Z$
and $\varpi_D$ are smooth over elements of
$T^1_{f}$ which project to
 $(\BC^{\times})^k\subset \BC^k$.
Furthermore, by the same argument as \cite[\S 6.1]{DF},
there is a real structure $\hat{\sigma}: {\cal Z}\to
{\cal Z}$ which interchanges
${\cal D}$ and $\bar{\cal D}$
and which restricts to the central fiber
as the given real structure $\sigma$ on $Z_0$.
This induces a complex conjugation on
$T^1_f$ compatible with the standard one
on $\BC^k$, and the fibers $Z_t$ over points $t$ of the real slice
which project to $(\br^{\times})^k\subset \br^k$
are   twistor
spaces \cite{DF,LS2}. Moreover,  these
twistor spaces contain  degree-1 divisors
$D_t$ which are disjoint from their images
$\bar{D}_t$ under the real structure; and
since $b_1(D_t)=b_1(M)$ is even, it therefore follows
 \cite{boyer,pont} that any such $Z_t$ is the twistor space of
of a scalar-flat K\"ahler metric on $D_t$.

To prove the Quotient Theorem, it thus
suffices to show that there are suitable real values
of $t\in T_f$ for which $D_t$ is biholomorphic to $M$.
In order to show that this is possible, we will use the
following:

\begin{lem} For $N$ as above, the natural map
$$H^1(\tau^0_{\z,\ddbar}) \to H^1(\tau^0_\ddbar)$$
is a surjection.
\label{zot}
\end{lem}

\begin{proof}
The same normalization sequence
of the preceeding lemma tells us that the
natural map
$$ H^1(\tau^0_{\z,D_0\bar{D}_0}) \to
   H^1(\Theta_{Z_-, D_-\bar{D}_-Q_-}) \oplus
   H^1(\Theta_{Z_+, D_+\bar{D}_+Q}) $$
is a surjection. On the other hand,
the analogous exact sequence
$$0 \longrightarrow \tau^0_\d \longrightarrow
 q_* \Theta_{D_-{\sqcup}D_-+, \ell_- \sqcup \ell_+}
 \longrightarrow i_* \Theta_\ell
 \longrightarrow 0 $$
tells us that the natural map
$$ H^1(\tau^0_\d) \longrightarrow
H^1(\Theta_{D_-, \ell_-} ) \oplus H^1(\Theta_{D_+,\ell_+})$$
is an isomorphism; treating   $\bar{D} _0$ similarly
then tells us that
$$H^1(\tau^0_{D_0\bar{D}_0}) = H^1(\Theta_{D_-,\ell_-})
\oplus H^1(\Theta_{\bar{D}_-, {\bar{\ell}_-}})
\oplus H^1(\Theta_{D_+,\ell_+})
\oplus H^1(\Theta_{\bar {D}_+, \bar{\ell}_+}).
$$
It therefore suffices to show  that the natural maps
$H^1(\Theta_{Z_\pm,{D_\pm\bar{D}_\pm}Q_\pm}) \longrightarrow
H^1(\Theta_{D_\pm,\ell_\pm}) \oplus
H^1(\Theta_{\bar{D}_\pm, \bar{\ell}_\pm})$ are surjective.
But as these maps occur in the long exact sequences
induced by
$$0 \longrightarrow
\Theta_{Z_\pm,Q_\pm} \otimes {\cal I}_{D_\pm\bar{D}_\pm}
\longrightarrow
\Theta_{Z_\pm,D_\pm\bar{D}_\pm Q_\pm}
\longrightarrow
\Theta_{D_\pm,\ell_\pm} \oplus \Theta_{\bar{D}_\pm, \bar{\ell}_\pm}
\longrightarrow  0 ,$$
we need merely show that
$H^2(\Theta_{Z_\pm, Q_\pm} \otimes {\cal I}_{D_\pm\bar{D}_\pm}) =0$.

Now for $Z_+$, one has a short exact sequence
$$ 0 \longrightarrow
{\cal V} \otimes {\cal I}_{D_\pm\bar{D}_\pm}
\longrightarrow
\Theta_{Z_\pm,Q_\pm} \otimes {\cal I}_{D_\pm\bar{D}_\pm}
\longrightarrow
\CO
\longrightarrow 0,$$
where
${\cal V}$ denotes the vertical, Q-relative tangent sheaf
 \cite[3.2]{LS2}. Since
$\pi_*^q\CO=0$, $q=1,2$, it follows that $H^j(Z_+, \CO )=
H^j(\bcp_1, \CO )$, and thus
$H^2(\Theta_{Z_+, Q}
 \otimes {\cal I}_{D_+\bar{D}_+}) = H^2(
{\cal V} \otimes {\cal I}_{D_\pm\bar{D}_\pm})$.
On the other hand \cite{LS2}, $\pi^1_{*}({\cal V} \otimes
 {\cal I}_{D_+\bar{D}_+})$ vanishes on the complement
of two points of $\bcp_1$,  so that $H^1(\pi^1_{*}({\cal V} \otimes
 {\cal I}_{D_+\bar{D}_+}))=0$, whereas $\pi^2_{*}({\cal V} \otimes
 {\cal I}_{D_+\bar{D}_+})$ vanishes outright. The claim therefore
follows from the Leray spectral sequence of
$\pi : Z_+\to \bcp_1$.

As for $Z_-$, an argument analogous to that of Lemma 2
shows that  \linebreak
$H^2(\Theta_{Z_-,Q_-} \otimes {\cal I}_{D_-\bar{D}_-})
\cong
[H^2( \Theta_{\tilde Z _{N},Q} \otimes {\cal I}_{D_N \bar{D}_N})]_{\Phi}$.
Indeed, for the branched covering map $\alpha$ and blowing down map $\beta$,
one can check in a similar way as we have done before that
$$  \alpha^j_{\ast}( \Theta_{\tilde{Z_N},Q} \otimes
{\cal I}_{\widetilde{D_N{\bar{D}_N}}}) =
\left\{ \begin{array}{ll}
(\Theta_{Z_-,Q_-} \otimes {\cal I}_{D_-\bar{D}_-}
\oplus [
\Theta_{Z_-,Q_-} \otimes {\cal I}_{D_-\bar{D}_-})
\otimes {\cal L}] &j=0\\
0&j\neq 0,\end{array}\right.$$
 while
$$  \beta^j_{\ast}(\Theta_{\tilde{Z_N},Q} \otimes
{\cal I}_{\widetilde{D_N{\bar{D}_N}}}) =
\left\{ \begin{array}{ll}  \Theta_{Z_N,L_{\Phi}}
 \otimes {\cal I}_{D_N\bar{D _N}}
&j=0\\0  &j\neq 0,\end{array}\right.$$
Hence
$H^2(\Theta_{\tilde Z _{N},Q}\otimes\i_{\tilde{D}_{N}\tilde{\bar{D}}_{N}})
 \cong H^2(\Theta_{Z_{N},L_{\Phi}}\otimes\i_{D_{N}\bar{D}_{N}}).$
But now the long exact sequence of
$$ 0 \longrightarrow
 \Theta_{Z_N,L_{\Phi}} \otimes {\cal I}_{D_N \bar{D}_N}
\longrightarrow
 \Theta_{Z_N} \otimes {\cal I}_{D_N \bar{D}_N}
\longrightarrow
\nu_{L_{\Phi},Z_N} \otimes {\cal I}_{D_N \bar{D}_N}
\longrightarrow
0$$
and the fact that $\nu_{L_{\Phi},Z_N} \otimes {\cal I}_{D_N \bar{D}_N}
\cong \CO(-1) \oplus \CO(-1) $
combine to imply that
$H^2( \Theta_{Z_N,L_{\Phi}} \otimes {\cal I}_{D_N \bar{D}_N})
=
H^2( \Theta_{Z_N} \otimes {\cal I}_{D_N \bar{D}_N})$,
so that
$H^2(\Theta_{Z_-,Q_-} \otimes {\cal I}_{D_-\bar{D}_-})
\cong
[H^2( \Theta_{\tilde Z _{N},Q} \otimes {\cal I}_{D_N \bar{D}_N})]_{\Phi}$,
as claimed. Since $H^2( \Theta_{Z_N} \otimes {\cal I}_{D_N \bar{D}_N})
=0$ by Theorem 2, the result follows.
\end{proof}

Now, according to Proposition \ref{surf},
$D_0$
has a versal family of deformations
 over a neighborhood of
$0\in H^1(D_-,\tau^0_{D_0})\times {\Bbb C}^k$
with the property that the  fiber over
$(0,  t_2)$, $t_2\in (\BC^{\times})^k$,
is biholomorphic to
$M$; moreover, the Kodaira-Spencer map of
this family at $0$ is compatible with the natural exact sequence
$$0\to H^1(\tau^0_\d)\to  T^1_\d  \to H^0(\tau^1_\d)\to 0,$$
for which it therefore provides a splitting.
By versality, the family $\varpi_D$ is therefore
induced by a map $$T^1_f\supset U\to H^1(\tau^0_\d)
\oplus \overline{H^1(\tau^0_\d)}\oplus {\Bbb C}^k
 \oplus \overline{{\Bbb C}^k}$$
 whose derivative at $0$ amounts to the natural
restriction homomorphism $T^1_f\to T^1_{\ddbar}$
and
which intertwines the real structure of
$T^1_f$ with the anti-linear map on the target which
interchanges the obvious pairs of factors.
The derivative  of the induced map
$U\to H^1(\tau^0_\d)
\oplus \overline{H^1(\tau^0_\d)}$
is therefore surjective at $0$ by Lemma \ref{zot};
and if we restrict this map to the real slice and then
project to the first factor, the resulting map
$(\Re T^1_f)\cap U\to H^1(\tau^0_\d)$ therefore also has
surjective derivative at $0$.
The inverse image of $0$ is thus a
$k$-dimensional real submanifold $V$ of $\Re T^1_f$,
and the image of $V$ in
$T^1_{\ddbar}=H^1(\tau^0_\d)
\oplus \overline{H^1(\tau^0_\d)}\oplus {\Bbb C}^k
 \oplus \overline{{\Bbb C}^k}$
is a neighborhood of $0$ in a diagonally imbedded
$\br^k\subset {\Bbb C}^k
 \oplus \overline{{\Bbb C}^k}$.
Since the generic element of $V$ therefore projects
to an element of $({\Bbb C}^{\times})^k$, it follows that
the
restriction of  $\varpi_Z$ to $V$
is a simultaneous real
smoothing of $(Z_0, D_0)$ for which  the
generic hypersurface $D_t$ is a copy of $M$.
Thus $M$ admits scalar-flat
K\"ahler metrics, and we have proved
 Theorem \ref{key}.
\hfill \rule{.5em}{1em}\mbox{}\bigskip

\begin{remark}
The above argument
actually proves a bit more; namely,
any small  deformation of $M$ also
admits scalar-flat K\"ahler metrics.
This ostensibly stronger statement,
however, is actually an immediate formal consequence
of the mere statement of Theorem \ref{key}
in light of   the deformation theory of
scalar-flat K\"ahler surfaces \cite{LS1}
together with
Theorem 2.
\end{remark}

\section{The Main Theorems}
As a first step toward proving  our main results,
we now apply the results of the last section to
scalar-flat K\"ahler metrics on some
specific surfaces.

\begin{propn}  If $\bcp_1\times\bcp_1$ is blown up
at $13$ suitably chosen points, the resulting
complex surface admits scalar-flat K\"ahler metrics. \label{pp}
\end{propn}
\begin{proof}
The strategy is to apply Theorem \ref{key}
when $N$ is a
a two-fold blow-up of
${\Bbb CP}_1\times \Sigma$, where $\Sigma$
is a compact complex curve of genus 2.
 We therefore
begin by constructing such a scalar-flat K\"ahler
metric on this  manifold which admits a suitable
involution $\Phi$. This will be done by
 careful use of the
hyperbolic ansatz construction of \cite{L}.

Let $\Sigma$ be a  of genus 2, and let $\phi : \Sigma \to \Sigma$
be the Weierstra\ss\ involution
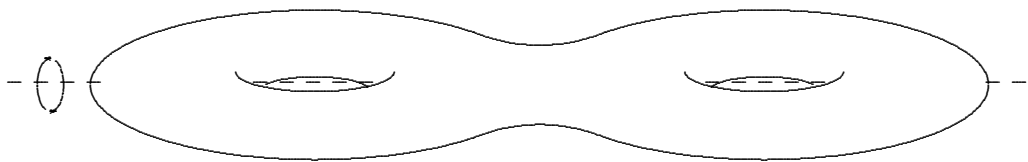
\begin{center}
\mbox{
\beginpicture
\setplotarea x from -80 to 290, y from 25 to 105
\ellipticalarc axes ratio 3:1  270 degrees from 100 80
center at 40 60
\ellipticalarc axes ratio 3:1  -270 degrees from 150 80
center at 210 60
\ellipticalarc axes ratio 4:1 -180 degrees from 70 65
center at 40 65
\ellipticalarc axes ratio 4:1 145 degrees from 60 59
center at 40 58
\ellipticalarc axes ratio 4:1 180 degrees from 180 65
center at 210 65
\ellipticalarc axes ratio 4:1 -145 degrees from 190 59
center at 210 58
\ellipticalarc axes ratio 1:2 160 degrees from -60 70
center at -60 60
\ellipticalarc axes ratio 1:2 160 degrees from -60 50
center at -60 60
\arrow <2pt> [.2,.7] from -60 70 to -62 70
\arrow <2pt> [.2,.7] from -60 50 to -58 50
{\setquadratic
\plot 100 80    125 75    150 80   /
\plot 100 40    125 45    150 40  /
}
{\setlinear
\setdashes
\plot  13 61 60 61   /
\plot  184 61  230 61 /
\plot  -80 61 -40 61   /
\plot  290 61 310 61   /
}
\endpicture
}
\end{center}
\noindent which realizes $\Sigma$  as a 2-sheeted
branched cover  $\pi: \Sigma\to {\Bbb CP}_1$; let
$\hat{q}\in \Sigma$ be one of the 6 fixed points of $\phi$,
and set $q=\pi(\hat{q})$. Let $h_{\Sigma}$
be the curvature -1 Hermitian metric on $\Sigma$, and notice
that  $\phi$ is an isometry of $h_{\Sigma}$.

We now equip the 3-manifold
$X:=\Sigma \times (-1,1)$ with the hyperbolic metric
$$ h=\frac{h_{\Sigma}}{(1-t^2)}+\frac{dt^2}{(1-t^2)^2}~ ,$$
where $t$ is the standard coordinate on $(-1,1)$.
Let $p_{\pm}=(\hat{q}, \pm \frac{1}{2})\in X$,
let $G_{\pm}$ be the hyperbolic Green's functions of $p_{\pm}\in X$,
and set $V=1+G_+ +G_-$. We then let $P$ be the principal
$S^1$-bundle on $X-\{ p_{\pm}\}$ with connection 1-form
$\theta$ such that
$$\star dV = d\theta$$
and such that the restriction of $(P,\theta)$ to the hypersurface
$t=0$ is the trivial bundle-with-connection $\Sigma\times S^1$.
We then endow $P$ with the Riemannian metric
$$g=(1-t^2)[Vh+V^{-1}\theta^2].$$
The metric space completion of $(P,g)$ is then a smooth compact  Riemannian
4-manifold $(N,g_N)$ of scalar-curvature zero, and  admits a complex
structure $J_N$ with respect to which $g_N$ is K\"ahler.
Moreover \cite{L}, the complex surface  $(N,J_N)$ is biholomorphic
to  ${\Bbb CP}_1\times \Sigma$ blown up at two points in the fiber over
$q\in \Sigma$.

Consider the map $\psi : \Sigma \times (-1,1)\to \Sigma \times (-1,1)$ given
by $(\zeta , t ) \mapsto (\phi (\zeta ) , -t )$. As
$\psi^*V=V$ and $\psi$ is an orientation-{\em reversing} isometry of
$X$, it follows that $\psi^*P\cong \bar{P}$
as a principal bundle-with-connection, where
$\bar{P}$ denotes $P$ equipped with the inverse $S^1$-action.
There is therefore a unique isometry $\Phi$ of $P$ which covers
$\psi$  and restricts to the hypersurface $t=0$
as $\phi \times c : \Sigma\times S^1\to \Sigma\times S^1$,
where $c : S^1\to S^1$ is the reflection $e^{i\vartheta}
\to e^{-i\vartheta}$.
This extends to the Riemannian completion $N$ as an involution of the
desired type. Indeed, it is not hard to see that $\Phi$
is induced by
 the involution $r\times \phi$ of ${\Bbb CP}_1\times \Sigma$,
where $r$ is $180^{\circ}$ rotation of ${\Bbb CP}_1=S^2$ about an axis,
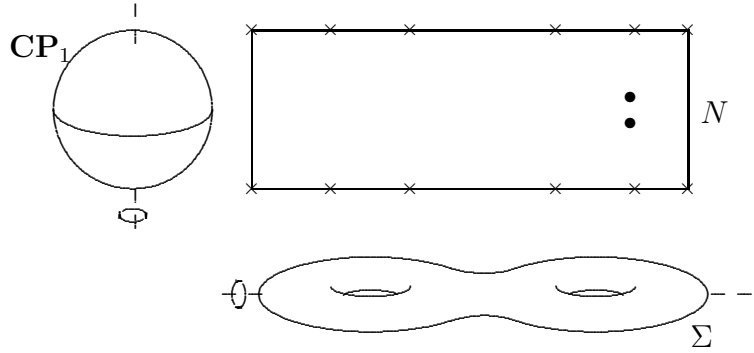
\begin{center}
\mbox{
\beginpicture
\setplotarea x from 0 to 290, y from -5 to 140
\putrectangle corners at 75 130 and 240 70
\put {$_{\times}$} [B1] at 75 130
\put {$_{\times}$} [B1] at 105 130
\put {$_{\times}$} [B1] at 135 130
\put {$_{\times}$} [B1] at 190 130
\put {$_{\times}$} [B1] at 220 130
\put {\circle*{4}} [B1] at 220 105
\put {\circle*{4}} [B1] at 220 95
\put {$_{\times}$} [B1] at 240 130
\put {$_{\times}$} [B1] at 75 70
\put {$_{\times}$} [B1] at 105 70
\put {$_{\times}$} [B1] at 135 70
\put {$_{\times}$} [B1] at 190 70
\put {$_{\times}$} [B1] at 220 70
\put {$_{\times}$} [B1] at 240 70
\put {$\Sigma$} [B1] at 245 10
\put {${\Bbb CP}_1$} [B1] at -5 120
\put {$N$} [B1] at 250 95
\circulararc 360 degrees from 30 130 center at 30 100
\ellipticalarc axes ratio 3:1 -180 degrees from 60 100
center at 30 100
\ellipticalarc axes ratio 3:1  270 degrees from 150 40
center at 120 30
\ellipticalarc axes ratio 3:1  -270 degrees from 175 40
center at 205 30
\ellipticalarc axes ratio 4:1 -180 degrees from 135 33
center at 120 33
\ellipticalarc axes ratio 4:1 145 degrees from 130 30
center at 120 29
\ellipticalarc axes ratio 4:1 180 degrees from 190 33
center at 205 33
\ellipticalarc axes ratio 4:1 -145 degrees from 195 30
center at 205 29
\ellipticalarc axes ratio 1:2 140 degrees from 70 35
center at 70 30
\ellipticalarc axes ratio 1:2 140 degrees from 70 25
center at 70 30
\ellipticalarc axes ratio 2:1 140 degrees from 35 60
center at 30 60
\ellipticalarc axes ratio 2:1 140 degrees from 25 60
center at 30 60
\arrow <2pt> [.1,.3] from 70 35 to 71 35
\arrow <2pt> [.1,.3] from 70 25 to 69 25
\arrow <2pt> [.1,.3] from 35 60 to 35 59
\arrow <2pt> [.1,.3] from 25 60 to 25 61
{\setquadratic
\plot 150 40    163 38    175 40   /
\plot 150 20    163 22    175 20  /
}
{\setlinear
\setdashes
\plot  27 140 27 125  /
\plot  27 70 27 55  /
\plot  60 30 80 30   /
\plot  245 30 265 30   /
}
\endpicture
}
\end{center}
which interchanges the two blown-up  points, which are
antipodal on the equator of an invariant $S^2$.

One may now apply Theorem~\ref{key},
but it remains for us
to understand the structure of the resulting $(M,J_M)$.
In order to do this, first blow up $N$ at the $12$ fixed
points of $\Phi$, and notice that we have the
following arrangement of curves in the blow-up $\tilde{N}$:
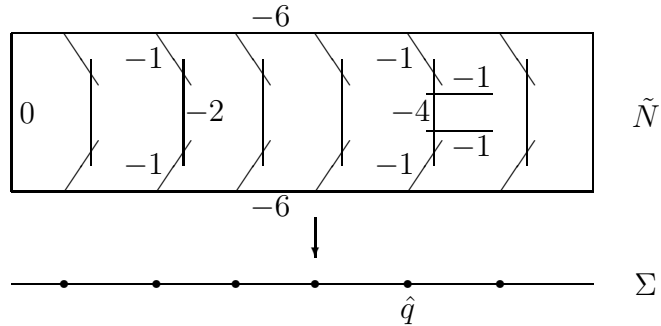
\begin{center}
\begin{picture}(240,110)(0,3)
\put(-10,110){\line(1,0){220}}
\put(-10,110){\line(0,-1){60}}

\put(10,110){\line(2,-3){13}}
\put(20,100){\line(0,-1){40}}
\put(10,50){\line(2,3){13}}

\put(45,110){\line(2,-3){13}}
\put(55,100){\line(0,-1){40}}
\put(45,50){\line(2,3){13}}
\put(40,100){\makebox(0,0){$-1$}}
\put(40,60){\makebox(0,0){$-1$}}
\put(63,80){\makebox(0,0){$-2$}}

\put(75,110){\line(2,-3){13}}
\put(85,100){\line(0,-1){40}}
\put(75,50){\line(2,3){13}}

\put(105,110){\line(2,-3){13}}
\put(115,100){\line(0,-1){40}}
\put(105,50){\line(2,3){13}}

\put(140,110){\line(2,-3){13}}
\put(150,100){\line(0,-1){40}}
\put(140,50){\line(2,3){13}}
\put(147,73){\line(1,0){25}}
\put(147,87){\line(1,0){25}}
\put(164,93){\makebox(0,0){$-1$}}
\put(164,67){\makebox(0,0){$-1$}}
\put(135,100){\makebox(0,0){$-1$}}
\put(135,60){\makebox(0,0){$-1$}}
\put(141,80){\makebox(0,0){$-4$}}

\put(175,110){\line(2,-3){13}}
\put(185,100){\line(0,-1){40}}
\put(175,50){\line(2,3){13}}

\put(88,44){\makebox(0,0){$-6$}}
\put(88,116){\makebox(0,0){$-6$}}
\put(230,80){\makebox(0,0){$\tilde{N}$}}
\put(230,15){\makebox(0,0){$\Sigma$}}

\put(-10,50){\line(1,0){220}}
\put(210,110){\line(0,-1){60}}
\put(105,40){\vector(0,-1){15}}
\put(-10,15){\line(1,0){220}}
\put(10,15){\circle*{3}}
\put(45,15){\circle*{3}}
\put(75,15){\circle*{3}}
\put(105,15){\circle*{3}}
\put(140,15){\circle*{3}}
\put(175,15){\circle*{3}}

\put(140,5){\makebox(0,0){$\hat{q}$}}
\put(-4,80){\makebox(0,0){$0$}}
\end{picture}
\end{center}
Now $\tilde{N}$ is a $2$-fold branched cover of $M$,
with ramification locus equal to the  twelve $(-1)$-curves
introduced by  the blowing up the fixed points. Descending
to $M$ will thus double the self-intersection of these branch curves, while
halving the self-intersection of any curve on which $\Phi$
acts non-trivially.
The corresponding picture of  $M$ is therefore as follows:

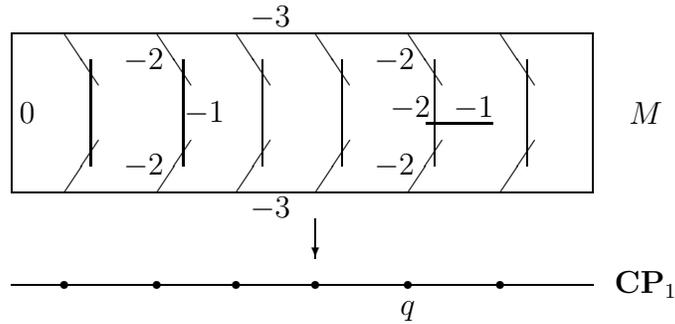
\begin{center}
\begin{picture}(240,120)(0,3)
\put(-10,110){\line(1,0){220}}
\put(-10,110){\line(0,-1){60}}

\put(10,110){\line(2,-3){13}}
\put(20,100){\line(0,-1){40}}
\put(10,50){\line(2,3){13}}

\put(45,110){\line(2,-3){13}}
\put(55,100){\line(0,-1){40}}
\put(45,50){\line(2,3){13}}
\put(40,100){\makebox(0,0){$-2$}}
\put(40,60){\makebox(0,0){$-2$}}
\put(63,80){\makebox(0,0){$-1$}}

\put(75,110){\line(2,-3){13}}
\put(85,100){\line(0,-1){40}}
\put(75,50){\line(2,3){13}}

\put(105,110){\line(2,-3){13}}
\put(115,100){\line(0,-1){40}}
\put(105,50){\line(2,3){13}}

\put(140,110){\line(2,-3){13}}
\put(150,100){\line(0,-1){40}}
\put(140,50){\line(2,3){13}}
\put(147,76){\line(1,0){25}}
\put(165,82){\makebox(0,0){$-1$}}
\put(135,100){\makebox(0,0){$-2$}}
\put(135,60){\makebox(0,0){$-2$}}
\put(141,82){\makebox(0,0){$-2$}}

\put(175,110){\line(2,-3){13}}
\put(185,100){\line(0,-1){40}}
\put(175,50){\line(2,3){13}}

\put(88,44){\makebox(0,0){$-3$}}
\put(88,116){\makebox(0,0){$-3$}}
\put(230,80){\makebox(0,0){$M$}}
\put(230,15){\makebox(0,0){$\bcp_1$}}

\put(-10,50){\line(1,0){220}}
\put(210,110){\line(0,-1){60}}
\put(105,40){\vector(0,-1){15}}
\put(-10,15){\line(1,0){220}}
\put(10,15){\circle*{3}}
\put(45,15){\circle*{3}}
\put(75,15){\circle*{3}}
\put(105,15){\circle*{3}}
\put(140,15){\circle*{3}}
\put(175,15){\circle*{3}}

\put(140,5){\makebox(0,0){$q$}}
\put(-4,80){\makebox(0,0){$0$}}
\end{picture}
\end{center}
\noindent Contracting  $13$ judiciously chosen exceptional  curves,

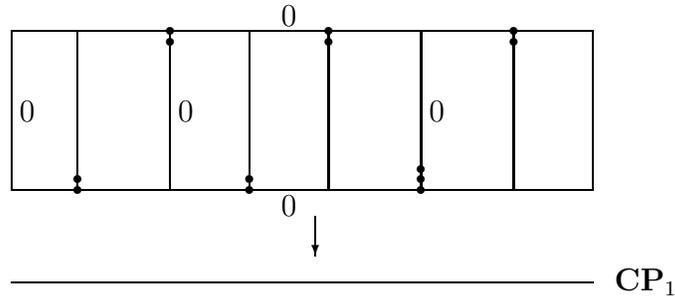
\begin{center}
\begin{picture}(240,120)(0,3)
\put(-10,110){\line(1,0){220}}
\put(-10,110){\line(0,-1){60}}

\put(15,110){\line(0,-1){60}}
\put(50,110){\line(0,-1){60}}
\put(80,110){\line(0,-1){60}}
\put(110,110){\line(0,-1){60}}
\put(145,110){\line(0,-1){60}}
\put(180,110){\line(0,-1){60}}

\put(15,50){\circle*{3}}
\put(80,50){\circle*{3}}
\put(145,50){\circle*{3}}
\put(15,54){\circle*{3}}
\put(80,54){\circle*{3}}
\put(145,54){\circle*{3}}
\put(145,58){\circle*{3}}

\put(50,110){\circle*{3}}
\put(110,110){\circle*{3}}
\put(180,110){\circle*{3}}
\put(50,106){\circle*{3}}
\put(110,106){\circle*{3}}
\put(180,106){\circle*{3}}

\put(95,44){\makebox(0,0){$0$}}
\put(95,116){\makebox(0,0){$0$}}
\put(230,15){\makebox(0,0){$\bcp_1$}}

\put(-10,50){\line(1,0){220}}
\put(210,110){\line(0,-1){60}}
\put(105,40){\vector(0,-1){15}}
\put(-10,15){\line(1,0){220}}

\put(-4,80){\makebox(0,0){$0$}}
\put(56,80){\makebox(0,0){$0$}}
\put(151,80){\makebox(0,0){$0$}}
\end{picture}
\end{center}
we get $\bcp_1\times\bcp_1$ as our minimal model. Invoking
Theorem \ref{key}, we thus conclude that the above iterated blow-up
of $\bcp_1\times\bcp_1$ admits scalar-flat K\"ahler metrics.
Moreover, it also follows from our smoothing argument
  that any sufficiently small deformation

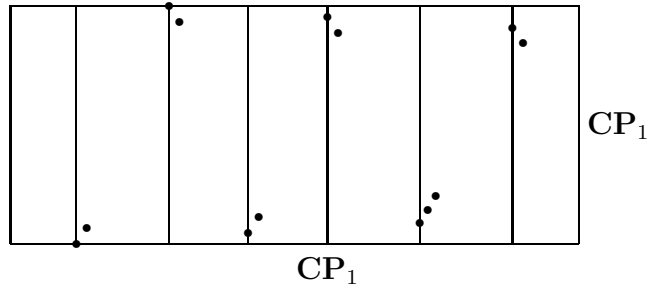
\begin{center}
\begin{picture}(240,120)(0,3)
\put(220,65){\makebox(0,0){$\bcp_1$}}
\put(110,10){\makebox(0,0){$\bcp_1$}}

\put(-10,110){\line(1,0){215}}
\put(-10,110){\line(0,-1){90}}
\put(205,110){\line(0,-1){90}}
\put(15,110){\line(0,-1){90}}
\put(50,110){\line(0,-1){90}}
\put(80,110){\line(0,-1){90}}
\put(110,110){\line(0,-1){90}}
\put(145,110){\line(0,-1){90}}
\put(180,110){\line(0,-1){90}}

\put(15,20){\circle*{3}}
\put(80,24){\circle*{3}}
\put(145,28){\circle*{3}}
\put(19,26){\circle*{3}}
\put(84,30){\circle*{3}}
\put(148,33){\circle*{3}}
\put(151,38){\circle*{3}}

\put(50,110){\circle*{3}}
\put(110,106){\circle*{3}}
\put(180,102){\circle*{3}}
\put(54,104){\circle*{3}}
\put(114,100){\circle*{3}}
\put(184,96){\circle*{3}}
\put(-10,20){\line(1,0){215}}
\end{picture}
\end{center}
of $M$ also admits scalar-flat K\"ahler metrics.
\end{proof}

\begin{cor}
If ${\Bbb CP}_2$ is blown up at 14 suitably chosen points,
the resulting complex surface admits scalar-flat
K\"ahler metrics.
\end{cor}
\begin{proof}
The blow-up of ${\Bbb CP}_2$ at two distinct points is
biholomorphic to the blow-up of ${\Bbb CP}_1\times
{\Bbb CP}_1$ at one point. Now use Proposition \ref{pp}.
\end{proof}

\begin{propn} Let $E\approx T^2$ be any elliptic curve.
\label{ep}
If  $E\times\bcp_1$  is blown up at $6$ suitably chosen points, the
 resulting
complex surface admits scalar-flat K\"ahler metrics.
\end{propn}
\begin{proof} Our  strategy is similar to that used in the
genus $0$ case, but we will now have to exercise  great
 care in order to
compensate for the  non-trivial Jacobi variety of  $E$.

We begin by choosing 4 distinct points $q_1, q_2, r_1, r_2\in E$
such that  $q_1+ q_2= r_1+ r_2$ as divisors, and a holomorphic line
bundle $L\to E$ which is a square-root of the divisor
$q_1+ q_2= r_1+ r_2$. (For example,
if $E$ is ${\Bbb C}/\Lambda$, where   lattice $\Lambda$ is generated by
$1$ and $\tau$, we may take $q_1, q_2, r_1, r_2$ to respectively be
the equivalence classes of $0$,  $(1+\tau)/2$, $1/2$, and $\tau/2$,
whereas $L$ may be taken to be the divisor of the
point represented by $(1+\tau)/4$.) Thus  $L^{\otimes 2}$ comes equipped with
 holomorphic sections $u_1$ and $u_2$ whose  zero
sets are respectively $\{q_1, q_2\}$ and $\{ r_1, r_2\}$,
and all these zeroes are simple.  We now let $\Sigma$ be the
2-fold branched cover of $E$ with branch locus $\{q_1, q_2\}$
associated to the bundle $L$; explicitly,
$$\Sigma=\{ \zeta \in L ~|~ \zeta\otimes \zeta =u_1  \pi_L (\zeta)\}$$
where $\pi_L:L\to E$ is the canonical projection. Let
$\phi :  \Sigma\to\Sigma$  be the
 involution with $2$ fixed points induced by multiplication by
$-1$ in $L$, and let  $\pi : \Sigma\to E$  denote the
canonical projection induced by $\pi_L$. Let $\pi^{-1}(\{ q_1, q_2\})=
\{ \hat{q}_1, \hat{q}_2\}$,  and let $\pi^{-1}(\{ r_1, r_2\})=
\{ \hat{r}_1, \hat{r}_2, \hat{r}_3, \hat{r}_4\}$.
 Notice that $\pi^*L$ is the divisor line bundle
of $\hat{q}_1+\hat{q}_2$, whereas its square $\pi^*L^{\otimes 2}$
is the divisor of $\hat{r}_1+ \hat{r}_2+ \hat{r}_3+ \hat{r}_4$.

As before, we equip the 3-manifold
$X:=\Sigma \times (-1,1)$ with the hyperbolic metric
$$ h=\frac{h_{\Sigma}}{(1-t^2)}+\frac{dt^2}{(1-t^2)^2}~ ,$$
but this time we set $p_j=(\hat{q}_j, 0)$, $j=1, \ldots , 4$.
Let $G_j$ be the hyperbolic Green's functions of $p_j\in X$, and
set $V=1+\sum_{j=1}^4G_j$. On $(\Sigma \times \{ 0\})-\{p_j\}
\cong \Sigma -\{\hat{q}_j\}$, let $P_0$ be the
flat $S^1$ bundle with ${\Bbb Z}_2$ monodromy corresponding to
the branched cover with branch-points
$\hat{r}_1 , \ldots , \hat{r}_4$ associated with
the line bundle $\pi^*L$.
Then \cite{L}
there exists
  a principal
$S^1$-bundle $P\to X-\{ p_{j}\}$ with connection 1-form
$\theta$ such that
$$\star dV = d\theta$$
and such that the restriction of $(P,\theta)$ to the hypersurface
$t=0$ is the flat bundle $P_0$.
 Indeed, the Chern-Weil
theorem guarantees that we can find a connection with
curvature $\star dV$ because $d \star dV=0$ on
$X-\{ p_j\}$ and $[\frac{1}{2\pi}]\in H^2(X-\{ p_j\}, {\Bbb Z})$.
Because $V$ is symmetric in $t$, a connection with
curvature $\star dV$
is automatically flat on $(\Sigma \times \{ 0\})-\{p_j\}$,
and its holonomy around $p_j$ is automatically
$e^{i\int_{D} \star dV} = e^{\frac{i}{2}\int_S \star dV}=
e^{\frac{i}{2}\int_B d\star dV}=e^{-i\pi}=-1$,
where $D$ is a disk in $t\geq 0$ that bounds a loop around
$p_j$ in $\Sigma \times \{ 0\}$, $S$ is the sphere
made up of $D$ and its reflection in $t$, and
 $B$ is the ball about $p_j\in X$ with boundary $S$.
Twisting by the pull-back of a flat connection
on $\Sigma$ now allows us to modify the  the restriction of
$(P,\theta)$ to $\Sigma \times \{ 0\}$ so as to
obtain  any given
flat connection with holonomy $-1$ round the points
$p_j$; and  $P_0$ fits the bill.

As before, we now
 endow $P$ with the Riemannian metric
$$g=(1-t^2)[Vh+V^{-1}\theta^2].$$
The metric space completion of $(P,g)$ is then a smooth compact  Riemannian
4-manifold $(N,g_N)$ of scalar-curvature zero, and  admits a complex
structure $J_N$ with respect to which $g_N$ is K\"ahler.
There is a smooth holomorphic curve  $\hat{\Sigma}\subset N$
which is a copy of the two-fold cover of  $\Sigma$
with branch points $\hat{r}_1 , \ldots , \hat{r}_4$
associated with $\pi^*L$, obtained by taking
closure of the  ${\Bbb Z}_2$ bundle from which $P_0$
was constructed. This allows one to observe  \cite{L}
that $(N,J_N)$ is obtained from ${\Bbb P}(\pi^*L\oplus {\cal O})\to \Sigma$
by blowing up the points $\hat{r}_1 , \ldots , \hat{r}_4$
on the zero section of
$\pi^*L\subset {\Bbb P}(\pi^*L\oplus {\cal O})$;
the key point is that  $\hat{\Sigma}\subset N$
corresponds to the proper transform of the
curve $\zeta^{\otimes 2}=u_3$, where $u_3\in
\Gamma (\Sigma ,{\cal O}(\pi^*L^{\otimes 2}))$ is
the standard section with
simple zeroes at $\hat{r}_1 , \ldots , \hat{r}_4$.

Consider the involution $\iota : \pi^*L\to \pi^*L$
which maps  $(\pi^*L)_z$ to $(\pi^*L)_{\phi(z)}$
by multiplication by $-1$ in $L_{\pi (z)}=L_{\pi (\phi (z))}$.
This involution extends ${\Bbb P}(\pi^*L\oplus {\cal O})$,
and so lifts to an involution  $\Phi : N\to N$
with exactly $4$ fixed points, corresponding to
$0$ and $\infty$  in the fibers over $\hat{q}_1$ and $\hat{q}_2$.
Since $\pi^*L$ is the divisor line bundle of
$\hat{q}_1+\hat{q}_2$, there is a section $u$ of
$\pi^*L$  with simple zeroes at $\hat{q}_1$ and $\hat{q}_2$,
and by averaging we may arrange that $\iota (u(z))=-u(\phi (z))$
for all $z\in \Sigma$; the proper transform $C_u\subset N$ of the
image of $u$ is then $\Phi$-invariant.

Now the involution $\Phi$ just  corresponds
 to the unique connection-preserving  involution of $P$
which covers $\phi \times 1: [\Sigma \times (-1,1)]\to
[\Sigma \times (-1,1)]$ and extends
 $\iota: \hat{\Sigma}\to \hat{\Sigma}$.
It therefore  preserves $g_N$, and we
may  apply Theorem~\ref{key}.
The remaining task is thus to analyze the
structure of the surface  $(M,J_M)$.

Since $\tilde{N}$ is obtained from $N$
by blowing up the $4$ fixed points
of $\Phi$, it contains the
following arrangement of curves:

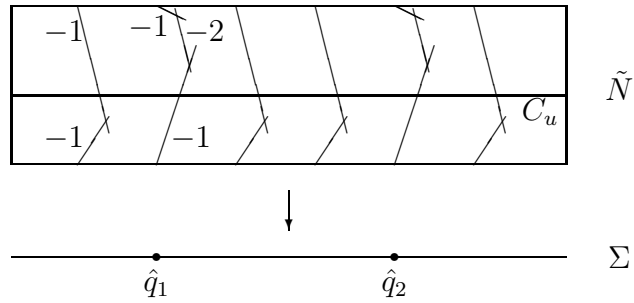
\begin{center}
\begin{picture}(240,110)(0,3)
\put(-10,110){\line(1,0){210}}
\put(-10,110){\line(0,-1){60}}
\put(-10,76){\line(1,0){210}}

\put(15,110){\line(1,-4){12}}
\put(15,50){\line(2,3){12}}
\put(10,100){\makebox(0,0){$-1$}}
\put(10,60){\makebox(0,0){$-1$}}

\put(45,110){\line(2,-1){11}}
\put(52,109){\line(1,-4){6}}
\put(45,50){\line(1,3){15}}
\put(42,102){\makebox(0,0){$-1$}}
\put(58,60){\makebox(0,0){$-1$}}
\put(63,100){\makebox(0,0){$-2$}}

\put(75,110){\line(1,-4){12}}
\put(75,50){\line(2,3){12}}

\put(105,110){\line(1,-4){12}}
\put(105,50){\line(2,3){12}}

\put(135,110){\line(2,-1){11}}
\put(142,109){\line(1,-4){6}}
\put(135,50){\line(1,3){15}}

\put(165,110){\line(1,-4){12}}
\put(165,50){\line(2,3){12}}

\put(190,70){\makebox(0,0){$C_u$}}

\put(220,80){\makebox(0,0){$\tilde{N}$}}
\put(220,15){\makebox(0,0){$\Sigma$}}

\put(-10,50){\line(1,0){210}}
\put(200,110){\line(0,-1){60}}
\put(95,40){\vector(0,-1){15}}
\put(-10,15){\line(1,0){210}}

\put(45,15){\circle*{3}}
\put(135,15){\circle*{3}}

\put(45,5){\makebox(0,0){$\hat{q}_1$}}
\put(135,5){\makebox(0,0){$\hat{q}_2$}}
\end{picture}
\end{center}

The corresponding picture of  $M$ is therefore:

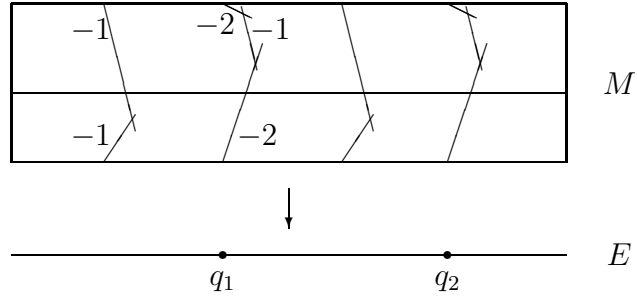
\begin{center}
\begin{picture}(240,110)(0,3)
\put(-10,110){\line(1,0){210}}
\put(-10,110){\line(0,-1){60}}
\put(-10,76){\line(1,0){210}}

\put(25,110){\line(1,-4){12}}
\put(25,50){\line(2,3){12}}
\put(20,100){\makebox(0,0){$-1$}}
\put(20,60){\makebox(0,0){$-1$}}

\put(70,110){\line(2,-1){11}}
\put(77,109){\line(1,-4){6}}
\put(70,50){\line(1,3){15}}
\put(67,102){\makebox(0,0){$-2$}}
\put(83,60){\makebox(0,0){$-2$}}
\put(88,100){\makebox(0,0){$-1$}}

\put(115,110){\line(1,-4){12}}
\put(115,50){\line(2,3){12}}

\put(155,110){\line(2,-1){11}}
\put(162,109){\line(1,-4){6}}
\put(155,50){\line(1,3){15}}

\put(220,80){\makebox(0,0){$M$}}
\put(220,15){\makebox(0,0){$E$}}

\put(-10,50){\line(1,0){210}}
\put(200,110){\line(0,-1){60}}
\put(95,40){\vector(0,-1){15}}
\put(-10,15){\line(1,0){210}}

\put(70,15){\circle*{3}}
\put(155,15){\circle*{3}}

\put(70,5){\makebox(0,0){$q_1$}}
\put(155,5){\makebox(0,0){$q_2$}}
\end{picture}
\end{center}

Contracting  $6$ exceptional  curves in the right order,

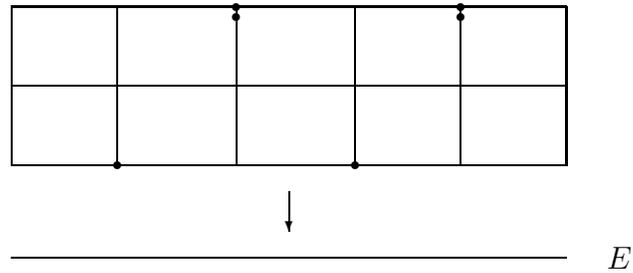
\begin{center}
\begin{picture}(240,110)(0,3)
\put(-10,110){\line(1,0){210}}
\put(-10,110){\line(0,-1){60}}
\put(-10,80){\line(1,0){210}}

\put(30,110){\line(0,-1){60}}
\put(75,110){\line(0,-1){60}}
\put(120,110){\line(0,-1){60}}
\put(160,110){\line(0,-1){60}}

\put(30,50){\circle*{3}}
\put(75,110){\circle*{3}}
\put(75,106){\circle*{3}}
\put(120,50){\circle*{3}}
\put(160,110){\circle*{3}}
\put(160,106){\circle*{3}}

\put(220,15){\makebox(0,0){$E$}}

\put(-10,50){\line(1,0){210}}
\put(200,110){\line(0,-1){60}}
\put(95,40){\vector(0,-1){15}}
\put(-10,15){\line(1,0){210}}
\end{picture}
\end{center}
we get $E\times\bcp_1$ as our minimal model.
Theorem \ref{key} thus tells us that the above iterated blow-up
of $E\times\bcp_1$
admits scalar-flat K\"ahler metrics, as do
 its sufficiently small deformations.
\end{proof}

\begin{remark} Notice that
the scalar-flat K\"ahler surface we have just constructed
admits a holomorphic ${\Bbb C}^{\times}$-action; namely, the
${\Bbb C}^{\times}$-action on $\bcp_1\times E$ induced by the
``earth-rotation'' of  $\bcp_1$ lifts to
$M$, which is obtained from $\bcp_1\times E$ by
iteratively blowing up fixed points of the action. Also notice,
however,  that
the induced
action is  not semi-free, since the
isotropy of any generic point on either of the second-level
blow-up curves is $\{ \pm 1\} \subset {\Bbb C}^{\times}$.
A theorem of  Lichnerowicz \cite{lich} now implies
that    $S^1\subset {\Bbb C}^{\times}$ acts isometrically
on this scalar-flat K\"ahler surface.
By combining elements of the proofs of  Propositions
\ref{pp} and \ref{ep}, one can similarly
construct scalar-flat K\"ahler metrics with non-semi-free
$S^1$-action on 14-point blow-ups of $\bcp_1\times \bcp_1$.

One might, in principle,  use the
Toda lattice equation \cite{l} to construct
the above metrics explicitly. However, Theorem 2 of \cite{Ltan}
insists that these solutions definitely
cannot arise from the hyperbolic ansatz case of that
construction. This is consistent with the fact that
the orbits with exceptional isotropy
give rise to  very peculiar
orbifold singularities of the associated  3-dimensional
Einstein-Weyl geometry.

Any of these metrics is a counter-example to
 the assertion \cite[Proposition 3.1]{LS1}
that a blown-up ruled surface of genus  $<2$ cannot
admit both a ${\Bbb C}^{\times}$-action and
a K\"ahler class of total scalar curvature $0$.
While the argument offered there does indeed work
 if one insists that
the  ${\Bbb C}^{\times}$-action be semi-free,
the recipe given for a section of the anti-canonical
bundle simply breaks down if, as in the present case,
 an isolated fixed point
is blown up at some stage in the iterated blow-up process.
 \end{remark}

\begin{thm} Let $(M,J)$ be a ruled surface--- i.e.
 suppose that $M$ is a
compact complex 2-manifold
for which there exists a holomorphic map  $M\to \Sigma$
with generic fiber $\bcp_1$ and range a
 Riemann surface $\Sigma$.
Then  $(M,J)$ has blow-ups $(\tilde{M},\tilde{J})$
which admit scalar-flat K\"ahler metrics.\label{punch}
\end{thm}
\begin{proof} Any ruled surface $M$ is bimeromorphic \cite{bpv}
to some product surface $\Sigma \times {\Bbb CP}_1$.
If $M_1$ is  any blow-up of $\Sigma \times {\Bbb CP}_1$,
it then follows that there is a blow-up $M_2$ of $M_1$
which is also a blow-up of $M$:
\setlength{\unitlength}{1ex}
\begin{center}\begin{picture}(20,17)(0,3)
\put(10,17){\makebox(0,0){$M_2$}}
\put(2,5){\makebox(0,0){$M$}}
\put(19,5){\makebox(0,0){$M_1$}}
\put(11,15.5){\vector(2,-3){6}}
\put(9,15.5){\vector(-2,-3){6}}
\end{picture}\end{center}
\noindent But \cite[Theorem 4.6]{KP} and Theorem 2
tell us  that any 1-point blow-up of a non-minimal
scalar-flat K\"ahler surface also admits
K\"ahler metrics;  induction then says the
same is also true of iterated blow-ups.
It therefore
suffices to find just one blow-up $M_1$ of each
$\Sigma \times {\Bbb CP}_1$ which
admits a scalar-flat K\"ahler metric.
If the genus of $\Sigma$ is at least 2, the
hyperbolic ansatz \cite{L} then explicitly
constructs such metrics on a 2-point
blow-up $M_1$ of $\Sigma \times {\Bbb CP}_1$.
For genera $0$ and $1$, on the other hand,
the necessary surfaces $M_1$ are constructed
in Propositions \ref{pp} and \ref{ep}.
\end{proof}
\setcounter{main}{0}
\begin{main} Let $(M,J)$ be a compact complex 2-manifold which admits
a K\"ahler metric for which the integral of the
scalar curvature is non-negative. Then precisely
one of the following holds:
\begin{itemize}
\item $(M,J)$ admits a Ricci-flat K\"ahler metric; or
\item any blow-up of  $(M,J)$ has blow-ups $(\tilde{M},\tilde{J})$
which admit
scalar-flat K\"ahler metrics. Moreover, any blow-up of such an
 $(\tilde{M},\tilde{J})$  admits scalar-flat K\"ahler metrics, too.
\end{itemize}
\end{main}
\begin{proof}
By a vanishing theorem of Yau
\cite{yau0}, any compact K\"ahler surface of non-negative
total scalar curvature is either ruled or has
 $c_1^{\Bbb R}=0$.
In the latter case, Yau's solution of the
Calabi conjecture \cite{yau} guarantees the existence of a
Ricci-flat metric on $M$.
If, on the other hand, $M$ is ruled, so is any blow-up $\hat{M}$
of $M$; and
Theorem \ref{punch} tells us that
some blow-up $\tilde{M}$ of $\hat{M}$ therefore
admits a scalar-flat K\"ahler metric.
Finally, any blow-up of a non-minimal
scalar-flat K\"ahler surface also admits
scalar-flat K\"ahler metrics, by \cite[Theorem 4.6]{KP}
and Theorem 2, and this may be appied to $\tilde{M}$ to
prove the last clause.
\end{proof}

\begin{main} Let $(M,J)$ be a compact complex 2-manifold which admits
a K\"ahler metric for which the integral of the
scalar curvature is positive. Then
 any blow-up of  $(M,J)$ has blow-ups $(\tilde{M},\tilde{J})$
which admit
 K\"ahler metrics of constant positive scalar curvature.
Moreover, any blow-up of such an
 $(\tilde{M},\tilde{J})$  also admits  such metrics.
\end{main}

\begin{proof}
The hypothesis says that $c_1\cdot [\omega ]>0$ for
some K\"ahler class $[\omega]$, so that $c_1^{\Bbb R}\neq 0$
and $M$ certainly cannot admit a Ricci-flat metric.
On the other hand,
a straightforward  inverse-function-theorem argument
\cite{lsim} shows that if a
compact complex manifold admits
a non-Ricci-flat scalar-flat K\"ahler metric,
it also admits K\"ahler metrics
of constant positive scalar curvature.
The result therefore follows from Theorem \ref{maj}.
\end{proof}

\section{Appendix: The  Vanishing Theorem}\label{van}
In this appendix, we demonstrate that Kodaira-Spencer theory
is unobstructed for the
twistor spaces of all non-minimal scalar-flat
K\"ahler surfaces. The proof is a direct extension of that
presented in \cite{LS1}, but allows for the possibility of
non-semi-free ${\Bbb C}^{\times}$-actions.

\setcounter{thm}{1}
\begin{thm}
Suppose that $N$ is a non-minimal compact complex surface
with scalar-flat K\"ahler metric $g_N$. Then its
twistor space satisfies $H^2(Z_N, \Theta)=H^2(Z_N, \Theta\otimes
\kappa^{-1/2})=H^2(Z_N, \Theta_{Z,D\bar{D}})=0$, where
$D$ and $\bar{D}$ are canonical divisors associated with
$J_N$ and $-J_N$.
\end{thm}
\begin{proof} The natural homomorphisms
\bea H^2(Z_N, \Theta\otimes
\kappa^{-1/2}) &\to & H^2(Z_N, \Theta_{Z,D\bar{D}}), \\
H^2(Z_N, \Theta_{Z,D\bar{D}})&\to & H^2(Z_N, \Theta)\eea
are  surjective because $H^2(N, \Theta_N )=
H^2(N, {\cal O}(\kappa_N^{-1}))=0$ for any ruled surface $N$.
On the other hand, careful inspection of
 the Penrose transform shows \cite[Theorem 2.7]{LS1} that
$H^2(Z_N, \Theta\otimes
\kappa^{-1/2})$ is canonically identified  with
the kernel of \be
d{\cal F}|_{[\omega]}: H^0(N, \Theta_N)\to A^*, \eel{fut}
where ${\cal F}$ is the Futaki invariant,
$[\omega]$ is the K\"ahler class,
and $A\subset H^{1,1}(N)=H^2(N)$ is the hyperplane
$$\{\alpha \in H^{1,1}(N)~|~c_1\cup \alpha =0\} .$$
Our goal here  will thus be to show that (\ref{fut}) is injective.

Because
\cite{lich} the automorphism group of
$N$ has a compact real form given by the isometry goup of
$g_N$, $H^0(N, \Theta_N)$ is spanned by vector fields
$\Xi$ which generate $\BC^{\times}$ actions which are free on an
open dense, and such that the $S^1$-action generated by
$\xi=\Im \Xi$ is isometric with respect to $g_N$, while
$\Re \Xi$ is globally a gradient vector field because
the contraction of $\Xi$ with any holomorphic 1-form
vanishes identically.
In fact, $H^0(N, \Theta_N)$ is at most 1-dimensional; for if
$\Xi_1$ and $\Xi_2$ are two such fields, we must have
$\Xi_1\wedge \Xi_2\equiv 0$ because \cite{LS1,yau0}
the existence of a scalar-flat K\"ahler metric
forces $H^0(N, \kappa^{-1})=0$. This implies that the
generic orbit of $\Xi_1$ has closure $F\cong \bcp_1$ which is
also the closure of a generic orbit of $\Xi_2$. If
$\Xi_1$ does not have the same zeroes on this
$\bcp_1$ as  $\Xi_2$, the isometry goup contains
an  $SU(2)$ which acts transitively on this 2-sphere,
and the orbits of this $SU(2)$ are all either
holomorphically embedded $\bcp_1$'s or points; but
the latter type of orbit is impossible, because an $S^2=\bcp_1$ orbit
near a fixed point would be contained in the domain of a holomorphic
chart, contradicting the fact that
every holomorphic function on $\bcp_1$ is constant.
Thus every orbit of this $SU(2)$  would be a
$\bcp_1$, and we would conclude that $N$ would be
a minimal ruled surface--- contradicting our hypotheses. Thus  $\Xi_1$
and $\Xi_2$ must have the same zeroes on the generic
orbit-curve $F$, and hence, since both generate
generically free $\BC^{\times}$-actions,
$\Xi_1
= \pm \Xi_2$. This shows that $h^0(N, \Theta_N)\leq 1$,
as claimed.

 As there is nothing to show if
$H^0(N, \Theta_N)=0$, we may therefore
assume henceforth that $H^0(N, \Theta_N)\cong\BC$ and is spanned by
a holomorphic vector field
$\Xi$ whose imaginary part $\xi$ is a
periodic Killing field, with generic minimal period $2\pi$.
By averaging,
we can  represent any other  K\"ahler class by a
metric $g$ which is $S^1$-invariant. With respect to such a
metric $g$, whose K\"ahler form we shall call $\omega$,
let  $t$ be the Hamiltonian function of $\xi$, conventionally
normalized so as to have range of the form $[-a,a]$, and let
$\Sigma$ be the stable quotient $N/\! /{\BC}^{\times}$.
If $t= a$ and/or $t= -a$ are isolated fixed points,
blow them up to obtain a new complex surface $\hat{N}$;
otherwise, set $N=\hat{N}$.
 We then have
 a map $\Pi : \hat{N}\to \Sigma\times [-a,a]$
with  $S^1$-orbits as fibers. If
$p_1, \ldots , p_m$ are the images in
 $\Sigma\times (-a,a)$ of the isolated
fixed points, and if
$X=[\Sigma\times (-a,a)]-\{ p_j\}$,
then on the open dense set $Y=\Pi^{-1}(X)\subset \hat{N}$
we may  express the given  K\"ahler metric $g$
in the form$$
g=w\hat{g}(t)+w~dt^{\otimes 2}+w^{-1}\theta^{\otimes 2}~,
$$
for some positive functions $w>0$ on $X$
and a family  orbifold metrics $\hat{g}(t)$ on
$\Sigma$. Here $\theta$ is the unique 1-form on $Y$
whose kernel is orthogonal to $\xi$ and such that
$\theta (\xi )=1$.
The orbifold points in $X$ of the $\hat{g}(t)$
exactly correspond those $S^1$-orbits
in $M$ which are non-trivial, but which have
period $2\pi/n$ for some $n$. Since these
also coincide with points at which the
$\BC^{\times}$-action has exceptional
isotropy, the orbifold points of $X$ exactly
consist of vertical line segments in
 $\Sigma\times (-a,a)$ which
join two of the $p_j$.
If we let $\jmath$ range over the
set of exceptional curves contained in fibers of
$\hat{N}\to \Sigma$, each such curve is the closure of
a ${\BC}^{\times}$-orbit,  and so we may define
  integers
$m_{\jmath}\geq 1$  as the order of
the associated isotropy group, and
real numbers $t^{-}_{\jmath}<t^{+}_{\jmath}$,  defined
to be  the  minimum and maximum
values, respectively,
 of the Hamiltonian
$t$ on the associated exceptional rational curve $E_{\jmath}$.
Also define $$m_{\jmath}(t)=\left\{ \begin{array}{ll}
m_{\jmath}&t^{-}_{\jmath}<t<t^{+}_{\jmath}\\
1&\mbox{otherwise. }\end{array}\right.$$
Finally, let $t_j$ be the  $t$ coordinate of $p_j$.

Because
$g$, $w$ and $dt$ are geometrically defined,  $\hat{g}(t)$
 is an invariantly defined, $t$-dependent  orbifold K\"ahler metric
on  $\Sigma$ for all $t\not\in \{ t_j\}$; and when
when $t\in \{ t_j\}$, it is
and defined everywhere
aside from a fine number of punctures.
Let $\hat{\omega} (t)$ be the K\"ahler form of $\hat{g}(t)$.
If $C^-$ and $C^+$ denote the ``repulsive'' and
``attractive'' curves
$t=\pm a$ in $\hat{N}$,  we then have
\cite{LS1}
\bea \left.\hat{\omega}\right|_{t=\pm a}&=&0 \\
\left.\frac{d }{d t} \hat{\omega} \right|_{t=-a} &=&
   \left.\hphantom{-}2\omega \right|_{C^-}  \\
\left.\frac{d }{d t} \hat{\omega} \right|_{t=a\hphantom{-}} &=&
 \left. - 2\omega \right|_{C^+}~ ,  \eea
while the density of scalar curvature may be written
globally on $Y\subset M$ as
$$ s\vol = [-2\hat{\rho} + \frac{d^2}{dt^2} \hat{\omega}]\wedge
dt\wedge \theta~ .$$
Here the Ricci form $\hat{\rho}(t)$ of $\hat{g}(t)$
satisfies
$$\frac{1}{2\pi}\int_{\Sigma} \hat{\rho}(t)= \chi (\Sigma ) -
\sum_{\jmath} (1-\frac{1}{m_{\jmath}(t)})$$
for $t\not\in \{t_j\}$ by
the Gauss-Bonnet Theorem for orbifolds \cite{Sat}.
Thus
\bea \int_{N} ts\vol &=& \int_Y ts\vol
\\&=&
\int_Y t[-2\hat{\rho} + \frac{d^2}{dt^2}
\hat{\omega}]\wedge dt\wedge \theta
\\&=&
2\pi\int_{\Sigma}\int_{-a}^at
\frac{d^2}{dt^2} \hat{\omega}~dt -4\pi \int_{-a}^at[\int_{\Sigma}
\hat{\rho}]~dt
\\&=&
2\pi \int_{\Sigma}\left( \left[t\frac{d}{dt} \hat{\omega}
\right]^a_{-a}-\int_{-a}^a
\frac{d\hat{\omega}}{dt} dt\right) \\&&~~~~
-4\pi \int_{-a}^a2\pi\left[ \chi (\Sigma ) -
\sum_{\jmath} (1-\frac{1}{m_{\jmath}(t)})\right]t~dt
\\&=&
2\pi \int_{\Sigma} \left[t\frac{d}{dt} \hat{\omega}-\hat{\omega}
\right]^a_{-a}
-8\pi^2
\sum_{\jmath} \frac{1}{m_{\jmath}}
\int_{t^-_{\jmath}}^{t^+_{\jmath}}t~dt
\\&=&
2\pi a\left[\int_{C^-}2\omega -
\int_{C^+}2\omega \right]
-4\pi^2
\sum_{\jmath} \frac{1}{m_{\jmath}}
[(t^+_{\jmath})^2-(t^-_{\jmath})^2]
\\&=&
\left[\int_{C^-}\omega -
\int_{C^+}\omega \right]\int_F\omega
+2\pi\sum_{\jmath}[(a-t_{\jmath}^+)-(t_{\jmath}^- -(-a))]
\int_{E_{\jmath}}\omega
\\&=&
\left[\int_{C^-}\omega -
\int_{C^+}\omega \right]\int_F\omega
 +\sum_{\jmath}\left[ \sum_{\imath \more \jmath}
m_{\imath}\int_{E_{\imath}}\omega - \sum_{\imath \less \jmath}
m_{\imath}\int_{E_{\imath}}\omega \right]
\int_{E_{\jmath}}\omega
\\&=&
\left[\int_{C^-}\omega -
\int_{C^+}\omega \right]\int_F\omega+
\sum_{\jmath;\imath\less \jmath}(m_{\jmath}-m_{\imath})
   \left[\int_{E_{\imath}}\omega\right] \int_{E_{\jmath}}\omega ~.
\eea
Here  $F$ is any smooth  fiber of $\hat{N}\to \Sigma$,
 and the partial ordering $\less$
on $\{ \jmath\}$ means
that, with respect to the flow of
$\Re \Xi$,
the first exceptional curve precedes
the second in some singular fiber of
$\hat{N}\to \Sigma$.

On the other hand,
the total scalar curvature of such a metric $g$ is
\bea \int_{N} s\vol &=& \int_Y s\vol
\\&=&
\int_Y [-2\hat{\rho} + \frac{d^2}{dt^2}
\hat{\omega}]\wedge dt\wedge \theta
\\&=&
2\pi\int_{\Sigma}\int_{-a}^a
\frac{d^2}{dt^2} \hat{\omega}~dt -4\pi \int_{-a}^a\int_{\Sigma}
\hat{\rho}~dt
\\&=&
2\pi \int_{\Sigma} \left[\frac{d}{dt} \hat{\omega}
\right]^a_{-a}  -4\pi \int_{-a}^a2\pi\left[ \chi (\Sigma ) -
\sum_{\jmath} (1-\frac{1}{m_{\jmath}(t)})\right]~dt
\\&=&
-2\pi \left[\int_{C^-}2\omega +
\int_{C^+}2\omega \right]
-16\pi^2a\chi (\Sigma )+8\pi^2\sum_{\jmath}
(1-\frac{1}{m_{\jmath}})(t^+_{\jmath}-t^-_{\jmath})
\\&=&
-4\pi\left[\int_{C^-}\omega +
\int_{C^+}\omega +
2(1-{\bf g})\int_F\omega
+\sum_{\jmath}
(1-m_{\jmath})\int_{E_{\jmath}}\omega\right]
{}~~,
\eea
where $\bf g$ denotes the genus of $\Sigma$.
When the above vanishes, we then have
\bea {\cal F}(\Xi ,  [\omega ])&:=&
-{\textstyle \frac{1}{2}}\int_{N} t~(s-s_0)\vol
=-{\textstyle \frac{1}{2}}\int_{N} ts\vol
 \\&=&
{\textstyle \frac{1}{2}}\left[\int_{C^+}
\omega -\int_{C^-}\omega\right]\int_F\omega
+{\textstyle \frac{1}{2}}
\sum_{\jmath;\imath\less \jmath}(m_{\imath}-m_{\jmath})
   \left[\int_{E_{\imath}}\omega\right] \int_{E_{\jmath}}\omega
 ~.
\eea
(Here $s_0=\int s~d\mu/\int d\mu$ denotes the average value of the
scalar curvature.) If the $\BC^{\times}$-action
generated by $\Xi$ is semi-free, all the $E_{\jmath}$ terms drop out,
and the formula simplifies to yield that of \cite[Theorem 3.2]{LS1}.

Let ${\cal F}:={\cal F}(\Xi ,  \cdot)$
and let ${\cal S}=\frac{-1}{4\pi}\int s~du=-c_1\cdot [\omega]$,
thought of as
functions on the K\"ahler cone $\subset H^{1,1}(N, {\Bbb R})$.
Letting $\Omega = [\omega ]$ denote the
K\"ahler class, we therefore have
\be {\cal S}(\Omega )
= \Omega (C^-+C^++ 2(1-{\bf g})F+\sum (1-m_{\jmath})E_{\jmath})
\eel{joe}
and
\be 2{\cal F}({\Omega})=
\Omega (F) \Omega (C^+-C^-)+\sum_{\jmath;\imath\less \jmath}
(m_{\imath}-m_{\jmath})\Omega (E_{\imath})  \Omega (E_{\jmath})\eel{sam}
provided that ${\cal S} (\Omega )=0$. Thus, when
 $\Omega \in H^{1,1}(N, {\Bbb R})$ is the K\"ahler
class of a scalar-flat K\"ahler metric $g_N$,
\be d{\cal S}= C^-+C^++ 2(1-{\bf g})F+\sum (1-m_{\jmath})E_{\jmath}
\eel{harry}
and
\be
2d{\cal F}|_{\Omega}\equiv \Upsilon
\bmod d{\cal S} ,\eel{mike}
where
\be \Upsilon := \Omega (F) [C^+-C^-]+
\Omega (C^+-C^-)F+\sum_{\jmath;\imath\less \jmath}
(m_{\imath}-m_{\jmath})[
\Omega (E_{\imath})E_{\jmath}+\Omega (E_{\jmath})E_{\imath}] ; \eel{moe}
here we have identified the  cotangent space
of $H^{1,1}(N,{\Bbb R})=H^2(N,{\Bbb R})$
with  $H_{2}(N,{\Bbb R})$, and each algebraic curve in $\hat{N}$
is used as short-hand for the homology class of
its image in $N$. Our goal is thus to show that $\Upsilon$
is never a multiple of $d\cal S$.

To do this, we begin by considering  the case in which
$N\neq \hat{N}$, which is to say that at least
one of the curves $C^{\pm}$ arises by blowing up
an isolated fixed point of the action.
 By reversing the sign of
$\Xi$ if necessary, we may assume that
$C^+$   arises in this way, and  let $N'$ denote the
surface obtained from $\hat{N}$ by contracting $C^+$.
Because the self-intersection
of the image of $F$ in  $N'$ is  $+1$, $N'$
 is an iterated blow-up of ${\Bbb CP}_2$,
and each generic fiber $F$ corresponds to a projective line.
Moreover, the image of $C^-$ in  ${\Bbb CP}_2$
meets each such projective line in a point, and
so  $C^-$ must be the proper transform of
a projective line; if $n$ is the number of
singular fibers of $\hat{N}\to \Sigma = {\Bbb CP}_1$,
we thus have  $(C^-)^2\leq 1-n$. But plugging
\be
F=\sum_{\mbox{\tiny fiber}} m_{\jmath}E_{\jmath} ~,\eel{pat}
$\Omega (C^+)=0$, and   ${\bf g}=0$ into (\ref{joe}),
the equation ${\cal S}(\Omega )=0$ becomes
\be 0=\Omega (C^-+ 2F+ \sum_{\jmath} (1-m_{\jmath})E_{\jmath})=
\Omega (C^-+ (2-n)F+ \sum E_{\jmath}) , \eel{chuck}
and hence  $n >  2$. Thus
$(C^-)^2\leq 1-n < -1$, so  $C^-$
cannot be a $(-1)$-curve, and
$N'=N$. Hence $F\cdot F=1$, and, invoking (\ref{chuck}),
\bea
F\cdot\Upsilon  &=& F\cdot \left[ - C^-\Omega (F) -\Omega (C^-) F
+ \sum_{\imath\less \jmath}
(m_{\imath}-m_{\jmath})[
\Omega (E_{\imath})E_{\jmath}+\Omega (E_{\jmath})E_{\imath}]\right]\\
&=& -\Omega (F+C^-- \sum_{\imath}(m_{\imath}-1)E_{\imath} ) \\
&=& \Omega (F) , \eea
so that $\Upsilon$ is certainly non-zero, and we need merely
show that it is linearly independent from $d{\cal S}$.

To do this, let $E_0$ be the first curve in some singular fiber,
chosen in such a way that $\Omega (E_0) < \Omega (F) /3$;
the latter is possible because the fact that $n \geq 3$
tells us that $F$ is homologous to the sum of $C^-$
and a collection of exceptional curves, at least 3 of which
are the first in their respective fibers. Let
$E_1$ denote the immediate successor of $E_0$;
the multiplicity of the latter curve is given by
$m_1=-E_0^2$, as may be proved by induction.
 Then
\bea
 E_0\cdot\Upsilon &=&
E_0\cdot \left[ -\Omega (C^-) F
-\Omega ( F) C^- + \sum_{\imath\less \jmath}
(m_{\imath}-m_{\jmath})[
\Omega (E_{\imath})E_{\jmath}+\Omega (E_{\jmath})E_{\imath}]\right]\\
&=& -\Omega (F) + E_0\cdot
[\sum_{0\less \jmath}
(1-m_{\jmath}) \Omega (E_{\jmath})E_0]\\
&&+ E_0\cdot [\sum_{1\less \jmath}
(m_1-m_{\jmath}) \Omega (E_{\jmath})E_1]+ E_0\cdot [\sum_{\imath\less 1}
(m_{\imath}-m_1)
\Omega (E_{\imath})E_1]\\
&=& -\Omega (F) + E_0^2  \sum_{\mbox{\tiny fiber}} (1-m_{\jmath})
\Omega (E_{\jmath})+ \sum_{\mbox{\tiny fiber}} (-E^2_0-m_{\jmath})
\\&&\hspace{2in} -(m_1-m_0)\Omega(E_0) + (m_0-m_1) \Omega (E_0)\\
&=& -(2+E^2_0) \Omega(F) + 2(1 + E_0^2)\Omega (E_0) . \eea
On the other hand,
$$ F\cdot d{\cal S} = F\cdot  \left[ C^-+ 2F + \sum_{\jmath}
(1-m_{\jmath} ) E_{\jmath}\right]=3$$
 and
$$E_0\cdot d{\cal S} = E_0\cdot \left[ C^-+ 2F + \sum_{\jmath}
(1-m_{\jmath} ) E_{\jmath}\right]= 1+ (1-m_1) = 2+E_0^2 , $$
so that
\bea
\left| \begin{array}{cc}
F\cdot \Upsilon & F\cdot d{\cal S}\\
E_0\cdot \Upsilon &E_0\cdot d{\cal S}
\end{array}\right|  &=& \left| \begin{array}{cc}
 \Omega (F)  & 3\\
-(2+E_0^2)\Omega (F)+ 2(1+E_0^2) \Omega (E_0) &2+E_0^2
\end{array}\right| \\ &&\\&=&
4(2+E_0^2) \Omega (F) - 6( 1+E_0^2) \Omega (E_0).
\eea
Now this last expression is certainly non-zero if
$E^2_0$ is $-1$ or $-2$; and if
$E_0^2 \leq -3$,  the inequality $-(E^2_0+1)\Omega (E_0)<
-(E^2_0+1) \Omega (F)/3$
yields
$$\left| \begin{array}{cc}
F\cdot \Upsilon & F\cdot d{\cal S}\\
E_0\cdot \Upsilon &E_0\cdot d{\cal S}
\end{array}\right| <
2(3+E_0^2) \Omega (F) \leq 0. $$
 Hence $\Upsilon \not\equiv
0\bmod d{\cal S}$ in
$H_2(N)$,  and the claim holds
whenever $N\neq \hat{N}$.

We now come to  the case in which $N=\hat{N}$.
Then, since $F\cdot F=F\cdot E_{\jmath}=0$,
  $$F\cdot d{\cal S}=F\cdot
\left[ C^++C^- + 2(1-{\bf g}) F+ \sum_{\jmath}
(1-m_{\jmath})E_{\jmath}\right]=
2,$$ and
$$ F\cdot \Upsilon = F\cdot
\left[ \Omega (F)
[C^+-C^-]+\Omega (C^+-C^-)F+(E_{\jmath} \mbox{ terms})
\right] = 0.$$
It therefore suffices to show that
$\Upsilon\neq 0$
in $H_2(N)$.
But
$$C^+\cdot \Upsilon = \Omega (C^+-C^-)+ (C^+)^2\Omega (F)
+\sum_{\imath} \Omega (E_{\imath})$$
and
$$C^-= \Omega (C^+-C^-) -(C^-)^2\Omega (F) +
\sum_{\jmath} (1-m_{\jmath}) \Omega (E_{\jmath}) , $$
so that
$$\frac{1}{2} (C^+-C^-) \cdot\Upsilon = \frac{(C^+)^2+ (C^-)^2}{2}
\Omega (F) + \sum_{\jmath} (m_{\jmath}-1) \Omega (E_{\jmath}) . $$
Now $-[(C^+)^2+ (C^-)^2]$ is precisely the number of times
one must blow up along fixed curves of the ${\Bbb C}^{\times}$-action
in order to obtain $N$ from a fiber-minimal model; thus
$-[(C^+)^2+ (C^-)^2]\geq n$, where $n$ is  the
number of singular fibers of $N\to \Sigma$.
If, on the other hand, the first and last exceptional
curves of every fiber were to have area $ \geq \Omega (F)/4$,
we would have $$\sum_{\jmath} (m_{\jmath} - 1) < n \left[ \Omega (F) - 2
\frac{\Omega (F)}{4}\right] = \frac{n}{2} \Omega (F), $$ and it would
therefore follow that
$\frac{1}{2} (C^+-C^-) \cdot\Upsilon < \frac{n}{2} \Omega (F) -
\frac{n}{2} \Omega (F)=0$, implying
$\Upsilon \neq 0$, as desired.

We may therefore  assume that either the first or the last
curve of some singular fiber has area $<  \Omega (F)/4$;
and, by reversing the sign of $\Xi$ if necessary,
 we may assume that the curve in question is actually
the {\em first} in its fiber.
We thus have an exceptional curve $E_0$ which meets $C^-$
and which satisfies $\Omega (E_0) < \Omega (F)/4$.
Now, by essentially the same calculation
we used in the  $N\neq \hat{N}$ case,
$$E_0\cdot \Upsilon = -(2+ E_0^2)\Omega (F)
+ 2(1+E_0^2) \Omega (E_0) ,$$
and in particular $\Upsilon \neq 0$
if $E^2_0$ is $-1$ or $-2$. If, on the other hand,
$E_0^2\leq -3$, the inequality
$(1+E^2_0) \Omega (E_0) > (1+E^2_0) \Omega (F)/4$ tells us  that
$$E_0\cdot \Upsilon >   -\frac{1}{2}(E_0^2 + 3)
 \Omega (F) \geq 0.$$
Thus $\Upsilon \not\equiv 0 \bmod d{\cal S}$,  and
$H^2(Z_N, \Theta)=H^2(Z_N, \Theta_{Z,D\bar{D}})=
H^2(Z, \Theta \otimes \kappa^{-1})=0$, as claimed.
 \end{proof}

\end{document}